\documentclass[12pt]{amsart}
\usepackage{amssymb,amsmath, amsthm,latexsym}

\theoremstyle{plain}

\newtheorem{lem}{Lemma}[section]
\newtheorem{theo}[lem]{Theorem}
\newtheorem{prop}[lem]{Proposition}
\newtheorem{cor}[lem]{Corollary}

\parskip=\medskipamount
\font\k=cmr7

  \newcommand {\di}{\mbox{\k disc}}
  \newcommand {\reg}{\mbox{\k reg}}
  \newcommand {\tm}{\mbox{\k temp}}
  \newcommand {\C}{{\mathbb C}}
  \newcommand {\N}{{\mathbb N}}
  \newcommand {\R}{{\mathbb R}}
  \newcommand {\Z}{{\mathbb Z}}
  \newcommand {\Q}{{\mathbb Q}}
  \newcommand {\A}{{\mathbb A}}
  \newcommand {\af}{{\mathfrak a}}
  \newcommand {\gf}{{\mathfrak g}}
  
  \newcommand {\nf}{{\mathfrak n}}
  \newcommand {\hf}{{\mathfrak h}}
  \newcommand {\ho}{{\mathfrak o}}
 \newcommand {\mX}{{\mathfrak X}}
 \newcommand {\mM}{{\mathfrak M}}
 \newcommand {\mN}{{\mathfrak N}}
 \newcommand {\mO}{{\mathfrak O}}
  
\renewcommand {\H}{{\mathcal H}}
  
  \newcommand {\Co}{{\mathcal C}}
  
\renewcommand {\L}{{\mathcal L}}
  
 \newcommand {\cP}{{\mathcal P}}
 \newcommand {\cF}{{\mathcal F}}
 \newcommand {\cL}{{\mathcal L}}
 \newcommand {\cA}{{\mathcal A}}
\newcommand {\ba}{\backslash}
 \newcommand {\ov}{\overline}
 \newcommand {\bt}{{\bf t}}
 \newcommand {\bs}{{\bf s}}
 
\renewcommand{\Im}{\operatorname{Im}}
\renewcommand{\Re}{\operatorname{Re}}

\newcommand{\tr}{\operatorname{tr}}
\newcommand{\Id}{\operatorname{Id}}
\newcommand{\Hom}{\operatorname{Hom}}
\newcommand{\vol}{\operatorname{vol}}

\newcommand{\supp}{\operatorname{supp}}

\newcommand{\SL}{\operatorname{SL}}
\newcommand{\GL}{\operatorname{GL}}
\newcommand{\SO}{\operatorname{SO}}
\newcommand{\rO}{\operatorname{O}}

\begin{document}
%Topmatter
\title[Trace formula]{\Large\bf On the spectral side of the Arthur
  trace formula}
\date{\today}

\author{Werner M\"uller}
\address{Universit\"at Bonn\\
Mathematisches Institut\\
Beringstrasse 1\\
D -- 53115 Bonn, Germany}
\email{mueller@math.uni-bonn.de}


%End topmatter
\maketitle

\setcounter{section}{-1}
\section{Introduction}

Let $G$ be a connected reductive algebraic group defined over $\Q$ and 
let $G(\A)$ be
the group of points  of $G$ with values in the ring of ad\`eles of $\Q$.
Then $G(\Q)$ embeds diagonally as a discrete subgroup of $G(\A)$. 
Let $G(\A)^1$ be the intersection of the kernels of the
maps $x\mapsto|\chi(x)|$, $x\in G(\A)$, where $\chi$ ranges over the group
$X(G)_\Q$ of characters of $G$ defined over $\Q$. Then the (noninvariant)
trace formula of Arthur is an identity
$$\sum_{\ho\in\mO}J_\ho(f)=\sum_{\chi\in\mX}J_\chi(f),\quad f\in C_0^\infty(G(\A)^1),$$
between distributions on $G(\A)^1$. The left hand side is  the
{\it geometric side} and the right hand side the {\it spectral side} of the 
trace formula. The distributions $J_{\ho}$ are parametrized by semisimple 
conjugacy in $G(\Q)$ and are closely related to weighted orbital integrals on
$G(\A)^1$. 

In this paper we are concerned with the spectral side of the
trace formula. The distribution $J_\chi$
are defined in terms of truncated Eisenstein series. They are parametrized 
by the set  of cuspidal data $\mX$ which consists of the Weyl group orbits of
pairs $(M_B,r_B)$, where $M_B$ is the Levi component of a standard parabolic
subgroup and $r_B$ is an irreducible cuspidal automorphic representation 
of $M_B(\A)^1$. 
In \cite{A4}, Arthur has derived an explicit formula for
 the distributions $J_\chi$ which expresses them in terms of generalized 
logarithmic derivatives of intertwining operators. So far, the resulting
integral-series is only known to converge conditionally. This suffices, 
for example,  for the comparison of trace formulas
which, at present, is the main application of the trace formula. 
However, with regard to potential applications of the trace formula in
spectral theory and geometry it would be highly desirable to know that 
 the spectral side of the trace formula is absolutely convergent. This
would also simplify the applications of the trace formula in the theory of 
automorphic forms \cite{Lb}. 

The problem of the absolute convergence of the spectral side of the trace
formula is the main issue of the present paper. We will not settle the 
problem, but we shall reduce it to a question
about  local components of automorphic representations. 

To describe the results in more detail we have to introduce some  notation. 
We fix a Levi component $M_0$ of a minimal parabolic subgroup $P_0$ of $G$.
Let $P$ be a parabolic subgroup 
of $G$, defined over $\Q$, with unipotent radical $N_P$. Let $M_P$ be the 
unique Levi component of $P$ which contains $M_0$. We denote the split 
component of the center of $M_P$ by $A_P$ and its Lie algebra by $\af_P$. 
For parabolic groups $P\subset Q$ there is a natural surjective map $\af_P\to
\af_Q$ whose kernel we will denote by $\af_P^Q$.
Let $\cA^2(P)$ be the space of square integrable automorphic forms on
$N_P(\A)M_P(\Q)\backslash G(\A)$.
Let $Q$ be another parabolic subgroup of $G$, defined over $\Q$, with Levi
component $M_Q$, split component $A_Q$ and corresponding Lie algebra $\af_Q$.
Let $W(\af_P,\af_Q)$ be the set of all linear isomorphisms from $\af_P$ to 
$\af_Q$
which are restrictions of elements of the Weyl group $W(A_0)$.  The theory of 
Eisenstein series associates to each $s\in W(\af_P,\af_Q)$ an intertwining 
operator
$$M_{Q|P}(s,\lambda): \cA^2(P)\rightarrow \cA^2(Q),\quad \lambda\in\af_{P,\C}^*,$$
which for $\Re(\lambda)$ in a certain chamber, can be defined by an 
absolutely convergent integral and admits an analytic continuation to a 
meromorphic function of $\lambda\in\af^*_{P,\C}$. Set
$$M_{Q|P}(\lambda):=M_{Q|P}(1,\lambda).$$
Let $\Pi(M_P(\A)^1)$ be the
set of equivalence classes of irreducible unitary representations of 
$M_P(\A)^1$. Let $\chi\in\mX$ and 
$\pi\in\Pi(M_P(\A)^1)$. Then $(\chi,\pi)$ singles out a certain subspace
$\cA^2_{\chi,\pi}(P)$ of $\cA^2(P)$ (see \S 1.6).
 Let $\ov \cA^2_{\chi,\pi}(P)$ be the Hilbert 
space completion of $\cA^2_{\chi,\pi}(P)$ with respect to the canonical inner
product. For each $\lambda\in\af_{P,\C}^*$ we have an induced representation
 $\rho_{\chi,\pi}(P,\lambda)$  of $G(\A)$ in $\ov \cA^2_{\chi,\pi}(P)$.

For each
Levi subgroup $L$ let $\cP(L)$ be the set of all parabolic subgroups with
Levi component $L$. If $P$ is a parabolic subgroup, let $\Delta_P$ denote the
set of simple roots of $(P,A_P)$. Let $L$ be a Levi subgroup which contains 
$M_P$. Set
\begin{equation*}
\begin{split}
\mM_L&(P,\lambda)=\\
&\lim_{\Lambda\to0}\left(\sum_{Q_1\in\cP(L)}
\vol(\af_{Q_1}^G/\Z(\Delta^\vee_{Q_1}))M_{Q|P}(\lambda)^{-1}
\frac{M_{Q|P}(\lambda+\Lambda)}{\prod_{\alpha\in\Delta_{Q_1}}\Lambda(\alpha^\vee)}\right),
\end{split}
\end{equation*}
where $\lambda$ and $\Lambda$ are constrained to lie in $i\af_L^*$, and for 
each $Q_1\in\cP(L)$, $Q$ is a group in $\cP(M_P)$ which is contained in $Q_1$.
Then $\mM_L(P,\lambda)$ is an unbounded operator which acts
on the Hilbert space $\ov \cA^2_{\chi,\pi}(P)$.  
In the special case that $L=M$ 
and $\dim\af_L^G=1$, the operator $\mM_L(P,\lambda)$ has a simple description. 
 Let $P$ be a parabolic subgroup with Levi component $M$. 
 Let $\alpha$ be the unique simple root of  $(P,A_P)$ and let $\tilde\omega$
be the element in $(\af_M^G)^*$ such that $\tilde\omega(\alpha^\vee)=1$. Let
$\ov P$ be the opposite parabolic group of $P$.  Then
$$\mM_L(P,z\tilde\omega)=-\vol(\af_M^G/\Z\alpha^\vee)
M_{\ov P|P}(z\tilde\omega)^{-1}\cdot
\frac{d}{dz}M_{\ov P|P}(z\tilde\omega).$$

Let $f\in C_c^\infty(G(\A)^1)$. Then Arthur
\cite[Theorem 8.2]{A4} proved that $J_\chi(f)$ equals the sum over Levi 
subgroups $M$ containing $M_0$, over $L$ containing $M$, over $\pi\in
\Pi(M(\A)^1)$, and over $s\in W^L(\af_M)_{\reg}$, a certain subset of the 
Weyl group, of the product of
\begin{equation*}
|W_0^M||W_0|^{-1}|\det(s-1)_{\af^L_M}|^{-1}|\cP(M)|^{-1},
\end{equation*}
a  factor to which we need not pay too much attention , and of
\begin{equation}\label{0.1}
\int_{i\af_L^*/i\af_G^*}\sum_{P\in\cP(M)}\tr(\mM_L(P,\lambda)M_{P|P}(s,0)
\rho_{\chi,\pi}(P,\lambda,f))\;d\lambda.
\end{equation}
So far, it is only known that $\sum_{\chi\in\mX}|J_\chi(f)|<\infty$ and the
 goal is to show that the integral--sum obtained by summing (\ref{0.1}) over 
$\chi\in\mX$ and $\pi\in\Pi(M(\A)^1)$ is absolutely convergent with respect 
to the trace norm. 

\smallskip
Given $\pi\in\Pi(M(\A))$ with $\pi=\otimes_v\pi_v$, let 
$J_{Q|P}(\pi_v,\lambda)$ be the local intertwining operator between the 
induced representations $I^G_P(\pi_{v,\lambda})$ and $I^G_Q(\pi_{v,\lambda})$.
By \cite[\S15]{CLL} and \cite{A7}
there exist normalizing factors $r_{Q|P}(\pi_v,\lambda)$ such 
that the normalized intertwining operators 
\begin{equation*}
R_{Q|P}(\pi_v,\lambda)=r_{Q|P}(\pi_v,\lambda)^{-1}
J_{Q|P}(\pi_v,\lambda)
\end{equation*}

satisfy the conditions of Theorem 2.1 of \cite{A7}. 

If $v<\infty$, let 
$K_v\subset G(\Q_v)$ be an open compact subgroup. Denote by 
$R_{Q|P}(\pi_v,\lambda)_{K_v}$ the restriction of $R_{Q|P}(\pi_v,\lambda)$ to
the subspace $\H_P(\pi_v)^{K_v}$ of $K_v$-invariant vectors in the 
Hilbert space $\H_P(\pi_v)$ of the induced representation.
If $v=\infty$, let $K_\infty\subset G(\R)$ be a maximal compact subgroup. Given
$\pi\in\Pi(M(\R))$ and $\sigma\in\Pi(K_\infty)$, let 
$R_{Q|P}(\pi,\lambda)_\sigma$ be the
restriction of $R_{Q|P}(\pi,\lambda)$ to the $\sigma$-isotypical subspace
of $\H_P(\pi)$. Let $\lambda_{\pi}$ and $\lambda_\sigma$ denote the 
Casimir eigenvalues of $\pi$ and $\sigma$, respectively.
 
For a given place $v$, let $\Pi_{\di}(M(\Q_v))$ be the subset 
consisting of all representations $\pi_v\in\Pi(M(\Q_v))$ such that there
exists an automorphic representation $\pi$ in the discrete spectrum of $M(\A)$
whose local component at the place $v$ is $\pi_v$. Finally, let 
$\Co^1(G(\A)^1)$ be the space of integrable rapidly decreasing functions on
$G(\A)^1$ (see \S1.3 for its definition). 
Then our main result is the following theorem, which reduces the problem of 
the absolute convergence of the spectral side of the Arthur trace formula to a
problem about local components of automorphic representations.

\begin{theo}\label{th0.1}
Suppose that for every $M\in\cL(M_0)$, $Q,P\in\cP(M)$ and every place $v$
the following holds.

\smallskip
\noindent
1) If $v<\infty$, then for every open compact subgroup $K_v\subset G(\Q_v)$
and every invariant differential operator $D_\lambda$ on $i\af_M^*$ there
exists $C>0$ such that
\begin{equation}\label{0.2}
\parallel D_\lambda R_{Q|P}(\pi_v,\lambda)_{K_v}\parallel\le C
\end{equation}
for all $\lambda\in i\af_M^*$ and all $\pi_v\in\Pi_{\di}(M(\Q_v))$.

\smallskip
\noindent
2) If $v=\infty$, then for all invariant differential operators 
$D_\lambda$
on $i\af_M^*$ there exist $C>0$ and $N\in\N$ such that
\begin{equation}\label{0.3}
\parallel D_\lambda R_{Q|P}(\pi,\lambda)_\sigma\parallel\le
C(1+\parallel\lambda\parallel+|\lambda_{\pi}|+|\lambda_\sigma|)^N
\end{equation}
for all $\lambda\in i\af_M^*$, $\sigma\in\Pi(K_\infty)$ and 
$\pi\in\Pi_{\di}(M(\R))$.

\smallskip
\noindent
Then for every $f\in\Co^1(G(\A)^1)$, the spectral side of the trace formula 
is absolutely convergent.
\end{theo}

We add some comments about the assumptions 1) and 2). 
It follows from results of Arthur \cite[p.51]{A5} and 
\cite[Lemma 2.1]{A8} that (\ref{0.2}) and (\ref{0.3}) hold uniformly for
tempered representations $\pi_v$. On the other hand, to establish (\ref{0.2}),
(\ref{0.3}) or (\ref{0.4}) is not a  problem of pure local harmonic 
analysis. One can not expect that these estimations  
will  hold uniformly for all $\pi_v\in\Pi(M(\Q_v))$. Let, for
example, $\dim\af_M/\af_G=1$ and suppose that for each $\epsilon>0$ there
exists $\pi_v\in\Pi(M(\Q_v))$ such that the normalized intertwining operator
$R_{\ov P|P}(\pi_v,\lambda)$ has a pole $\lambda_0$ with $|\Re(\lambda_0)|\le
\epsilon$. Then it is certainly not possible to obtain uniform estimates for
derivatives of $R_{\ov P|P}(\pi_v,\lambda)$ along the imaginary axis. 
An example where this actually occurs is $\GL_n$.

Especially the uniformity in $\sigma$ in (\ref{0.3}) seems to be difficult to 
achieve. Of course, this condition can be relaxed in various ways. 
If we relax (\ref{0.3}) by not requesting 
uniformity in $\sigma$, we get the following weaker version of 
Theorem \ref{0.1} which suffices for many purposes. Let $K=\prod_v K_v$ be a 
maximal compact subgroup of $G(\A)$ which is admissible relative to $M_0$
(see \S1.2).

\begin{theo}\label{th0.2}
Suppose that in Theorem \ref{0.1} in place of condition 2) the following 
condition holds:

\smallskip
\noindent
2') If $v=\infty$, then for all invariant differential operators 
$D_\lambda$
on $i\af_M^*$ and all $\sigma\in\Pi(K_\infty)$ there exist $C>0$ and $N\in\N$ 
such that
\begin{equation}\label{0.4}
\parallel D_\lambda R_{Q|P}(\pi,\lambda)_\sigma\parallel\le
C(1+\parallel\lambda\parallel+|\lambda_{\pi}|)^N
\end{equation}
for all $\lambda\in i\af_M^*$ and $\pi\in\Pi_{\di}(M(\R))$.

\smallskip
\noindent
Then for every $K$-finite $f\in\Co^1(G(\A)^1)$, the spectral side of the 
trace formula is absolutely convergent.
\end{theo}

At the moment we don't know how to prove any of the conditions (\ref{0.2}), 
(\ref{0.3}) and (\ref{0.4}) in general. However, for $G=\GL_n$, considered as 
an algebraic group over a number field, we are able to prove (\ref{0.2}) 
and (\ref{0.4}). The method relies on work of Luo, Rudnick and Sarnak 
\cite{LRS}
who established nontrivial bounds towards the generalized Ramanujan conjecture.
For $\GL_n$ any local component of a cuspidal automorphic representation is
equivalent to a full induced representation $I^G_P(\tau,\bs)$ where $\tau$ is 
tempered and the parameters $\bs=(s_1,...,s_r)$ satisfy $s_1>s_2>\cdots>s_r$ and $|s_i|<1/2$. If
$\pi_v$ is unramified, it follows from Theorem 2 of \cite{LRS} that the 
$s_i$'s satisfy the nontrivial bound
\begin{equation}\label{0.4a}
|s_i|<\frac{1}{2}-\frac{1}{n^2+1},\quad i=1,...,r.
\end{equation}
Using the method of \cite{LRS}, one can show that (\ref{0.4a}) holds also for
the ramified components. Furthermore, 
using the work of M{\oe}glin and Waldspurger \cite{MW} on the residual 
spectrum, one can
show that similar nontrivial bounds exist for the continuous parameters
of any local component of an automorphic representation in the discrete
spectrum of $\GL_m(\A)$. These bounds are the essential ingredients in the
proof of (\ref{0.2}) and (\ref{0.4}) in the case of $\GL_n$. 
Details will appear in a forthcoming paper with B. Speh \cite{MS}.

Now we shall explain the main steps of the proof of Theorem 0.1. First observe that
$M_{P|P}(s,0)$ is unitary. Therefore, in order to estimate the trace norm of
(\ref{0.1}),  it suffices to estimate the integral 
\begin{equation}\label{0.5}
\int_{i\af_L^*/i\af_G^*}\parallel\mM_L(P,\lambda)
\rho_{\chi,\pi}(P,\lambda,f)\parallel_1\;d\lambda.
\end{equation}

To deal with this integral, we introduce a certain normalization of 
intertwining operators. For $\pi\in \Pi(M(\A))$ let 
$\cA^2_\pi(P)$ be the space of square integrable automorphic forms of type 
$\pi$ (see \S 1). Let $M_{Q|P}(\pi,\lambda)$ denote the restriction of the 
intertwining 
operator $M_{Q|P}(\lambda)$ to $\cA^2_\pi(P)$. Let $\pi=\otimes_v\pi_v$ 
and let  $r_{Q|P}(\pi_v,\lambda)$ be the normalizing factor for the local
intertwining operator considered above.  Suppose
$\pi=\otimes  \pi_v$ occurs in the discrete spectrum of $M(\A)$, which is
equivalent to $\cA^2_\pi(P)\neq 0$, then the Euler product
\begin{equation*}
r_{Q|P}(\pi,\lambda)=\Pi_vr_{Q|P}(\pi_v,\lambda)
\end{equation*}

converges absolutely in some chamber and $r_{Q|P}(\pi,\lambda)$ admits a 
meromorphic continuation to $\af^*_{M,\C}$. Using this meromorphic 
function, we introduce the normalized global intertwining operator by
\begin{equation}\label{0.6}
N_{Q|P}(\pi,\lambda)=r_{Q|P}(\pi,\lambda)^{-1}M_{Q|P}(\pi,\lambda).
\end{equation}

By definition, the operator $N_{Q|P}(\pi,\lambda)$ is equivalent to the 
direct sum of finitely many copies of
$\otimes_v R_{Q|P}(\pi_v, \lambda)$.

Let $\mM_L(P,\pi,\lambda)$ be the restriction of $\mM_L(P,\lambda)$ to the
subspace $\cA^2_\pi(P)$. 
It follows from Arthur's theory of $(G,M)$ families \cite[p.1329]{A4} that
\begin{equation*}
\mM_L(P,\pi,\lambda)=\sum_S{\mN}^\prime_S(P,\pi,\lambda)\nu^S_L(P,\pi,\lambda),
\end{equation*}

where the sum runs over all parabolic subgroups $S$ containing $L$, the
operator 
${\mN}^{\prime}_S(P,\pi,\lambda)$ is built out of normalized intertwining
operators on the local groups $G(\Q_v)$ and $\nu^S_L(P,\pi,\lambda)$ is 
a scalar valued function which is defined in terms of normalizing factors.
This reduces the estimation of the integral (\ref{0.5}) to two separate
problems, one involving ${\mN}^\prime_S(P,\pi,\lambda)$ and the other one
$\nu^S_L(P,\pi,\lambda)$.

First we are dealing with $\nu_L^S(P,\pi,\lambda)$. By Proposition 7.5 of 
\cite{A4}, $\nu^S_L(P,\pi,\lambda)$ can be expressed in terms of logarithmic 
derivatives of normalizing factors associated with maximal
parabolic subgroups in certain Levi subgroups. Therefore we may assume that 
$\dim( \af_P/\af_G)=1$. Let $\alpha$ be the unique simple root of $(P,A)$. 
Then there exists a meromorphic function $\widetilde r_{{\overline P}|P}
(\pi,z)$ of one variable such that 
$r_{{\overline P|P}}(\pi,\lambda)=\widetilde r_{{\overline P}|P}(\pi,\lambda
(\alpha^\vee ))$, and our problem is to derive estimates, which are uniform 
with respect to $\pi$, of integrals of the form
\begin{equation}\label{0.7}
\int_{\R} \left|\widetilde r_{\overline P|P}(\pi,iu)^{-1}\frac{d}{du}
\widetilde r_{{\overline P}|P}(\pi,iu)\right|(1+u^2)^{-N}du.
\end{equation}

To deal with this integral, we note that 
$\widetilde r_{{\overline P}|P}(\pi,z)$ is a meromorphic function of order 
 $n=16\dim G+2$. This  follows from (\ref{0.6}), since  
by Theorem 0.1 of \cite{Mu3}, the matrix coefficients of 
$M_{Q|P}(\pi,\lambda)$ are meromorphic functions of order $\le n$ and
by Theorem 2.1 of \cite{A7}, the normalized local intertwining operators 
$R_{Q|P}(\pi_v,\lambda)$ are rational functions of 
$q_v^{-\lambda(\alpha^\vee)}$, if $v<\infty$, and of 
$\lambda(\alpha^\vee)$, if $v=\infty$. 
Thus there exist entire functions $r_i(\pi,z)$, $i=1,2$, of order $ \le n$ 
such that $\widetilde r_{{\overline P}|P}(\pi,z)=r_1(\pi,z)/r_2(\pi,z)$.
Using the representation of $r_i(\pi,z)$ as a Weierstra\ss{} product, we reduce
the estimation of the integral (\ref{0.7}) to the
 the estimation of the number of poles, counted with their
order, of $\widetilde r_{{\overline P}|P}(\pi,z)$ in a circle of radius
$R>0$. By (\ref{0.6}), this problem is closely related to the estimation
of the number of poles, counted with their order, of matrix coefficients of
$M_{Q|P}(\pi,\lambda)$ in a circle of radius $R>0$. The latter 
problem has been settled in \cite[Proposition 6.6]{Mu3}. 
Together with Proposition 7.5 of \cite{A4}, these estimates imply estimates 
for the corresponding integrals involving $\nu^S_L(P,\pi,\lambda)$. In this
way we get Theorem \ref{th5.4} which is our main technical result.

Next consider $\mN^\prime_S(P,\pi,\lambda)$. Given an open compact subgroup
$K_f$ of $G(\A_f)$ and  $\sigma\in\Pi(K_\infty)$, let 
$\cA^2_\pi(P)_{K_f,\sigma}$ be the subspace of $\cA^2_\pi(P)$ consisting
of all automorphic forms which are $K_f$-invariant and transform under
$K_\infty$ according to $\sigma$. Let 
$\mN^\prime_S(P,\pi,\lambda)_{K_f,\sigma}$ be the restriction of 
$\mN^\prime_S(P,\pi,\lambda)$ to the subspace $\cA^2_\pi(P)_{K_f,\sigma}$.
Now observe that for any $h\in\Co^1(G(\A)^1)$ there exists an open compact
 subgroup $K_f$ of $G(\A_f)$ such that $h$ is left and right invariant under 
$K_f$. Then the estimation of $\parallel\mN^\prime_S(P,\pi,\lambda)
\rho_{\chi,\pi}(P,\lambda,h)\parallel_1$ can be reduced to the estimation of
$\parallel \mN^\prime_S(P,\pi,\lambda)_{K_f,\sigma}\parallel$ where 
$\sigma$ runs over $\Pi(K_\infty)$. 
 By Arthur's theory of $(G,M)$-families, the estimation of the norm of the
finite rank operators $\mN_S^\prime(P,\pi,\lambda)_{K_f,\sigma}$  
can be reduced to the estimation of derivatives of finitely many normalized
local intertwining operators $R_{Q|P}(\pi_v,\lambda)_{K_v}$, $v<\infty$, and
$R_{Q|P}(\pi_\infty,\lambda)_{\sigma}$. Combined with Theorem
\ref{th5.4} this implies Theorem \ref{0.1}. The proof of Theorem \ref{0.2} is
similar. 

The paper is organized as follows. In \S 1 we collect some 
preliminary facts. In \S 2 we discuss briefly normalized local and
global intertwining operators. The local normalizing factors are studied in
some detail in \S 3. We recall the definition of the normalizing factors
and we prove some results that we need in the next section. In \S 4 we
investigate the poles of the global normalizing factors. This section is mainly
 based on results obtained in \cite{Mu3}. In \S 5 we establish Theorem 
\ref{th5.4} which is the main result about generalized logarithmic derivatives 
of global normalizing factors. In \S 6 we study the absolute
convergence of the spectral side of the trace formula and we prove our main
results, Theorem \ref{0.1} and  Theorem \ref{th0.2}. In \S 7 we discuss the
example of $\GL_n$ and we sketch a method to prove (\ref{0.2}) and (\ref{0.4}).

\smallskip
{\bf Acknowledgments} Most of this work was done during the author's stay
at the Institute for Advanced Study in Princeton. The author is very grateful 
to the IAS for hospitality and financial support by a grant from the 
Ellentuck Fund. He also would like to thank J. Cogdell,  
D. Ramakrishnan and F. Shahidi for some very useful discussions. 
 Especially he is very grateful to
W. Hoffmann and the referee for their comments, suggestions and corrections
which helped to improve the paper considerably.

%\newpage

\section[Preliminaries]{Preliminaries}
\setcounter{equation}{0}

We shall follow partially the notation introduced by Arthur 
\cite {A1}-\cite{A4}.

\subsection{}

Let $G$ be a reductive  algebraic group defined  over 
${\Q}$. As in \cite{A4}, \cite{A5}, we shall fix a subgroup $M_0$ of $G$,
defined over ${\Q}$, which is a Levi component of some minimal parabolic
subgroup of $G$ defined over $\Q$.  In this paper, a {\it parabolic
subgroup} will mean a parabolic subgroup of $G$, defined over ${\Q}$,
and a {\it Levi subgroup} of $G$ will mean a subgroup of $G$ which contains
$M_0$ and is the Levi component of some parabolic subgroup of $G$. It is a
reductive subgroup of $G$ which is defined over $\Q$. If $M\subset L$ are Levi
subgroups, we denote the set of Levi subgroups of $L$ which contain $M$ by
${\cL}^L(M)$. Furthermore, let ${\cF}^L(M)$ denote the set of
parabolic subgroups of $L$ defined over $\Q$ which contain $M$, and let
${\cP}^L(M)$ be the set of groups in ${\cF}^L(M)$ for which $M$
is a Levi component. If $L=G$, we shall denote these sets by ${\cL}(M)$, ${\cF}(M)$ and ${\cP}(M)$, respectively. Suppose that 
$P\in{\cF}^L(M)$. Then 
$$P=N_PM_P,$$
where $N_P$ is the unipotent radical of $P$ and $M_P$ is the unique Levi
component of $P$ which contains $M$.

Suppose that $M\subset M_1\subset L$ are
Levi subgroups of $G$.  If $Q\in\cP^L(M_1)$ and $R\in{\cP}^{M_1}(M)$, there is a
unique group $Q(R)\in{\cP}^L(M)$ which is contained in $Q$ and whose
intersection with $M_1$ is $R$. 

Let $A_P$ be the split component of the center of $M_P$. $A_P$ is defined over
$\Q$. Let
$X(M_P)_{\mathbb Q}$ be the group of characters of $M_P$ defined over 
$\Q$. Then 
\begin{equation*}
{\mathfrak a}_P=\mbox{Hom}(X(M_P)_{\mathbb Q},{\mathbb R})
\end{equation*}
is a real vector space whose dimension equals that of $A_P$. Its dual space is
\begin{equation*}
{\mathfrak a}^*_P= X(M_P)_{\mathbb Q}\otimes {\mathbb R}.
\end{equation*}
We shall often denote $A_P$, $\af_P$ and $\af_P^*$ by $A_M$, $\af_M$ and
$\af_M^*$, respectively,  since they depend
only on $M$. 
Also, we shall write $A_0=A_{M_{0}}$, ${\mathfrak a}_0={\mathfrak a}_{M_{0}}$
and ${\mathfrak a}^*_0={\mathfrak a}^*_{M_{0}}$.

Let $P\in{\cF}(M_0)$. We shall denote the roots of $(P,A_P)$ by $\Phi_P$, the
reduced roots by $\Sigma^r_P$, and
the simple roots by $\Delta_P$. They are elements in $X(A_P)_\Q$ and are
canonically embedded in $\af_P^*$.

Let $P$ and $Q$ be groups in $\cF(M_0)$ with $P\subset Q$. Then there is a
canonical surjection $\af_P\rightarrow \af_Q$ and a canonical injection ${\mathfrak a}^*_Q\hookrightarrow {\mathfrak a}^*_P$.
The kernel of the first map will be denoted by
${\mathfrak a}^Q_P.$ Then
${\mathfrak a}^Q_P$
is dual to 
${\mathfrak a}^*_P /{\mathfrak a}^*_Q.$

The group $M_Q\cap P$ is a parabolic subgroup of $M_Q$ with unipotent radical 
$$N_P^Q=N_P\cap M_Q.$$
Let $\Delta_P^Q$ be the set of simple roots of $(M_Q\cap P,A_P)$. Then
$\Delta_P^Q$ is a subset of $\Delta_P$. We may identify $\af_Q$ with the
subspace 
$$\{H\in\af_P\mid \alpha(H)=0,\;\alpha\in\Delta_P^Q\}.$$
Furthermore, to $\Delta_P^Q$ one can associate a basis $\{\alpha^\vee\mid
\alpha\in\Delta_P^Q\}$ of $\af_P^Q$ and $\hat\Delta_P^Q$ is defined to be the
corresponding dual basis of $(\af_P^Q)^*$ \cite{A2}. Then $\Delta_P^Q$ and
$\hat\Delta_P^Q$ are naturally embedded subsets of $\af_0^*$. Let

\begin{equation*}
{\mathfrak a}^+_P=\{ H\in{\mathfrak a}_P \mid \alpha (H)>0\quad
\mbox{for all}\quad \alpha   \in\Delta_P\},
\end{equation*}

and

\begin{equation*}
({\mathfrak a}^*_P)^+ =\{ \Lambda\in{\mathfrak a}^*_P\mid \Lambda(\alpha^\vee)>0
\quad \mbox{for all}\quad \alpha\in\Delta_P\}.
\end{equation*}

We shall denote  the restricted Weyl group of $(G,A_0)$ by $W_0$. It acts on 
${\mathfrak a}_0$ and ${\mathfrak a}^*_0$ in the usual way. For every $s\in
W_0$ we shall fix a representative $w_s$ in the intersection of $G(\Q)$ with
the normalizer of $A_0$. $w_s$ is determined modulo $M_0(\Q)$. 
If $P_1$ and $P_2$ are parabolic subgroups, let 
$W({\mathfrak a}_{P_1},{\mathfrak a}_{P_2})$ denote the set of distinct
isomorphisms  from 
${\mathfrak a}_{P_1}$ onto ${\mathfrak a}_{P_2}$ obtained by restricting
elements of $W_0$ to $\af_{P_1}$. $P_1$ and $P_2$ are said to be {\it
  associated} if $W({\mathfrak a}_{P_1},{\mathfrak a}_{P_2})$ is not empty.

\subsection{}

We fix an embedding of $G$ into $\GL_n,$ defined over $\Q$. For a given place
$v$ of $\Q,$ let $G({\Q}_v)$ be the group of ${\Q}_v$-rational points of $G$.
Let $\A$ be the ring of ad\`eles of $\Q$ and let $G(\A)$ be the corresponding
ad\`ele-valued group. If $f$ stands for the set of finite places of $\Q$ and
${\A}_f$ is the corresponding ring of finite ad\` eles, then

\begin{equation*}
G(\A)=G(\R)\times G({\A}_f).
\end{equation*}

For any prime $p,$ let

\begin{equation*}
G(\Z_p)=\GL_n(\Z_p)\cap G(\Q_p).
\end{equation*}

This is an open compact subgroup of $G(\Q_p)$. We shall fix a maximal compact
subgroup

\begin{equation*}
K=\prod_v K_v
\end{equation*}

of $G(\A)$ which is admissible relative to $M_0$ in the sense of
\cite{A5}. For any such $K$  the following properties hold:

\begin{enumerate}
\item[1)] For almost all primes $p$, one has
$K_p=G(\Z_p)$.

\item[2)] For every finite $p$,   $K_p$ is a special maximal compact
  subgroup. This implies that $G(\Q_p)=P_0(\Q_p)\cdot K_p$ for all
  $P_0\in\cF(M_0)$. 

\item[3)] The Lie algebras of $K_{\R}$ and $A_0(\R)$ are orthogonal with
  respect to the Killing form.
\end{enumerate}

Given $P\in\cF(M_0)$, let
$$M_P(\A)^1=\bigcap_{\chi\in X(M_P)_\Q}\ker(|\chi|).$$
This is a closed subgroup of $M_P(\A)$, and $M_P(\A)$ is the direct product of
$M_P(\A)^1$ and $A_P(\R)^0$. By the 
assumptions on $K$, $G(\A)=P(\A)K$. Therefore, any $x\in G(\A)$ can be written
as
\begin{equation}\label{1.1}
x=namk,\;    n\in N_P(\A),\;m\in M(\A)^1,\;
  a\in A_P(\R)^0,\; k\in K.
\end{equation}
Let 
\begin{equation*}
H_P\colon G(\A) \to {\mathfrak a}_P
\end{equation*}
be the associated height function as defined in \cite {A2}. Then $H_P(x)$ is
the image of $a\in A_P(\R)^0$ in the decomposition (\ref{1.1}) with respect to
the isomorphism $A_P(\R)^0\cong{\mathfrak a}_P$.

We shall fix a Euclidean norm $\parallel \cdot\parallel$
on ${\mathfrak a}_0$
which is invariant under the action of the Weyl group of 
$(G,A_0).$ On each space
${\mathfrak a}^Q_P \, ,\, P\subset Q,$
we take as Haar measure the Euclidean measure associated
to the restriction of 
$\parallel \cdot\parallel$ to
${\mathfrak a}^Q_P.$ 
We then normalize the Haar measures on 
$K$, $G(\A)$, $N_P(\A)$, $M_P(\A)$, $A_P(\R)^0$,  
$M_P(\A)^1$, etc. as in \cite {A2}.

\subsection{}

Let $\Xi$ and $\sigma$ be the functions that enter the definition of 
Harish-Chandra's Schwartz space on $G(\R)$ \cite[p.156]{Wa2} and extend them
to functions on $G(\A)$ in the obvious way. For any place $v$, let $G(\Q_v)^1$
denote the intersection of $G(\Q_v)$ with $G(\A)^1$. Let 
${\mathcal U}(\gf(\R)^1\otimes\C)$ be the universal enveloping algebra of the
complexification of the Lie algebra of $G(\R)^1$.  Let $K_f$ be an open 
compact subgroup of $G(\A_f)^1$. Then the double coset space 
$K_f\backslash G(\A)^1/K_f$ is a discrete union of countably many copies of 
$G(\R)^1$. In particular it is a differentiable manifold. Suppose that $\Omega$
is a  subset of $G(\A)^1$ such that $K_f\cdot\Omega\cdot K_f=\Omega$ and 
$K_f\backslash\Omega/K_f$ is the disjoint union of finitely many copies of 
$G(\R)^1$.
 Let $\Co^1(G(\A)^1;\Omega,K_f)$ be the
space of all functions $h\colon G(\A)^1\to\C$ satisfying the following 
conditions:
\begin{enumerate}
\item $h$ is bi-invariant under $K_f$, $\supp h\subset \Omega$, and 
$h: K_f\ba\Omega/K_f\to\C$ is a smooth function.
\item For all $D_1,D_2\in{\mathcal U}(\gf(\R)^1\otimes\C)$ and all $r\in\N$,
we have
\begin{equation*}
\begin{split}
\parallel h&\parallel_{D_1,D_2,r}\\
&:=\sup_{x\in G(\A)^1}\left((1+\sigma(x))^r\Xi^{-2}(x)|D_1\ast h\ast D_2(x)|
\right)<\infty.
\end{split}
\end{equation*}
\end{enumerate} 
Let $\Co^1(G(\A)^1;\Omega,K_f)$ be  equipped with the topology defined by 
the semi-norms $\parallel \cdot\parallel_{D_1,D_2,r}$.
Let $\Co^1(G(\A)^1)$ be the topological direct limit over all pairs 
$(\Omega,K_f)$ of the spaces $\Co^1(G(\A)^1;\Omega,K_f)$.

\subsection{}
Let $H$ be any algebraic group over $\Q$ and let $F$ be a local field. We 
shall denote by
$\Pi(H(\A))$ (resp. $\Pi(H(F))$, $\Pi(K)$, etc.) the 
 set of equivalence classes of irreducible 
unitary representations of $H(\A)$ (resp.  $H(F)$, $K$, etc.).

\subsection{}
Given a unitary character $\xi$ of $A_P(\R)^0$ , let 
$L^2(M_P(\Q)\ba M_P(\A))_\xi$ be the space of all measurable functions $\phi$ 
on $M_P(\Q)\ba M_P(\A)$ such that $\phi(xm)=\xi(x)\phi(m)$ for all
$x\in A_P(\R)^0$, $m\in M_P(\A)$, and $\phi$ is square integrable on 
$M_P(\Q)\ba M_P(\A)^1$. Let $\Pi_{\di}(M_P(\A))_\xi$ denote the subspace of
all $\pi\in\Pi(M_P(\A))$ which are equivalent to a subrepresentation of the
regular representation of $M_P(\A)$ on $L^2(M_P(\Q)\ba M_P(\A))_\xi$. Set
$$\Pi_{\di}(M_P(\A))=\bigcup_{\xi\in \Pi(A_P(\R)^0)}\Pi_{\di}(M_P(\A))_\xi.$$
Recall that  $\Pi(M_P(\A)^1)$ can be canonically identified
with the set of orbits  under the action of $i{\mathfrak a}^*_P$ defined by
\begin{equation*}
\pi\to \pi_\lambda =e^{\lambda(H_P(\cdot))} \pi,\quad
\pi\in\Pi(M_P(\A)),\quad \lambda \in i{\mathfrak a}^*_P.
\end{equation*}
Since $M_P(\A)$ is the direct product of $M_P(\A)^1$ and $A_P(\R)^0$, any
representation of $M_P(\A)^1$ corresponds to a representation of $M_P(\A)$
which is trivial on $A_P(\R)^0$. We identify these two representations and 
in this way we obtain an embedding of $\Pi(M_P(\A)^1)$ in $\Pi(M_P(\A))$.

Given $\pi\in\Pi(M_P(\A))$ with $\pi=\otimes_v\pi_v$, set 
$\pi_f=\otimes_{v<\infty}\pi_v.$
For an open compact subgroup $K_f\subset G(\A_f)$, let 
$$K_{M,f}=M_P(\A_f)\cap K_f.$$
Set
\begin{equation}\label{1.2}
\Pi_{\di}(M_P(\A);K_f)=\left\{\pi\in \Pi_{\di}(M_P(\A))\;\big|\;\pi_f^{K_{M,f}}
\not=\{0\}\right\}.
\end{equation}
Let $\Pi_{\di}(M_P(\A)^1;K_f)$ be the intersection of $\Pi_{\di}(M_P(\A);K_f)$
with the subspace $\Pi_{\di}(M_P(\A)^1)$ of $\Pi_{\di}(M_P(\A))$.
\subsection{}

Let $P=NM$ be a parabolic subgroup and let $\phi$ 
be a measurable, locally integrable function on $ N(\Q)\backslash G(\A)$. Then the constant term $\phi_P$ of $\phi$
along $P$ is defined for almost every $g$ by
\begin{equation}\label{1.4}
\phi_P(g)=\int_{N(\Q)\backslash N(\A)}\phi(ng)\;dn.
\end{equation}
This is a measurable, locally integrable function on $N(\A)\backslash
G(\A)$.

\subsection{}

Let $P$ be a parabolic subgroup. Then we denote by $\cA^2(P)$ the space of
automorphic forms on 
$N_P(\A) M_P(\Q)\ba G(\A)$
which are square integrable on 
$M_P(\Q)\backslash M_P(\A)^1\times K$.
This is the space of smooth functions
\begin{equation*}
\phi\colon N_P(\A) M_P(\Q)\ba G(\A)\rightarrow \C
\end{equation*}
which satisfy the following conditions: 
\begin{enumerate}
\item[i)] The span of the set of functions 
$$x\mapsto(z\phi)(xk),\quad x\in G(\A),$$
indexed by $k\in K$ and
$z\in Z({\mathfrak g}_{\C})$, is finite dimensional.

\item[ii)]

\begin{equation*} 
\parallel \phi\parallel^2=\int_K 
\int_{M_P(\Q)\backslash M_P(\A)^1} |\phi(mk)|^2 \;dm\; dk<\infty.
\end{equation*}
\end{enumerate}

Furthermore, an automorphic form $\phi\in \cA^2(P)$ is called cuspidal, if 
the following additional condition holds:
\begin{enumerate}
\item[iii)] For all standard parabolic subgroups $Q\subsetneqq P$, $\phi_Q=0$. 
\end{enumerate}

The subspace of
all cuspidal automorphic forms in $\cA^2(P)$ will be denoted by 
$\cA^2_0(P)$.

\subsection{} 
 
Given $\pi\in\Pi_{\di}(M_P(\A))_\xi$, let $\cA^2_\pi(P)$ be the subspace of 
$\cA^2(P)$ consisting of all functions $\phi$ such that for every 
$x\in G(\A)$, the function
$$\phi_x(m)=\phi(mx)\quad m\in M_P(\A),$$
belongs to the $\pi$-isotypical subspace of $L^2(M_P(\Q)\ba M_P(\A))_\xi$.
If $\pi\in\Pi(M_P(\A))$ is not contained in $\Pi_{\di}(M_P(\A))$, we put
$\cA^2_\pi(P)=0$. Let $K_f$ be an open compact subgroup of $G(\A_f)$. Then we
denote by $\cA^2_\pi(P)_{K_f}$ the subspace of all $K_f$-invariant functions
in $\cA^2_\pi(P)$. Furthermore, if $\sigma\in\Pi(K_\infty)$, then we denote
by $\cA^2_\pi(P)_{K_f,\sigma}$ the $\sigma$-isotypical subspace of  
$\cA^2_\pi(P)_{K_f}$.

\subsection{}

Let $\mX$ be the set of $W_0$ conjugacy classes of pairs
$(M_B,r_B)$, where $B$ is a  parabolic subgroup and $r_B$ is an
irreducible cuspidal automorphic representation of 
$M_B(\A)^1$. Let 
\begin{equation*}
L^2(M_P(\Q)\backslash M_P(\A)^1)=\bigoplus_{\chi\in{\mathfrak X}}
L^2(M_P(\Q)  \backslash M_P(\A)^1)_\chi
\end{equation*}
be the decomposition of $L^2(M_P(\Q)\ba M_P(\A)^1)$ introduced by Arthur in
\cite[Section 3]{A2}. 
Given $\chi\in\mX$, let $\cA^2_{\chi,\pi}(P)$ be the subspace of 
$\cA^2_\pi(P)$ consisting of all function $\phi$ such that for each $x\in
 G(\A)$, the restriction of $\phi_x$ to $M_P(\A)^1$  belongs to
 $L^2(M_P(\Q)\backslash M_P(\A)^1)_{\chi}$. 

If we identify $\Pi(M_P(\A)^1)$ with a subset of $\Pi(M_P(\A))$, then 
$\cA^2_{\chi,\pi}(P)$ is well defined for any $\pi\in\Pi(M_P(\A)^1)$. 
This is a space of functions on $N_P(\A)M_P(\A)A_P(\R)^0\ba G(\A)$.
The direct sum
$$\bigoplus_{\pi\in\Pi(M_P(\A)^1)}\cA^2_\chi(P,\pi)$$
is the space that was denoted by $A^2(P,\chi)$ in \cite{Mu3}.

\subsection{}
Let $\ov\cA^2(P)$ be the Hilbert space completion of 
$\cA^2(P)$. For any $\lambda\in \af_{P,\C}^*$ we have an induced 
representation $\rho(P,\lambda)$ of $G(\A)$ on 
$\ov\cA^2(P)$ which is defined by 
\begin{equation}\label{1.3}
(\rho(P,\lambda,y)\phi)(x)=\phi(xy)e^{(\lambda+\rho_P)(H_P(xy))}
e^{-(\lambda+\rho_P)(H_P(x))},
\end{equation}
for elements $x,y\in G(\A)$ and $\phi\in\ov\cA^2(P)$. The Hilbert space
completions $\ov\cA^2_\pi(P)$ and $\ov\cA^2_{\chi,\pi}(P)$ of the subspaces
$\cA^2_\pi(P)$ and $\cA^2_{\chi,\pi}(P)$, respectively, are invariant under
$\rho(P,\lambda)$ and we shall denote the restriction of $\rho(P,\lambda)$
to $\cA^2_\pi(P)$ (resp. $\cA^2_{\chi,\pi}(P)$) by $\rho_\pi(P,\lambda)$
(resp. $\rho_{\chi,\pi}(P,\lambda)$).

\subsection{}
Given any
irreducible unitary representation
$\pi$ of $M_B(\A)^1$, let $\lambda_\pi$ be the eigenvalue of the Casimir 
operator
of $M_B(\R)$, acting in the G{\aa}rding space 
$\H^\infty_{\pi_\infty }$
of the Archimedean constituent $\pi_\infty$ of $\pi$. For
$\chi\in{\mathfrak X}$ and $(M_B,r_B)\in\chi$, the Casimir eigenvalue
$\lambda_{r_B}$ depends only on the class $\chi$ and we denote it by
$\lambda_\chi$.

\section[Normalized intertwining operators]{Normalized  intertwining 
operators}
\setcounter{equation}{0}

Let $M,M_1\in\cL(M_0)$, $P\in\cP(M)$ and $P_1\in\cP(M_1)$. For each 
$s\in W(\af_M,\af_{M_1})$, $\phi\in \cA^2(P)$,  and $\lambda\in\af_{P,\C}^*$ 
such that $\Re(\lambda)\in(\af_P^*)^++\rho_P$, let 
$M_{P_1|P}(s,\lambda)\phi$  be defined by

\begin{equation}\label{2.1}
\begin{split}
M&_{P_1|P}(s,\lambda)\phi(x)=e^{-(s\lambda+\rho_{P_1})(H_{P_1}(x))}\\
&\int_{N_1(\A)\cap w_s N(\A)w_s^{-1}\backslash N_1(\A)} \phi(w_s^{-1}n_1x)
e^{(\lambda+\rho_P)(H_P(w_s^{-1}n_1x))}\;dn_1
\end{split}
\end{equation}
for $x\in G(\A)$.
The integral is absolutely convergent for $\lambda$ as above and admits
an analytic continuation to a meromorphic function of $\lambda\in\af_{P,\C}^*$
with values in the space of linear operators from $\cA^2(P)$ to $\cA^2(P_1)$.
This operator is the global intertwining operator
$$M_{P_1|P}(s,\lambda)\colon \cA^2(P)\to \cA^2(P_1).$$
Let $\pi\in\Pi_{\di}(M(\A))$ and  $\chi\in\mX$. Then $M_{P_1|P}(s,\lambda)$
maps the subspace $\cA_\pi^2(P)$ (resp. $A_{\chi,\pi}^2(P)$) to 
$A_{s\pi}^2(P_1)$ (resp.  $\cA^2_{\chi,s\pi}(P_1)$). 
The main functional equations are
\begin{equation}\label{2.2}
M_{P_2\mid P}(ts,\lambda)=M_{P_2\mid P_1}(t,s\lambda)M_{P_1|P}(s,\lambda)
\end{equation}
for $t\in W(\af_1,\af_2)$ and $s\in W(\af,\af_1)$.

By (1.4) and (1.5) of \cite{A4}, most of the considerations concerning
 intertwining operators can be reduced to the case  where $P_1$ and
$P$ have the same Levi component $M$, and $s=1$.

Thus, from now on we shall assume that $P,Q\in\cP(M)$ and we put
$$M_{Q|P}(\lambda):=M_{Q|P}(1,\lambda),\quad\lambda\in \af_{M,\C}^*.$$ 

Given $\pi\in\Pi(M(\A))$, let
$$M_{Q|P}(\pi,\lambda)\colon\cA^2_\pi(P)\to\cA^2_\pi(Q)$$
be the restriction of $M_{Q|P}(\lambda)$ to $\cA^2_\pi(P)$. We shall now 
express this operator in terms of local intertwining operators. 
Let $\pi_\lambda$ be the representation
of $P(\A)$ which is defined by 
$$\pi_\lambda(nm)=e^{\lambda(H_M(m))}\pi(m),\quad n\in N_P(\A),
\;m\in M_P(\A).$$
Let $(I_P^G(\pi_\lambda),\H_P(\pi))$ be the induced representation of $G(\A)$.
Similarly let $(I_Q^G(\pi_\lambda),\H_Q(\pi))$ be the representation of 
$G(\A)$ induced from $Q(\A)$. 
Let $\xi$ be a unitary character of $A_M(\R)^0$
and suppose that $\pi\in\Pi_{\di}(M(\A))_\xi$. We extend $\xi$ by 1 to a 
character of $M(\Q)A_M(\R)^0$. Then there is a canonical isomorphism
\begin{equation}\label{2.3}
j_P:\H_P(\pi)\otimes\Hom_{M(\A)}(\pi,I_{M(\Q)A_M(\R)^0}^{M(\A)}(\xi))\to
\ov\cA^2_\pi(P)
\end{equation}
of $G(\A)$-modules where $G(\A)$ acts on the left by 
$I_P^G(\pi_\lambda)\otimes\Id$. A similar isomorphism $j_Q$ exists with 
respect to
$Q$. Let $\H_P^0(\pi)$ (resp. $\H_Q^0(\pi)$) be the subspace of elements which
are right $K$-finite and left $Z(\gf_\C)$-finite. Using (\ref{2.3}), it
follows that $M_{Q|P}(\pi,\lambda)$ induces an intertwining operator
$$J_{Q|P}(\pi,\lambda)\colon\H^0_P(\pi)\to\H^0_Q(\pi)$$
such that
$$j_Q\circ (J_{Q|P}(\pi,\lambda)\otimes\Id)= M_{Q|P}(\pi,\lambda)\circ j_P.$$
It follows from (\ref{2.1}) that for $\Re(\lambda)\in(\af_P^*)^+ +\rho_P$,
 this operator is defined by the following  convergent integral 
\begin{equation}\label{2.4}
\begin{split}
J_{Q|P}  (\pi,\lambda)\phi(x)&=
e^{-(\lambda+\rho_Q)(H_Q(x))}\\
&\cdot\int_{N_Q(\A)\cap N_P(\A)\ba N_Q(\A)}
 \phi(nx)
e^{(\lambda+\rho_P)(H_P(nx))}\;dn.
\end{split}
\end{equation}
where $x\in G(\A)$ and $\phi\in\H_P^0(\pi)$. 

Let $v$ be any place of $\Q$ and let $(\pi_v,V_v)\in\Pi(M(\Q_v))$. Given
$\lambda\in{\af}^*_{M,\C}$, let $\pi_{v,\lambda}$ be the representation of
  $P(\Q_v)$ on $V_v$ defined by
\begin{equation*}
\pi_{v,\lambda}(n_v m_v)=\pi_v(m_v)e^{\lambda(H_M(m_v))},
\quad n_v\in N(\Q_v),\;m_v\in M(\Q_v).
\end{equation*}

Let $(I_P^G(\pi_{v,\lambda}),{\H}_P(\pi_v))$ denote induced representation. The
Hilbert space is the space of measurable functions
\begin{equation*}
\phi_v: N(\Q_v)\setminus G(\Q_v)\to V_v
\end{equation*}

such that
\begin{enumerate}
\item[1.] 
$\phi_v(m_v x_v)=\pi(m_v)\phi_v(x_v),\quad m_v\in M(\Q_v),x_v\in G(\Q_v);$

\item[2.] $\parallel \phi_v\parallel^2=\int_{K_{v}}\parallel
  \phi_v(k)\parallel^2_{V_v} dk< \infty$.
\end{enumerate}

Let ${\H}^0_P(\pi_v)\subset{\H}_P(\pi_v)$ be the subspace of $K_v$-finite
functions. Then the local intertwining operator
\begin{equation*}
J_{Q|P} (\pi_v,\lambda):{\H}^0_P(\pi_v)\to {\H}^0_Q(\pi_v)
\end{equation*}
is defined by
\begin{equation}\label{2.5}
\begin{split}
J_{Q|P} & (\pi_v,\lambda)\phi_v(x_v)=
e^{-(\lambda+\rho_Q)(H_Q(x_v))}\\
&\cdot\int_{N_Q(\Q_v)\cap N_P(\Q_v)\ba N_Q(\Q_v)}
 \phi_v(n_vx_v)
e^{(\lambda+\rho_P)(H_P(n_vx_v))}dn_v.
\end{split}
\end{equation}

The integral converges absolutely for $\Re(\lambda)\in(\af_P^*)^+ +\rho_P$  and
can be continued to a meromorphic function of $\lambda\in{\af}^*_{M,\C}$ with
values in the space of linear operators from ${\H}^0_P(\pi_v)$ to
${\H}^0_Q(\pi_v)$ \cite{Sh1}.

Now let $\pi\in \Pi(M(\A))$. Then $\pi$ is a restricted tensor product 
\begin{equation*}
\pi=\otimes_v\pi_v
\end{equation*}
where almost all $(\pi_v,V_v)$ are unramified, i.e., $\dim V_v^{K_v\cap
  M(\Q_v)}=1$ for almost all $v$.  Moreover, we have
\begin{equation*}
(I_P^G(\pi_\lambda),{\H}_P(\pi))\cong \left(\otimes_v I_P^G(\pi_{v,\lambda}),
\otimes_v{\H}_P(\pi_v)\right).
\end{equation*}
Let $\phi\in{\H}^0_P(\pi)$ and suppose that $\phi=\otimes_v\phi_v$. Observe
that each $\phi_v$ belongs to ${\H}^0_P(\pi_v)$ and for almost all $v$, 
$\phi_v$ is $K_v$--invariant. Comparing (\ref{2.4}) and (\ref{2.5}),
 it follows that
\begin{equation}\label{2.6}
J_{Q|P}(\pi,\lambda)\phi=\otimes_v(J_{Q|P}(\pi_v,\lambda)\phi_v)
\end{equation}
whenever the product on the right converges.

It is possible to normalize  local intertwining operators. Let $v$ be any 
valuation of $\Q$ and let $\pi_v\in\Pi(M(\Q_v))$.  It is proved in
\cite{A7}, \cite{CLL} that  there exist scalar valued 
meromorphic functions $r_{Q|P}(\pi_v,\lambda)$ of 
$\lambda\in\af_{P,\C}^*$  such that the
normalized intertwining operators
\begin{equation}\label{2.7}
R_{Q|P}(\pi_v,\lambda)=r_{Q|P}(\pi_v,\lambda)^{-1}J_{Q|P}(\pi_v,\lambda)
\end{equation}
satisfy the conditions $(R_1)-(R_8)$ of Theorem 2.1 of
\cite{A7}.
We recall some  of the properties satisfied by the normalized intertwining 
operators.

\begin{enumerate}
\item[(R.1)] If $S\in\cP(M)$, then
\begin{equation}\label{2.8}
R_{S|P}(\pi_v,\lambda)=R_{S|Q}(\pi_v,\lambda)R_{Q|P}(\pi_v,\lambda).
\end{equation}

\item[(R.2)]
\begin{equation}\label{2.9}
R_{Q|P}(\pi_v,\lambda)^*=R_{P|Q}(\pi_v,-\overline\lambda).
\end{equation}

\item[(R.3)]
Let $L\in\cL(M)$, $S\in\cP(L)$, and $Q,Q^\prime\in\cP^L(M)$. Then
\begin{equation}\label{2.10}
\left(R_{S(Q^\prime)|S(Q)}(\pi_v,\lambda)\phi\right)_k=R_{Q^\prime|Q}(\pi_v,
\lambda)\phi_k
\end{equation}
for any $\phi\in\H_P^0(S(R)(\Q_v))$ and $k\in K_v$.

\item[(R.4)]Let $v$ be a finite place. Suppose that $\pi_v$ is unramified, and that
$K_v$ is  very special. Let $\phi_{v}\in\H_P(\pi_v)$ be a function 
such that $\phi_{v}(k)=\phi_{v}(1)$ for all $k\in K_v$. Then in the compact
picture of the induced representation, one has
\begin{equation}\label{2.11}
R_{Q|P}(\pi_v,\lambda)\phi_{v}=\phi_{v}.
\end{equation}

\end{enumerate}

The functions $r_{Q|P}(\pi_v,\lambda)$ are called {\it normalizing
factors}. They satisfy similar properties.
We recall some of them. Given $P\in\cP(M)$, let $\Sigma_P^r$ be 
the set of  reduced roots of $(P,A_M)$. 

\begin{enumerate}
\item[(r.1)] For $\beta\in\Sigma_P^r$, let $M_\beta\in\cL(M)$ be such that
$$\af_{M_\beta}=\{H\in\af_M\mid\beta(H)=0\}$$
and let $P_\beta$ be the unique group in $\cP^{M_\beta}(M)$ whose simple 
root is $\beta$. Then
\begin{equation}\label{2.12}
r_{Q|P}(\pi_v,\lambda)=\prod_{\beta\in\Sigma^r_Q\cap\Sigma^r_{\overline P}}
r_{\overline P_\beta|P_\beta}(\pi_v,\lambda),
\end{equation}
Note that each $r_{\overline P_\beta|P_\beta}(\pi_v,\lambda)$ depends only
on the projection $\lambda(\beta^\vee)$.
\item[(r.2)]
If $\pi_v$ is an irreducible constituent of an induced representation 
$I_R^M(\sigma_v)$, where $\sigma_v\in\Pi_2(M_1(\Q_v))$, $R\in\cP^M(M_1)$, and
$M_1\subset M$, then
\begin{equation}\label{2.13}
r_{Q|P}(\pi_v,\lambda)=r_{Q(R)|P(R)}(\sigma_v,\lambda).
\end{equation}
\item[(r.3)] 
\begin{equation}\label{2.14}
r_{Q|P}(\pi_v,\lambda)r_{P|Q}(\pi_v,\lambda)
=J_{Q|P}(\pi_v,\lambda)J_{P|Q}(\pi_v,\lambda).
\end{equation}
\item[(r.4)]
\begin{equation}\label{2.15}
\overline{r_{Q|P}(\pi_v,\lambda)}=r_{P|Q}(\pi_v,-\overline\lambda).
\end{equation}
\item[(r.5)] 
If $v$ is a finite place of $\Q$, then $r_{Q|P}(\pi_v,\lambda)$ is a rational
function in the variables
$\{q_v^{-\lambda(\widetilde\beta)}\mid\beta\in\Sigma^r_Q\cap\Sigma^r_{\overline P}\}$, where the $\widetilde\beta$'s are suitably normalized ''coroots''. 
If $v=\infty$, then $r_{Q|P}(\pi_v,\lambda)$ is a rational function in the
variables
$\{\lambda(\beta^\vee)\mid\beta\in\Sigma^r_Q\cap\Sigma^r_{\overline P}
\}.$

\end{enumerate}

Now we return to global intertwining operators. Let 
$\pi\in\Pi_{\di}(M(\A))$. For $\phi\in\H_P^0(\pi)$ with $\phi=
\otimes_v\phi_v$ set 
\begin{equation}\label{2.16}
R_{Q|P}(\pi,\lambda)\phi=\otimes_v(R_{Q|P}(\pi_v,\lambda)\phi_v).
\end{equation}
Since  $\phi_v$ is $K_v$-
invariant for almost all $v$, it follows from (\ref{2.11}) that the right 
hand side  is actually a finite 
product and therefore, it converges for all $\lambda\in\af_{M,\C}^*$ which
are not poles of the local intertwining operators. In this way we get a 
 a meromorphic operator valued function 
$$R_{Q|P}(\pi,\lambda):\H_P^0(\pi)\to\H_Q^0(\pi)$$
of $\lambda\in\af_{M,\C}^*$. Using the isomorphism (\ref{2.3}) and the 
corresponding one for $Q$, we obtain a meromorphic operator valued function
$$N_{Q|P}(\pi,\lambda)\colon\cA^2_\pi(P)\to\cA^2_\pi(Q)$$
of $\lambda\in\af_{M,\C}^*$ such that
\begin{equation}\label{2.17}
j_Q\circ N_{Q|P}(\pi,\lambda)= (R_{Q|P}(\pi,\lambda)\otimes\Id)\circ j_P.
\end{equation}
Furthermore, put
\begin{equation}\label{2.18}
r_{Q|P}(\pi,\lambda)=\prod_vr_{Q|P}(\pi_v,\lambda).
\end{equation}
By (R.4) it follows that for $\phi$ as above, we have
$$J_{Q|P}(\pi_v,\lambda)\phi_v=r_{Q|P}(\pi_v,\lambda)\phi_v$$
for almost all $v$. Therefore, the infinite product (\ref{2.18}) converges
in the domain of absolute convergence of the infinite product (\ref{2.6})
and for $\lambda$ in this domain we have
\begin{equation}\label{2.19}
M_{Q|P}(\pi,\lambda)=r_{Q|P}(\pi,\lambda) N_{Q|P}(\pi,\lambda).
\end{equation}
Since both $M_{Q|P}(\pi,\lambda)$ and $N_{Q|P}(\pi,\lambda)$  are meromorphic
functions of $\lambda\in\af_{M,\C}^*$, it follows that $r_{Q|P}(\pi,\lambda)$
admits a meromorphic continuation to $\af_{M,\C}^*$. The meromorphic function
$r_{Q|P}(\pi,\lambda)$ is the global normalizing factor and $N_{Q|P}(\pi,
\lambda)$ is the normalized global intertwining operator.

Using (\ref{2.12}), (\ref{2.14}), (\ref{2.15}) and the functional equations 
(\ref{2.2}), it follows that $r_{Q|P}(\pi,\lambda)$
has the following properties

\begin{enumerate} 
\item 
\begin{equation}\label{2.20}
r_{Q|P}(\pi,\lambda)r_{P|Q}(\pi,\lambda)=1.
\end{equation}

\item
\begin{equation}\label{2.21}
\overline{r_{Q|P}(\pi,\lambda)}=r_{P|Q}(\pi,-\overline\lambda).
\end{equation}

\item
For each $\beta\in\Sigma^r_Q\cap\Sigma^r_{\overline P}$ let $P_\beta$ 
be as in (\ref{2.12}). Then

\begin{equation}\label{2.22}
r_{Q|P}(\pi,\lambda)=\prod_{\beta\in \Sigma^r_Q\cap\Sigma^r_{\overline P}}
r_{\overline P_\beta|P_\beta}(\pi,\lambda).
\end{equation}
Note that $r_{\overline P_\beta|P_\beta}(\pi,\lambda)$ depends only on the 
projection $\lambda(\beta^\vee)$.
\end{enumerate}

\section[Local normalizing factors]{Local normalizing factors}
\setcounter{equation}{0}
In this section we shall investigate the local normalizing factors in more 
detail. In particular, we shall study their logarithmic derivatives. To begin
with, we  recall the construction of the normalizing factors.

First assume that $v$ is a finite valuation. Then the existence of 
normalizing factors such that Theorem 2.1 of \cite{A7} holds has been 
verified by Langlands  in \cite{CLL}, Lecture 15. Let $\pi_v\in
\Pi(M(\Q_v))$. The local normalizing factors $r_{Q|P}(\pi_v,\lambda)$ have to 
satisfy (2.1)-(2.3) in \cite{A7}. Therefore, it suffices to
define them when $\dim(\af_M/\af_G)=1$ and $\pi_v$ is square integrable
modulo $A_G$. Assume for the moment that these conditions are satisfied.
Let $P\in\cP(M)$ and let $\alpha$ be the unique simple root 
of $(P,A_M)$. Then Langlands has shown that there exists a rational function 
$V_P(\pi_v,z)$ of one variable such that
\begin{equation}\label{3.1}
r_{\ov P|P}(\pi_v,\lambda)=V_P(\pi_v,q_v^{-\lambda(\widetilde\alpha)}),
\end{equation}

where $\widetilde\alpha\in\af_M$ is independent of $\pi_v$. We recall the
definition of $V_P(\pi_v,z)$. Suppose that $P_v$ is a parabolic subgroup of
$G$ defined over $\Q_v$ and let $M_v$ be a Levi component of $P_v$ over $\Q_v$.
Denote by $A_{M_v}$ the split component of the center of $M_v$. Set 
$$\af_{M_v}=\Hom(X(M_v)_{\Q_v},\R)$$
and
$$\af_{M_v}^*=X(M_v)_{\Q_v}\otimes\R.$$
Let 
$$H_{M_v}:M_v(\Q_v)\to\af_{M_v}$$
be defined by 
$$\hskip70pt q_v^{\langle H_{M_v}(m_v),\chi\rangle}=|\chi(m_v)|_v,\quad\chi\in X(M_v)_{\Q_v},\;m_v\in M_v(\Q_v).$$
 Given 
$\pi\in\Pi(M_v(\Q_v))$ and $\lambda\in i\af^*_{M_v}$, let $\pi_\lambda$ denote 
the representation defined by 
$$\pi_\lambda(m_v)=\pi(m_v)e^{\lambda(H_{M_v}(m_v))},\quad m_v\in M_v(\Q_v).$$
Let
$$\af_{M_v,\pi}^\vee=\{\lambda\in i\af_{M_v}^*\mid \pi_\lambda\cong\pi\}$$
denote the stabilizer of $\pi$ with respect to this action of $i\af_{M_v}^*$.
 Then $\af_{M_v,\pi}^\vee$ is a lattice in $i\af_{M_v}^*$
and the orbit $\ho_\pi$ of $\pi$ is equal to $i\af_{M_v}^*/\af_{M_v,\pi}^\vee$.
Let 
$$\af_{M_v,\Q_v}=H_{M_v}(M_v(\Q_v)),
\quad \widetilde\af_{M_v,\Q_v}=H_{M_v}(A_{M_v}(\Q_v)).$$
Then $\af_{M_v,\Q_v}$ and $\widetilde\af_{M_v,\Q_v}$ are lattices in 
$\af_{M_v}$. Given a real vector space $V$ and a closed subgroup $V_1$ of $V$,
$$V_1^\vee=\Hom(V_1,2\pi i\Z)\subset i V^*.$$
let us agree to set
Let $\af_{M_v,\pi}\subset\af_{M_v}$ be the dual lattice to 
$\af_{M_v,\pi}^\vee$. Then
$$\widetilde\af_{M_v,\Q_v}\subset\af_{M_v,\pi}\subset\af_{M_v,\Q_v}.$$
Set 
$$L_{M_v}=\left(\af_{M_v,\Q_v}+\af_{G_v}\right)/\af_{G_v},\quad
\widetilde L_{M_v}=\left(\widetilde\af_{M_v,\Q_v}+\af_{G_v}\right)/\af_{G_v},$$
and
$$L(\pi)=\left(\af_{M_v,\pi}+\af_{G_v}\right)/\af_{G_v}.$$
Then $L_{M_v}$, $\widetilde L_{M_v}$, and $L(\pi)$ are lattices in 
$\af_{M_v}^{G_v}=\af_{M_v}/\af_{G_v}$.

Suppose that $P_v$ is a maximal parabolic subgroup, that is 
$\dim\af_{M_v}^{G_v}=1$. Then there exists $\alpha(\pi)\in\af_{M_v}$ such that
$$L(\pi)=\frac{\log q}{2\pi}\Z(\alpha(\pi)).$$

In  \cite{Si1} Silberger has shown that for a supercuspidal representation 
$\pi$ there exists a rational function  $\widetilde U_{P_v}(\pi,z)$ such that 
the 
Plancherel measure $\mu(\pi,\lambda)$ satisfies
\begin{equation}\label{3.2}
\mu(\pi,\lambda)=\widetilde U_{P_v}\left(\pi,q^{-\lambda(\alpha(\pi))}\right).
\end{equation}
Let $\widetilde\alpha\in\af_{M_v}$ be such that
\begin{equation*}
L_{M_v}=\frac{\log q}{2\pi}\Z(\widetilde\alpha).
\end{equation*}
Since $L(\pi)\subset L_{M_v}$, there exists $k(\pi)\in\Z$ such that
$\alpha(\pi)=k(\pi)\widetilde\alpha$. Let
\begin{equation}\label{3.3}
U_{P_v}(\pi,z)=\widetilde U_{P_v}\left(\pi,z^{k(\pi)}\right).
\end{equation}
Then
$$\mu(\pi,\lambda)=U_{P_v}\left(\pi,q^{-\lambda(\widetilde\alpha)}\right).$$
Now suppose that $P_v$ is arbitrary, but $\pi$ is still supercuspidal. For
each reduced root $\alpha\in\Sigma_r(P_v,A_v)$ let $A_\alpha$ denote the 
largest subtorus of $A_v$ which lies in the kernel of the root character of 
$\alpha$. Let $M_\alpha$ denote the centralizer of $A_\alpha$. Let ${}^*P_\alpha=P_v\cap M_\alpha$. Then ${}^*P_\alpha=M_vN_\alpha$. Let 
$\mu_\alpha(\pi,\lambda)$ be the Plancherel measure with respect to 
$(M_\alpha,{}^*P_\alpha)$.
According to \cite{HC3}, Theorem 24, there exist  constants 
$\gamma=\gamma(G/M)$ and $\gamma_\alpha=\gamma(M_\alpha/M)$, 
$\alpha\in\Sigma_r(P_v,A_v)$,  such that
\begin{equation}\label{3.4}
\gamma^{-2}\mu(\pi,\lambda)=\prod_{\alpha\in\Sigma_r(P_v,A_v)}
\gamma_\alpha^{-2}\mu_\alpha(\pi,\lambda),
\end{equation}
Hence if
$\{\widetilde\alpha\mid\alpha\in\Delta_{P_v}\}$ is a set of generators of
the lattice $L_{M_v}$, then $\mu(\pi,\lambda)$ is a rational function in the
variables $\{q^{-\lambda(\widetilde\alpha)}\mid\alpha\in\Delta_{P_v}\}$.
Finally, by Theorem 1 of \cite{Si2}, this can be extended to all discrete
series representations of $M_v(\Q_v)$. 

Now let $P=MN$ be a maximal parabolic subgroup of $G$ defined over $\Q$.
Then $X(M)_\Q\subset X(M)_{\Q_v}$ induces an embedding $\af_M^*\subset
\af^*_{M_v}$ and by the above, it follows that there exists $\widetilde\alpha
\in\af_M$ and a rational function $U_P(\pi,z)$ such that
\begin{equation}\label{3.5}
\mu(\pi,\lambda)=U_P\left(\pi,q^{-\lambda(\widetilde\alpha)}\right),
\end{equation}
for all $\pi\in\Pi_2(M(\Q_v))$, $\lambda\in\af_{M,\C}^*$.
As shown by Langlands \cite{CLL}, the rational function $U_P(\pi,z)$ has the 
form
$$U_P(\pi,z)=a\frac{\prod_{i=1}^r(1-\alpha_iz)(1-\overline\alpha_i^{-1}z)}
{\prod_{i=1}^r(1-\beta_iz)(1-\overline\beta_i^{-1}z)},$$
where the $\alpha_i$'s and $\beta_i$'s satisfy $|\alpha_i|\le 1$,
$|\beta_i|\le 1$, $i=1,...,r$, and $a$ is a certain constant.
Then the 
rational function $V_P(\pi,z)$ in (\ref{3.1}) is defined  by
\begin{equation}\label{3.6}
V_P(\pi,z)=b\frac{\prod_{i=1}^r(1-\beta_iz)}{\prod_{i=1}^r(1-\alpha_iz)}
\end{equation}
for a suitable constant $b$.
In particular, it follows that $2r$ is the number of poles of $U_P(\pi,z)$. 
For our applications we need a bound for $r$. This is done in the following 
lemma.
\begin{lem}\label{l3.1}
Let $M\in\L(M_0)$ be such that $\dim(\af_M/\af_G)=1$.
There exists $C>0$ such that for all $P\in\cP(M)$ and all 
$\pi\in\Pi_2(M(\Q_v))$ the number of poles of the rational function 
$V_P(\pi,z)$ is less than or equal to $C$.
\end{lem}

\begin{proof} 
Let $P_v$ be a maximal parabolic subgroup of $G$ defined over $\Q_v$ and let
$\pi\in\Pi(M_v(\Q_v))$ be supercuspidal. By Theorem 1.6 of \cite{Si1}, the
rational function $\widetilde U_{P_v}(\pi,z)$ in (\ref{3.2}) has at most 4 
poles.
Now observe that $\widetilde L_{M_v}\subset L(\pi)\subset L_{M_v}$ and 
$L_{M_v}/\widetilde L_{M_v}$ is finite. This implies that the number of poles
of the rational function $U_{P_v}(\pi,z)$ defined by (\ref{3.3}) is bounded
by a constant which is independent of $\pi$. The general case is reduced to
this one using the product formula (\ref{3.4}) and Theorem 1 of \cite{Si2}. 
\end{proof}

Using (2.1)-(2.3) of \cite{A7}, the local normalizing factors can be
defined for all $M\in\cL(M_0)$, $P,Q\in\cP(M)$ and $\pi_v\in\Pi(M(\Q_v))$.

Next suppose that $v=\infty$. In this case  the existence of  normalizing 
factors such that Theorem 2.1 in \cite{A7} holds has been established
by Arthur \cite{A7}. The definition is as follows. Let ${}^LM$ be the $L$-group
of $M$ and 
 let 
$\rho=\tilde\rho_{Q|P}$ be the contragredient representation of the
adjoint representation $\rho_{Q|P}$ of ${}^LM$ on the complex vector 
space ${}^L\nf_Q\cap{}^L\nf_P\backslash {}^L\nf_Q$. Let $L(s,\pi,\rho)$
be the $L$-factor attached to $\pi$ and $\rho=\tilde\rho_{Q|P}$. Then Arthur 
has shown in \cite{A7} that the functions
\begin{equation}\label{3.7}
r_{Q|P}(\pi,\lambda):=\frac{L(0,\pi_\lambda,\rho)}{L(1,\pi_\lambda,\rho)}
\end{equation}
satisfy all properties required by normalizing factors.  We briefly recall the
 definition of the $L$-function
and refer to \cite[pp.33-35]{A7} for more details.

To any  $\pi\in\Pi(M(\R))$ and $\lambda\in\af_{P,\C}^*$,
there corresponds a map
$$\phi_\lambda:W_\R\to {}^LM$$
from  the Weil group of $\R$ to the $L$-group of $M$, which is uniquely
determined by $\pi_\lambda$ up to conjugation by ${}^LM^0$ \cite{L3}.
 Let
\begin{equation}\label{3.8}
\rho\cdot\phi_\lambda=\bigoplus_\tau\tau_\lambda
\end{equation}
be the decomposition of $\rho\cdot\phi_\lambda$ into irreducible
representations of $W_\R$. Then by definition

$$L(s,\pi_{\lambda},\rho)=L(s,\rho\cdot\phi_\lambda)=
\prod_\tau L(s,\tau_\lambda).$$

So it remains to describe the $L$-factors 
$L(0,\tau_\lambda)L(1,\tau_\lambda)^{-1}$. To this end
let $T\subset M$ be a maximal torus over $\R$ whose real split component is
$A_M$. Let $\langle\lambda,\lambda^\vee\rangle$ denote the canonical 
pairing $X^*(T)\times X_*(T)\to\Z$ between the space $X^*(T)$ of characters and
the space $X_*(T)$ of one-parameter subgroups of $T$. 
Let $\Sigma_P(G,T)$ be the set of roots of $(G,T)$ which restrict to
roots of  $(P,A_M)$. The Galois group Gal$(\C/\R)$ acts on  $\Sigma_P(G,T)$.
Let $\ov \sigma$ be the action of the nontrivial element of Gal$(\C/\R)$.
The eigenspaces of $\tilde \rho_{Q|P}(\phi_\lambda(\C^*))$ are the root
spaces of $\{-\alpha^\vee|\alpha\in\Sigma_P(G,T)\}$ and the irreducible
constituents $\tau_\lambda$ of $\tilde\rho_{Q|P}\cdot\phi_\lambda$ 
correspond to orbits of $\ov\sigma$ in $\Sigma_P(G,T)$. Furthermore, the
map $\phi\colon W_\R\to {}^LM$  determines elements 
$\mu,\nu\in X^*(T)\otimes\C$ with $\mu-\nu\in X^*(T)$. 
Let 
$$\Gamma_\C(z):=2(2\pi)^{-z}\Gamma(z)\;\mbox{and}\;\Gamma_\R(z)=\pi^{-z/2}
\Gamma(z/2).$$
If a two-dimensional
constituent $\tau_\lambda$ corresponds to a pair $\{\alpha,\ov\sigma\alpha\}$ 
of complex roots, then $\tau_\lambda$ is induced from the quasi-character
$$z\mapsto z^{\langle\mu+\lambda,\alpha^\vee\rangle}
\ov z^{\langle\nu+\lambda,\alpha^\vee\rangle}$$
of $\C^*$. Replacing $\alpha^\vee$ by $\ov\sigma\alpha^\vee$ if necessary,
we can assume that $\langle\ov\sigma\mu-\mu,\alpha^\vee\rangle$ is a 
nonpositive integer. 
Then
\begin{equation}\label{3.9}
\frac{L(0,\tau_\lambda)}{L(1,\tau_\lambda)}=\frac{\Gamma_\C(\langle \mu+\lambda,\alpha^\vee\rangle)}{\Gamma_\C(\langle \mu+\lambda,\alpha^\vee\rangle+1)}.
\end{equation}
The one-dimensional constituents $\tau_\lambda$ correspond to real roots 
$\alpha_0$ in $\Sigma_P(G,T)$. There is at most one of these. If $\alpha_0$ 
exists, then $\tau_\lambda$ is induced from the quasi-character of $\R^*$ 
$$x\mapsto \left(\frac{x}{|x|}\right)^{-N_0}|x|^{\langle\mu+\lambda,\alpha_0^\vee\rangle},$$
where $N_0\in\{0,1\}$.
In this case
\begin{equation}\label{3.10}
\frac{L(0,\tau_\lambda)}{L(1,\tau_\lambda)}=\frac{\Gamma_\R(\langle \mu+\lambda,\alpha^\vee_0\rangle+N_0)}{\Gamma_\R(\langle \mu+\lambda,\alpha^\vee_0\rangle+N_0+1)}.
\end{equation}

{\bf Remark:}
 It has been conjectured by Langlands \cite[p.282]{L1} that for 
any local field, intertwining operators can be normalized by $L$-functions.
For $\GL(n)$ this was proved by Shahidi \cite{Sh2}. Namely, let $P$ be a 
standard maximal parabolic subgroup of $\GL(n)$. Then $\af_{P,\C}^*\cong\C^2$
and $M_P=\GL(n_1)\times\GL(n_2)$ for some decomposition $n=n_1+n_2$. Let $F$
be a local non-Archimedean field and let $\psi$ be a non-trivial additive \
character of $F$. Let $\pi_1\otimes\pi_2$ be an irreducible unitary 
representation of $M_P(F)=\GL(n_1,F)\times\GL(n_2,F)$. Let $L(z,\pi_1\times
\check{\pi}_2)$ and $\varepsilon(z,\pi_1\times\check{\pi}_2,\psi)$ be the 
Rankin-Selberg $L$-function and the $\varepsilon$-factor defined by Jacquet,
Piatetski-Shapiro and Shalika \cite{JPS}. Then the normalizing factor 
$r_{\ov P|P}(\pi,s)$, $s=(s_1,s_2)$, can be 
chosen to be
$$r_{\ov P|P}(\pi,s)=\frac{L(s_1-s_2,\pi_1\times\check{\pi}_2)}
{\varepsilon(s_1-s_2,\pi_1\times\check{\pi}_2,\psi)L(1+s_1-s_2,\pi_1\times
\check{\pi}_2)}.$$
In \cite{Sh4}, this has been generalized to quasi-split groups and generic
representations $\pi$.

\smallskip
We can now estimate the logarithmic derivatives of the normalizing factors.
First we consider the case of a finite place $v$. Let $q_v$ be the number of 
elements of the residue field of $\Q_v$ (which in our case is a prime number).
\begin{lem}\label{l3.2}
Let $\alpha\in\ov{(\af_{M}^*)^+}$.
There exist $C,c>0$ such that for every finite valuation  $v$ and every 
$\pi_v\in\Pi(M(\Q_v))$ we have
\begin{equation}\label{3.11}
\bigg|r_{Q|P}(\pi_v,z\alpha)^{-1}\cdot\frac{d}{dz}r_{Q|P}(\pi_v,z\alpha)
\bigg|
\le Cq_v^{-2}
\end{equation}
for $\Re(z)\ge c$.
\end{lem}

\begin{proof} 

First we assume that $\dim\af_M/\af_G=1$ and $\pi_v\in\Pi(M(\Q_v))$ is square 
integrable modulo 
$A_M(\Q_v)$.  Let $\alpha$ be the unique simple root of $(P,A)$.
Then $r_{\ov P|P}(\pi_v,\lambda)$ is given by (\ref{3.1}). Let $\lambda=z\alpha$,
$z\in\C$. 
Then $\lambda(\alpha^\vee)=z$ and by (\ref{3.6}), it follows that 
\begin{equation*}
\begin{split}
r_{\ov P|P}(\pi_v,z\alpha)^{-1}\frac{d}{dz}r_{\ov P|P}(\pi_v,z\alpha)=&
\log( q_v)\,q_v^{-z}\\
&\cdot\sum_{i=1}^r\left\{\frac{\beta_i}
{1-\beta_iq_v^{-z}}-\frac{\alpha_i}
{1-\alpha_iq_v^{-z}}\right\}.
\end{split}
\end{equation*}
Recall that 
the $\alpha_i$'s and $\beta_i$'s satisfy 
$|\alpha_i|\le 1$, $|\beta_i|\le 1$, $i=1,...,r$. Moreover, by  
Lemma \ref{l3.1} there exists $C_1\ge0$, which is independent of $\pi_v\in
\Pi(M(\Q_v))$, such that  $r\le C_1$. Therefore, for
$\Re(z)>3$ we obtain
\begin{equation}\label{3.12}
\Big|r_{\ov P|P}(\pi_v,z\alpha)^{-1}\frac{d}{dz}r_{\ov P|P}(\pi_v,z\alpha)\Big|\le
C_1\frac{\log(q_v)q_v^{-\Re(z)}}{1-q_v^{-\Re(z)}}<C_2q_v^{-2}.
\end{equation}
Now let $M\in\cL(M_0)$ be arbitrary, but still assume that $\pi_v$ is square 
integrable modulo $A_M(\Q_v)$. Let $P,Q\in\cP(M)$. For each $\beta\in\Sigma_P^r$
let $M_\beta\in\cL(M)$ be such that
$$\af_{M_\beta}=\{H\in\af_M\mid\beta(H)=0\}.$$
Then $\dim\af_M/\af_{M_\beta}=1$. Let $P_\beta$ be the unique group in 
$\cP^{M_\beta}(M)$ whose simple root is $\beta$. Furthermore, let $\alpha\in
(\af_M^*)^+$, $\nu\in\af_{M,\C}^*$ and $z\in\C$. Then by (\ref{2.12}) we get
\begin{equation*}
\begin{split}
&r_{Q|P}(\pi_v,z\alpha+\nu)^{-1}\cdot\frac{d}{dz}r_{Q|P}(\pi_v,z\alpha+\nu)\\
&=\sum_{\beta\in\Sigma^r_P\cap\Sigma^r_{\ov Q}}r_{\ov P_\beta|P_\beta}(\pi_v,
\langle z\alpha+\nu,\beta^\vee\rangle\beta)^{-1}
\cdot\frac{d}{dz}r_{\ov P_\beta|P_\beta}(\pi_v,
\langle z\alpha+\nu,\beta^\vee\rangle\beta).
\end{split}
\end{equation*}
By assumption we have $\langle\alpha,\beta^\vee\rangle\ge 0$ for every 
$\beta\in\Sigma_P^r$. If 
$\langle\alpha,\beta^\vee\rangle=0$, then the corresponding logarithmic 
derivative vanishes. Suppose that $a:=\langle\alpha,\beta^\vee\rangle>0$. Let
$c_0>0$ be such that $\parallel\beta\parallel\le c_0$ for all $\beta\in
\Sigma^r_P$.
Then
$$\Re\left(z\langle\alpha,\beta^\vee\rangle+
\langle\nu,\beta^\vee\rangle\right)\ge a\Re(z)-c_0
\parallel\nu\parallel$$
and it follows from  
(\ref{3.7}) that there exist $C,c>0$, depending on $\alpha$, $c_0$ and 
$\parallel\nu\parallel$, such that 
\begin{equation}\label{3.13}
\Big|r_{Q|P}(\pi_v,z\alpha+\nu)^{-1}\cdot\frac{d}{dz}r_{Q|P}(\pi_v,z\alpha+\nu)
\Big|\le Cq_v^{-2}
\end{equation}
for all $\pi_v\in\Pi_2(M(\Q_v))$ and $\Re(z)\ge c$.

Next assume that $\pi_v$ is tempered. Then $\pi_v$ is an irreducible 
constituent of an induced representation $I_R^M(\tau_v)$, where 
$R\in\cP^M(M_1)$, $M_1\subset M$   and $\tau_v\in\Pi_2(M_R(\Q_v))$. Then 
$I_P^G(\pi_{v,\lambda})$ is canonically
isomorphic to a subrepresentation of $I_{P(R)}^G(\tau_{v,\lambda})$ and by
(\ref{2.13}) we have
$$r_{Q|P}(\pi_v,\lambda)=r_{Q(R)|P(R)}(\tau_v,\lambda),$$
where $P(R)\subset P$, $Q(R)\subset Q$. Now recall that there is a 
canonical inclusion $\af_P^*\subset \af_{P(R)}^*$ and with respect to this
inclusion, we have $(\af_P^*)^+\subset\ov{(\af_{P(R)}^*)^+}$. Thus 
$\alpha$ can be identified with an element of $\ov{(\af_{P(R)}^*)^+}$. Hence 
 (\ref{3.11}) holds for all tempered $\pi_v\in\Pi(M(\Q_v))$.

Now  let $\pi_v$ be an arbitrary representation in $\Pi(M(\Q_v))$. Then  
 $\pi_v$ is the Langlands quotient $J_R^M(\tau_v,\mu)$ of a representation 
$I^M_R(\tau_v,\mu)$, where
$M_R$ is an admissible Levi subgroup of $M$, $\tau_v\in\Pi(M_R(\Q_v))$ is a
tempered representation, and $\mu$ is a point in the chamber of 
$(\af_R^M)^*=\af_R^*/\af_M^*$ attached to $R$ \cite{Si3}. Set 
$\Lambda=\mu+\lambda$. 
Then, as
explained in \cite[p.30]{A7}, we have
\begin{equation}\label{3.14}
r_{Q|P}(\pi_v,\lambda)= r_{Q(R)|P(R)}(\tau_v,\Lambda).
\end{equation}
Let $\rho_v\in\af_M^*$ be defined by
$$\delta_P(a)^{1/2}=q_v^{\rho_v(H(a))},\quad a\in A_M.$$
Then it follows from Theorem 3.3 of Chapter XI of \cite{BW} that
$$\langle\mu,\beta^\vee\rangle\le \langle\rho_v,\beta^\vee\rangle,\quad 
\beta\in\Phi(R,A_R).$$
Since $\mu$ belongs to $\af_R^*/\af_M^*$, it follows that
 $\parallel\mu\parallel\le \parallel\rho_v\parallel$. Let 
$\alpha\in(\af_P^*)^+$. As observed above,
$\alpha$ can be identified with an element of $\ov{(\af_{P(R)}^*)^+}$. Hence,
combining (\ref{3.13}) and (\ref{3.14}) the desired estimation (\ref{3.11})
follows.
\end{proof}

Next we consider the infinite place. Let $\pi\in\Pi(M(\R))$ and let 
$\phi\colon W_\R\to{}^LM$ be the map associated to $\pi$. Let 
  $\mu,\nu\in X^*(T)\otimes\C$ be the elements determined by the map $\phi$
(see \cite{L3}, \cite[p.34]{A7}).
To indicate the dependence on $\pi$, we shall write $\mu_\pi$ and $\nu_\pi$. 
We note that there is a canonical  injection of the space
$$\af_{M,\C}^*=X^*(M)_\Q\otimes\C$$
into $X^*(T)\otimes\C$. 
\begin{lem}\label{l3.3}
Let $\beta\in(\af_M^*)^+$. 
There exist $C,c>0$ such that 
\begin{equation}\label{3.15}
\bigg|r_{Q|P}(\pi,z\beta)^{-1}\cdot
\frac{d}{dz}r_{Q|P}(\pi,z\beta)
\bigg|
\le C
\end{equation}
for for all $\pi\in\Pi(M(\R))$ and all $z\in\C$ with $\Re(z)\ge c$.
\end{lem}

\begin{proof}
First assume that  $\pi\in\Pi(M(\R))$ is tempered. 
As explained  above, the normalizing factor $r_{Q|P}(\pi,\lambda)$ is a  
product of finitely many meromorphic functions each of which is either of the 
form (\ref{3.9}) or (\ref{3.10}). So it suffices to consider the logarithmic 
derivative of the Gamma factors. Recall that for $\Re(z)>0$ the following
formula holds 
\begin{equation*}
\frac{\Gamma^\prime(z+1)}{\Gamma(z+1)}=\frac{1}{2z}+\log z-\int_0^\infty
\left\{\frac{1}{2}-\frac{1}{u}+\frac{1}{e^u-1}\right\}e^{-uz}\,du
\end{equation*}
\cite[p.248]{Wh}. Let $0\le a\le 1$ and  $\Re(z)>2$. Then we get
\begin{equation*}
\begin{split}
\bigg|\frac{\Gamma^\prime(z)}{\Gamma(z)}
-\frac{\Gamma^\prime(z+a)}{\Gamma(z+a)}\bigg|&\le
\frac{2}{\Re(z)}+\log\bigg|1+\frac{a}{z-1}\bigg|+\frac{\pi}{2}\\
&\quad\,+2\int_0^\infty\bigg|\frac{1}{2}
-\frac{1}{u}+\frac{1}{e^u-1}\bigg|e^{-u\Re(z)/2}\,du.
\end{split}
\end{equation*}
Hence there exists $C>0$ such that 
\begin{equation}\label{3.16}
\bigg|\frac{\Gamma^\prime(z)}{\Gamma(z)}-\frac{\Gamma^\prime(z+a)}{\Gamma(z+a)}\bigg|\le C\quad\mbox{for}\quad\Re(z)>2.
\end{equation}

Let $\beta\in(\af_{M}^*)^+$, $\nu\in\af_{M,\C}^*$ and  
$\alpha\in\Sigma_P(G,T)$. Since 
$\alpha\upharpoonright\af_M$ is a root of $(P,A)$, it follows that 
$\langle\beta,\alpha^\vee\rangle>0$. Let $\alpha\in\Sigma_P(G,T)$ be such that
$\langle\ov\sigma\mu_\pi-\mu_\pi,\alpha^\vee\rangle\le0$. Then we have
$\Re\langle\mu_\pi,\alpha^\vee\rangle\ge0$. Hence
$$\Re\langle\mu_\pi+z\beta+\nu,\alpha^\vee\rangle\ge
\langle\beta,\alpha^\vee\rangle\Re(z)
-\parallel\alpha\parallel\cdot\parallel\nu\parallel.$$ 
Using (\ref{3.16}) together with (\ref{3.9}) and (\ref{3.10}), it follows that 
there exist constants  $C,c>0$ such that
\begin{equation}\label{3.17}
\bigg|\frac{L(1,\tau_{z\beta+\nu})}{L(0,\tau_{z\beta+\nu})}\cdot
\frac{d}{dz}\left(\frac{L(0,\tau_{z\beta+\nu})}{L(1,\tau_{z\beta+\nu})}
\right)\bigg|\le C
\end{equation}
for $\Re(z)\ge c(1+\parallel\nu\parallel)$ and all $\pi\in\Pi(M(\R))$, where 
$\tau_\lambda$ and 
$\pi_\lambda$ are related by (\ref{3.8}). 

Now let $\pi$ be an arbitrary representation in $\Pi(M(\R))$. Then there
exist a parabolic subgroup $R$ of $M$, a tempered representation $\tau$ of
$M_R(\R)$ and a point $\xi$ in the positive chamber of $(\af_R^*/\af_M^*)$
attached to $R$ such that $\pi$ is equivalent to the Langlands quotient
$J_R(\tau,\xi)$ \cite{L3}. Set $\Lambda=\xi+\lambda$. Then, as explained in 
\cite[p.30]{A7}, we have
$$r_{Q|P}(\pi,\lambda)=r_{Q(R)|P(R)}(\tau,\Lambda).$$
Moreover, by Theorem 5.2 of Chapter IV of \cite{BW} it follows that
$$|\Re\langle\xi,\alpha^\vee\rangle|\le 4\parallel\rho_P\parallel\quad
\mbox{for all}\quad \alpha\in\Delta_P.$$
Together with (\ref{3.17}) this implies the claimed result.
\end{proof}

\section[Poles of normalizing factors]{Poles of global normalizing factors}
\setcounter{equation}{0}

Let $M\in\cL(M_0)$ and $P,Q\in\cP(M)$. Let $\pi\in\Pi_{\di}(M(\A))$ with 
$\pi=\otimes_v\pi_v$. Then by \S 2  the 
infinite product 
$$r_{Q|P}(\pi,\lambda)=\prod_v r_{Q|P}(\pi_v,\lambda)$$
is absolutely convergent in  some chamber and admits an analytic extension to 
a  meromorphic function of $\lambda\in\af_{M,\C}^*$. In  this section we 
shall study the 
poles of $r_{Q|P}(\pi,\lambda)$.

Recall that a function $f\colon \C^N\to\C$ is called a meromorphic function
of order $p\ge0$, if $f$ can be written as a quotient $f=g_1/g_2$ of two 
entire functions $g_i:\C^N\to\C$, $i=1,2$, satisfying 
$$|g_i(z)|\le C e^{c\parallel z\parallel^p},\quad z\in\C^N,\;i=1,2,$$
for certain constants $C,c>0$. With this definition we have the following
proposition.

\begin{prop}\label{p4.1}  Let $n=\dim G(\R)/K_\infty$. For all 
$\pi\in\Pi_{\di}(M(\A))$, the normalizing factor 
$r_{Q|P}(\pi,\lambda)$ is a meromorphic function of $\lambda\in\af_{M,\C}^*$
of order $\le n+2$.
\end{prop}

\begin{proof} 
By (\ref{2.22}) we may assume that $\dim(\af_M/\af_G)=1$. Let
$P\in\cP(M)$ and let $\alpha$ be the unique simple root of $(P,A_M)$. Let 
$\pi\in\Pi_{\di}(M(\A))$. Then  $\cA^2_\pi(P)\not=\{0\}$ and
we have  to consider the intertwining operator 
$M_{\overline P|P}(\pi,\lambda)$. Recall that 
$M_{\overline P|P}(\pi,\lambda)$
is unitary for $\lambda\in i\af_{M}^*$. In particular, 
$M_{\overline P|P}(\pi,\lambda)$ is regular at $\lambda=0$. Put
$$M(\pi,\lambda)=M_{P|\overline P}(\pi,0)M_{\overline P|P}(\pi,\lambda),
\quad\lambda\in\af_{M,\C}^*.$$
Next consider the normalized intertwining operator $N_{\ov P|P}(\pi,\lambda)$
which is defined by (\ref{2.19}).
It follows from (\ref{2.8}), (\ref{2.9}), (\ref{2.16}) and (\ref{2.17}) 
that
$$N_{P|\overline P}(\pi,0)^*N_{P|\overline P}(\pi,0)=\Id.$$
Hence $N_{P|\overline P}(\pi,\lambda)$  is regular at $\lambda=0$ and 
$N_{P|\ov P}(\pi,0)$ is invertible. Put
$$N(\pi,\lambda)=N_{P|\overline P}(\pi,0)N_{\overline P|P}(\pi,\lambda),\;\quad
\lambda\in\af_{M,\C}^*.$$
Furthermore by (\ref{2.20}) and (\ref{2.21}) we get
\begin{equation}\label{4.1}
|r_{\overline P|P}(\pi,\lambda)|=1,\quad\lambda\in i\af_M^*.
\end{equation}
Thus $r_{\overline P|P}(\pi,\lambda)$ is also regular at $\lambda=0$ and
$r_{\overline P|P}(\pi,0)\not=0$. By
(\ref{2.19}) we get
\begin{equation}\label{4.2}
M(\pi,\lambda)=r_{P|\overline P}(\pi,0)r_{\overline P|P}(\pi,\lambda)
N(\pi,\lambda).
\end{equation}
Now observe that there 
exists an open compact subgroup $K_f\subset G(\A_f)$ such that 
$\cA^2_\pi(P)_{K_f}\not=\{0\}$.
 Hence there exists $\sigma
\in\Pi(K_\infty)$
such that $\cA^2_\pi(P)_{K_f,\sigma}\not=\{0\}$ 
(cf. section 1.8 for the definition).
Put
$$d=\dim \cA^2_\pi(P)_{K_f,\sigma}$$
and
$$c(\pi,\sigma)=r_{P|\overline P}(\pi,0)^{-d}.$$
Then $|c(\pi,\sigma)|=1$. Let $M(\pi,\lambda)_{K_f,\sigma}$  
(resp. $N(\pi,\lambda)_{K_f,\sigma}$) denote the restriction 
of $M(\pi,\lambda)$ (resp. $N(\pi,\lambda)$) to the subspace 
$\cA^2_\pi(P)_{K_f,\sigma}$. Then we have
$\det N(\pi,\lambda)_{K_f,\sigma}\not\equiv0$ and by  (\ref{4.2}) we get
\begin{equation}\label{4.3}
r_{\ov P|P}(\pi,\lambda)^d=c(\pi,\sigma)\,
\frac{\det M(\pi,\lambda)_{K_f,\sigma}}
{\det N(\pi,\lambda)_{K_f,\sigma}}.
\end{equation}
Thus it suffices to prove that both the numerator and the denominator on the
right hand side are meromorphic functions of order $\le n+2$. 
As for the numerator, it follows from Theorem 0.1 of \cite{Mu3} that
 $\det M(\pi,\lambda)_{K_f,\sigma}$ is a meromorphic function of 
$\lambda\in\af_{M,\C}^*$ of order $\le n+2$. In fact, in \cite{Mu3} we only
dealt with the case of the trivial character $\xi$. However, all the results
of \cite{Mu3} can be extend without any difficulty  to the case of a nontrivial
character $\xi$. It remains to consider the
denominator. By (\ref{2.11}), (\ref{2.16}) and (\ref{2.17}) there exists a
finite set $S_\pi$ of finite places of $\Q$ such that
\begin{equation}\label{4.4}
\begin{split}
\det N(\pi,\lambda)_{K_f,\sigma}&=
\det\left(R_{P|\overline P}(\pi_\infty,0)_{\sigma}
R_{\overline P|P}(\pi_\infty,\lambda)_{\sigma}\right)\\
&\quad\cdot\prod_{v\in S_\pi}\det\left(R_{P|\overline P}(\pi_v,0)_{K_v}
R_{\overline P|P}(\pi_v,\lambda)_{K_v}\right),
\end{split}
\end{equation}
where $R_{\overline P|P}(\pi_\infty,\lambda)_{\sigma}$
denotes the restriction of $R_{\overline P|P}(\pi_\infty,\lambda)$ to the 
$\sigma$-isotypical subspace $\H_P(\pi_\infty)_{\sigma}$ of
$\H_P(\pi_\infty)$ and 
$R_{\overline P|P}(\pi_v,\lambda)_{K_v}$ denotes the restriction of 
$R_{\overline P|P}(\pi_v,\lambda)$ to the subspace $\H_P(\pi_v)^{K_v}$ of
$K_v$-invariant functions. By Theorem 2.1 of \cite{A7}, 
$R_{\overline P|P}(\pi_\infty,\lambda)$ is a rational function of 
$\lambda(\alpha^\vee)$ and for each finite place $v$, 
$R_{\overline P|P}(\pi_v,\lambda)$ is a rational function of 
$q_v^{-\lambda(\alpha^\vee)}$. Therefore 
$\det(R_{P|\overline P}(\pi_\infty,0)_{\sigma}
R_{\overline P|P}(\pi_\infty,\lambda)_{\sigma})$  
 is a rational functions of $\lambda(\alpha^\vee)$ and for each $v<\infty$,
$\det(R_{P|\overline P}(\pi_v,0)_{K_v}
R_{\overline P|P}(\pi_v,\lambda)_{K_v})$ is a rational function of
$q_v^{-\lambda(\alpha^\vee)}$. Since  the function $z\in\C\mapsto q^{-z}$
is entire and of order 1, it follows that 
$\det N(\pi,\lambda)_{\sigma,K_f}$ is 
a meromorphic function of $\lambda  \in\af_{M,\C}^*$ of order $\le 1$. By
(\ref{4.3}) it follows that $r_{\overline P|P}(\pi,\lambda)^d$ and hence
$r_{\overline P|P}(\pi,\lambda)$ is a meromorphic function of 
$\lambda  \in\af_{M,\C}^*$ of order $\le n+2$.
\end{proof}

{\bf Remark.} 
Assume that $G$ is a quasi-split connected reductive group over
a number field $F$ with ring of ad\`eles $\A_F$. Let $P=MN$ be a maximal 
parabolic subgroup of $G$. Let $\pi$ be a globally generic cuspidal 
representation of $M(\A_F)$. Then it follows from \cite{Sh4} that the 
intertwining operator $M_{\ov P|P}(\pi,\lambda)$ can be normalized by
automorphic $L$-functions. Furthermore in \cite{GS}, Gelbart and Shahidi 
proved that the $L$-functions occurring in the normalizing factor are 
meromorphic functions of order 1. Therefore, one should 
expect that the normalizing factor $r_{Q|P}(\pi,\lambda)$ is of order 1
in general.

Now assume that $\dim \af_M/\af_G=1$. Our next goal is to 
 estimate the number of poles of $r_{\ov P|P}(\pi,\lambda)$ in a 
circle of radius $R>0$. For this purpose  we have to introduce some notation.

Let $\Pi_{\di}(M(\A);K_f)$ be the space of representations defined by 
(\ref{1.2}). For every $\pi\in\Pi_{\di}(M(\A);K_f)$ we have 
$\cA^2_\pi(P)^{K_f}\not=\{0\}$. 

Let $\pi$ be an irreducible unitary representation of $M(\R)$ and let
$I_P^G(\pi)$ be the induced representation of $G(\R)$. Recall that
among all $K_\infty$-types $\tau_{\Lambda^\prime}$ occurring in $I_P^G(\pi)$, 
the minimal
$K_\infty$-types of $I_P^G(\pi)$ are those $\tau_\Lambda$ for which 
$|\Lambda^\prime+2\rho_K|^2$ is minimizing at $\Lambda^\prime=\Lambda$. Let 
$W_P(\pi)$ be the
set of minimal $K_\infty$-types of $I_P^G(\pi)$. Then $W_P(\pi)$ is a non 
empty finite subset of $\Pi( K_\infty)$. Let $\lambda_\pi$  be the Casimir
eigenvalues of $\pi$ and for any $\tau\in
\Pi(K_\infty)$, let $\lambda_\tau$ be the Casimir eigenvalue of $\tau$. Put
\begin{equation}\label{4.5}
\Lambda_\pi:=\min_{\tau\in W_P(\pi)}\sqrt{\lambda_\pi^2+\lambda_\tau^2}.
\end{equation}

If $\pi\in\Pi(M(\A))$, put
$$\Lambda_\pi:=\Lambda_{\pi_\infty}.$$

For a given pole $\eta$ of $r_{\overline P|P}(\pi,\lambda)$, let $m(\eta)$
denote its order. Set
$$n_P(\pi,R)=\sum_{|\eta|\le R}m(\eta),$$
where the sum runs over all poles of $r_{\ov P|P}(\pi,\lambda)$.
\begin{prop}\label{p4.2} 
Let $m=\dim G$ and let $K_f$ be an open compact subgroup
of $G(\A_f)$.  There exists $C>0$ such that for all $R>0$ and
all $\pi\in\Pi_{\di}(M(\A);K_f)$ we have
$$n_P(\pi,R)\le C(1+R^2+\Lambda_\pi^2)^{8m}.$$
\end{prop} 

\begin{proof}
Let $\pi\in\Pi_{\di}(M(\A);K_f)$. Then there exists 
$\sigma\in\Pi(K_\infty)$ such that $\cA^2_\pi(P)_{K_f,\sigma}
\not=\{0\}$. Put
$$\overline N(\pi,\lambda):=N_{P|\overline P}(\pi,\lambda)
N_{\overline P|P}(\pi,0),\quad \lambda\in\af_{M,\C}^*.$$
Then
$$\overline N(\pi,\lambda)N(\pi,\lambda)=\Id$$
and by (\ref{4.3}) we get
\begin{equation}\label{4.6}
r_{\overline P|P}(\pi,\lambda)^d=c(\pi,\sigma)\,
\det (M(\pi,\lambda)_{K_f,\sigma})\cdot
\det (\overline N(\pi,\lambda)_{K_f,\sigma}).
\end{equation}
Thus it suffices to estimate the number of poles of the functions 
occurring on the right. It follows from Proposition 6.6 and
Lemma 6.1  of \cite{Mu3}, that the number  of
poles, counted with their order, of $\det M(\pi,\lambda)_{K_f,\sigma}$ 
in the disc $|\lambda|\le R$ is bounded by
$$C(1+R^2+\lambda_\pi^2+\lambda_\sigma^2)^{8m},$$
where $C>0$ is independent of $\pi$ and $\sigma$. As noted above, in 
\cite{Mu3} we only dealt with the case of the trivial character $\xi$. 
However, everything can be extended to a nontrivial character $\xi$ 
without any difficulty.

It remains to consider $\det \ov N(\pi,\lambda)_{K_f,\sigma}$. 
For any place $v$ let 
$$\ov R(\pi_v,\lambda)=R_{P|\ov P}(\pi_v,\lambda)R_{\ov P|P}(\pi_v,0).$$
By (\ref{2.16}) and (\ref{2.17}) we have
$$\ov N(\pi,\lambda)\circ j_P=j_P
\circ\left((\otimes_v\ov R(\pi_v,\lambda))\otimes\Id\right),$$
and there exists a finite set $S_\pi$
of finite places, which depends only on $\pi$ and $K_f$, such that
$$\ov R(\pi_v,\lambda)_{K_v}=\Id$$
for all $v\notin S_\pi\cup\{\infty\}$. Thus
\begin{equation}\label{4.7}
\det\left(\ov N(\pi,\lambda)_{K_f,\sigma}\right)=
\det \left(\ov R(\pi_\infty,\lambda)_{\sigma}\right)
\cdot\prod_{v\in S_\pi}\det\left( \ov R(\pi_v,\lambda)_{K_v}\right).
\end{equation}
Let $n_P(\pi_v,R)$, $v<\infty$, (resp. $n_P(\pi_\infty,R)$) denote the number 
of poles, counted with the order, of 
$\det\left(\ov R(\pi_v,\lambda)_{K_v}\right)$ (resp.
$\det\left(\ov R(\pi_\infty,\lambda)_{\sigma}\right)$) in the disc 
$|\lambda|< R$. Then we have to estimate $n_P(\pi_v,R)$ for
any $v\le\infty$.

Let $v<\infty$ and let $\pi_v$ be any irreducible unitary representation of
$M(\Q_v)$. Let $\widetilde\alpha\in\af_M$ be as  (\ref{3.5}).
By Theorem 2.2.2 of \cite{Sh1} there exists a polynomial 
$Q_v(z)$, $Q_v(0)=1$, such that
$$Q_v\left(q_v^{-\lambda(\widetilde\alpha)}\right)J_{P|\overline P}(\pi_v,\lambda)$$
is a holomorphic and non-zero operator. Moreover, the degree of the 
polynomial $Q_v$ is independent of $\pi_v$. Let
$$d_v=\dim\H_P(\pi_v)^{K_v}.$$
Then it follows from (\ref{2.7}) and the definition of $\ov R(\pi_v,\lambda)$
that
$$r_{P|\ov P}(\pi_v,\lambda)^{d_v}
Q_v\left(q_v^{-\lambda(\widetilde\alpha)}\right)^{d_v}
\det\left(\ov R(\pi_v,\lambda)_{K_v}\right)$$
is a holomorphic function on $\af_{M,\C}^*$.
By (r.5) there exist polynomials $P_1(z)$ and $P_2(z)$ such that  
$$r_{P|\ov P}(\pi_v,\lambda)=\frac{P_1(q_v^{-\lambda(\widetilde\alpha)})}
{P_2(q_v^{-\lambda(\widetilde\alpha)})}.$$
Thus it suffices to estimate the number of zeros of 
$P_1(q_v^{-\lambda(\widetilde\alpha)})$ and 
$Q_v(q_v^{-\lambda(\widetilde\alpha)})$,
respectively, in a circle of radius $R>0$. 
First observe that  for every $z\in\C$ the number of solutions of $q_v^{-s}=z$ 
in the disc $|s|\le R$ is bounded by $1+(2\pi)^{-1}\log(q_v)R$. Furthermore,
 the degree of the polynomial $Q_v$ is bounded by some constant $c_v>0$
which is independent of $\pi_v$. Using Lemma \ref{l3.1} and 
(2.1)--(2.3) of \cite{A7},  it follows that
the degree of the polynomial $P_1(z)$ is also bounded by a constant which
is independent of $\pi_v$. This implies that there exists $C_v>0$ such that
$$n_P(\pi_v,R)\le C_v\dim\left(\H_P(\pi_v)^{K_v}\right)(1+R)$$
for all $\pi_v\in\Pi(M(\Q_v))$ and $R>0$. It remains to estimate the
dimension of $\H_P(\pi_v)^{K_v}$. Suppose that $\pi_v$ is the component at 
$v$ of a representation $\pi\in\Pi_{\di}(M(\A));K_f)$. Then there
exists $\xi\in\Pi(A_M(\R)^0)$ such that $\pi\in\Pi_{\di}(M(\A))_\xi$. Let 
$\H_P(\pi)_{\sigma}^{K_f}$ be the $\sigma$-isotypical subspace of 
$\H_P(\pi)^{K_f}$. 
By (\ref{2.3}) it follows that
$$\H_P(\pi)_{\sigma}^{K_f}
\otimes\Hom_{M(\A)}(\pi,I_{M(\Q)AM(\R)^0}^{M(\A)}(\xi))
\cong\cA_\pi^2(P)_{K_f,\sigma}.$$
Moreover we have
$$\H_P(\pi)_ {\sigma}^{K_f}\cong\H_P(\pi_\infty)_{\sigma}
\otimes\bigotimes_{v<\infty}\H_P(\pi_v)^{K_v}$$
and $\dim\H_P(\pi_v)^{K_v}=1$ for $v\notin S_\pi$. Thus it follows that
$$\dim\H_P(\pi_v)^{K_v}\le\dim A_\pi^2(P)_{K_f,\sigma}.$$
The right hand side can be estimated  by Lemma 6.1 of \cite{Mu3}. 
It follows that

\begin{equation}\label{4.8}
n_P(\pi_v,R)\le C_v(1+\lambda_\pi^2+\lambda_{\sigma}^2)^{3m}(1+R).
\end{equation}

\smallskip
Now let $v=\infty$ and let $\pi_\infty\in\Pi(M(\R))$. Set
$$\ov J(\pi_\infty,\lambda)
=J_{P|\ov P}(\pi_\infty,\lambda)J_{\ov P|P}(\pi_\infty,0).$$
Then
\begin{equation}\label{4.9}
\ov R(\pi_\infty,\lambda)=\left(r_{P|\ov P}(\pi_\infty,\lambda)
r_{\ov P|P}(\pi_\infty,0)\right)^{-1}\ov J(\pi_\infty,\lambda).
\end{equation}
Let $K_{M,\infty}=K_\infty\cap M(\R)$ and let
$$\sigma|_{K_{M,\infty}}=\bigoplus_{\tau\in\Pi(K_{M,\infty})}
n_\tau\tau.$$
Set
$$[\sigma:\pi_\infty]=\sum_{\tau\in\Pi(K_{M,\infty})}n_\tau[\tau:
\pi_\infty|_{K_{M,\infty}}].$$
By Corollary 4.7 of \cite{VW}, there exist complex numbers $a_i(\pi_\infty)$,
 $i=1,...,r$, and $b_i(\pi_\infty,\sigma)$, $i=1,...,
r[\sigma:\pi_\infty]$, with $r=r(\pi_\infty)$ depending only on 
$\pi_\infty$,  and a constant $C\in\C$, such that

\begin{equation}\label{4.10}
\begin{split}
\det&\ov J(\pi_\infty,\lambda)\\
&=C\;\frac{\prod_{i=1}^{r(\pi_\infty)}
\Gamma\left(\langle\lambda,\alpha^\vee\rangle/
(4\langle\rho_P,\alpha^\vee\rangle)-
a_i(\pi_\infty)\right)^{[\sigma:\pi_\infty]}}
{\prod_{i=1}^{r(\pi_\infty)[\sigma:\pi_\infty]}
\Gamma\left(\langle\lambda,\alpha^\vee\rangle/
(4\langle\rho_P,\alpha^\vee\rangle)-b_i(\pi_\infty)\right)}.
\end{split}
\end{equation}
  
\begin{lem}\label{l4.3}
There exists $c>0$ such that $r(\tau)\le c$ for all $\tau\in\Pi(M(\R))$.
\end{lem}

\begin{proof} Let $b_\tau$ be the polynomial which is associated to $\tau$ 
by Theorem 1.5 of \cite{VW}. Then $r(\tau)$ is the degree of $b_\tau$
\cite[p.228]{VW}. So we have to estimate the degree of $b_\tau$. The 
polynomial $b_\tau$ is obtained from a more general 
polynomial $b_{\tau,\lambda}$ occurring in Theorem 2.2 of \cite{VW} by
choosing $\lambda=4\rho_P$. The polynomial $b_{\tau,\lambda}$ is 
associated to $\tau$ and a finite dimensional representation $(\eta,F)$
of $G$ satisfying the conditions (1)--(3) in \cite[p.210]{VW}. Then 
$\lambda$ is the action of $\af_M$ on $F^{\nf}$. It follows from the
constructions on pp. 217-219 in \cite{VW}, that the polynomial 
$b_{\tau,\lambda}$ is the product of the denominator of $\beta$ and the
denominator of the element $\check{z}_\nu$, defined on p.219. 
Let $\Omega\in Z(\gf_\C)$ be the Casimir element, and let 
$\chi_\Lambda$
be the infinitesimal character of $\tau$. Then it follows 
 that $b_{\tau,\lambda}$ equals
\begin{equation*}
\begin{split}
b_{\tau,\lambda}(\nu)=&\prod_{\mu\in\Pi(F)-\{\lambda\}}
\left(\chi_{\Lambda+\nu}(\Omega)
-\chi_{\Lambda+\nu+\lambda-\mu}(\Omega)\right)^{r(\mu)}\\
&\cdot\prod_{\mu\in\Pi(F)-\{\lambda\}}
\left(\chi_{\Lambda+\nu+\lambda}(\Omega)-\chi_{\Lambda+\nu+\mu}(\Omega)\right)^{r(\mu)}
\end{split}
\end{equation*}

Here $\Pi(F)$ denotes the set of weights of $F$ with respect to a fixed Cartan
subalgebra $\hf$ of $\gf$ and $r(\mu)$ is the multiplicity of a given weight 
$\mu$. From this description of $b_{\tau,\lambda}$ it follows that
$$r(\tau)\le 2(\dim F-1).$$
 Finally, 
$b_\tau$ is obtained by choosing $F$ to be the representation described
in example 2.1 of \cite{VW}.
\end{proof}

Now recall that the  poles of the Gamma function $\Gamma(z)$ are simple 
poles at $z=0,-1,-2,...$ and $1/\Gamma(z)$ is entire. Then it follows from 
(\ref{4.10}) together with Lemma \ref{l4.3} that there exists a constant 
$C_1>0$,
 independent of $\pi_\infty$, such that the number of zeros, counted with their
order, of $\det \ov J(\pi_\infty,\lambda)$ in the disc $|\lambda|\le R$ is 
bounded by
\begin{equation}\label{4.11}
C_1[\sigma:\pi_\infty](1+R).
\end{equation}

By Theorem 8.1 of \cite{Kn} and remark 1 following Theorem 8.4
in \cite{Kn}, we have 
$$[\sigma:\pi_\infty]\le\sum_{\tau\in\Pi(K_{M,\infty})}n_\tau\dim\tau\le
\dim\sigma.$$
Furthermore, by Weyl's dimension formula, there exists $C_2>0$ such
that $\dim\sigma\le C_2(1+\lambda_{\sigma}^2)^p$, where 
$p=1/2\dim K$.
Thus (\ref{4.11}) is bounded by
$$C_2(1+\lambda_{\sigma}^2)^p(1+R).$$

It remains to consider the normalizing factor 
$r_{P|\ov P}(\pi_\infty,\lambda)$. It is given by formula (\ref{3.7}). 
Let $\phi_\lambda:W_\R\to {}^LM$ be the map associated to 
$(\pi_{\infty})_\lambda$.
Let
$q$ be the number of irreducible constituents occurring in the 
decomposition (\ref{3.8}) of $\rho\cdot\phi_\lambda$.
Then $q$ is bounded independently of $\pi$ and it follows from the description 
of the $L$-factors in \S 3 that $r_{P|\ov P}(\pi_\infty,\lambda)$ is a 
product of $q$ meromorphic functions of the form (\ref{3.9}) or (\ref{3.10}).
From the form of these functions it follows immediately  that the number of 
poles, counted with their order, of $r_{P|\ov P}(\pi_\infty,\lambda)$ in the 
disc $|\lambda|\le R$ is bounded by $C(1+R)$. Putting our estimated together,
we have proved that there exists $C>0$, depending on $K_f$, such that
\begin{equation}\label{4.12}
n_P(\pi,R)\le C(1+R^2+\lambda_\pi^2+\lambda_{\sigma}^2)^{8m}
\end{equation}
for all $R\ge0$, and all $\pi\in\Pi_{\di}(M(\A))$ and $\sigma\in\Pi(K_\infty)$ 
such that $A_\pi^2(P)_{K_f,\sigma}\not=\{0\}$. 

Let $\pi\in\Pi_{\di}(M(\A);K_f)$. Let $\tau$ be a minimal 
$K_\infty$-type of $I_P^G(\pi_\infty)$.  Choose $\sigma\in\Pi(K_\infty)$ such
that $\sigma=\tau$.
Then (\ref{4.12}) applied to $\sigma$ together with
the definition of $\Lambda_\pi$ implies the proposition.
\end{proof}

\begin{cor}\label{c4.4} 
Let $m=\dim G$ and $n=16m+2$. There exists $C>0$, depending on $K_f$,
such that for each $\pi\in\Pi_{\di}(M(\A);K_f)$ we have
$$\sum_{\rho\not=0}\frac{m(\rho)}{|\rho|^n}\le C(1+\Lambda_\pi^2)^{8m},$$
where $\rho$ runs over the poles of $r_{\ov P|P}(\pi,\lambda)$. 
\end{cor}

\section[Logarithmic derivatives]{Logarithmic derivatives of global 
normalizing factors}
\setcounter{equation}{0}
In this section we shall study generalized  logarithmic derivatives of the 
global normalizing factors. First we assume that $M\in\cL(M_0)$
is such that $\dim\af_M/\dim\af_G=1$. Let $P\in\cP(M)$ and  let $\alpha$ be 
the unique simple root of $(P,A)$. Let $\pi\in\Pi_{\di}(M(\A))$ with 
$\pi=\otimes_v\pi_v$. By property (r.5) satisfied by the local normalizing
factors, it follows that for each place $v$, there exists a  meromorphic 
function 
$\tilde r_{\ov P|P}(\pi_v,z)$ of one complex variable $z$ such that the local
normalizing factor $r_{\ov P|P}(\pi_v,\lambda)$ is given by
$$r_{\ov P|P}(\pi_v,\lambda)=\tilde r_{\ov P|P}(\pi_v,\lambda(\alpha^\vee)).$$
Let 
$$\tilde r_{\ov P|P}(\pi,z):=\prod_v\tilde r_{\ov P|P}(\pi_v,z).$$
The infinite product is absolutely convergent in the half-plane 
$\Re(z)>\rho_P(\alpha^\vee)$,  admits a meromorphic continuation to 
$\C$ and the 
global normalizing factor is given by 
$$r_{\ov P|P}(\pi,\lambda)=\tilde r_{\ov P|P}(\pi,\lambda(\alpha^\vee)),
\quad\lambda\in\af_{M,\C}^*.$$ 
Our present goal is to estimate  the logarithmic derivative of 
$\tilde r_{\ov P|P}(\pi,z)$ along the imaginary axis.

To begin with, observe that by Lemma \ref{l3.2} and 
Lemma \ref{l3.3}  there exist $C,c>0$ such that
\begin{equation}\label{5.1}
\begin{split}
\bigg|\tilde r_{\ov P|P}&(\pi,z)^{-1}\frac{d}{dz}\tilde r_{\ov P|P}(\pi,z)
\bigg|\\
&\le
\sum_{v\le\infty}\bigg|\tilde r_{\ov P|P}(\pi_v,z)^{-1}
\frac{d}{dz}\tilde r_{\ov P|P}(\pi_v,z)\bigg|\le C
\end{split}
\end{equation}
for all $\pi\in\Pi_{\di}(M(\A))$ and $\Re(z)\ge c$. Using  (\ref{2.20}) and 
(\ref{2.21}), it follows that the function $\widetilde r_{\ov P|P}(\pi,z)$
satisfies
\begin{equation}\label{5.2}
\widetilde r_{\ov P|P}(\pi,z)\ov{\widetilde r_{\ov P|P}(\pi,-\ov{z})}=1.
\end{equation}
Hence  we get
$$\bigg|\tilde r_{\ov P|P}(\pi,z)^{-1}\frac{d}{dz}
\tilde r_{\ov P|P}(\pi,z)\bigg|
=\bigg|\tilde r_{\ov P|P}(\pi,-\ov z)^{-1}
\left(\frac{d\tilde r_{\ov P|P}}{dz}\right)(\pi,-\ov z)
\bigg|,$$
and together with (\ref{5.1}) we obtain the following proposition.
\begin{prop}\label{p5.1}
There exist $C,c>0$ such that 
$$\bigg|\tilde r_{\ov P|P}(\pi,z)^{-1}\frac{d}{dz}\tilde r_{\ov P|P}(\pi,z)\bigg|
\le C$$
for all $\pi\in\Pi_{\di}(M(\A))$ and  all $z\in\C$ with  $|\Re(z)|\ge c$.
\end{prop}

In order to get estimates for the logarithmic derivative on the imaginary
axis, we shall use the partial fraction decomposition of the meromorphic
function
${\widetilde r}_{\overline P|P}(\pi,z)^{-1}
(d/dz)({\widetilde r}_{\overline P|P}(\pi,z))$,
which allows us to treat the sum of the principal parts separately.
Let $n=16\dim G+2$. Then it follows from Corollary \ref{c4.4} that 
${\widetilde r}_{\overline P|P}(\pi,z)$ is a meromorphic function of order
$\le n$. Thus there exist entire functions $r_1(\pi,z)$ and $r_2(\pi,z)$ 
of order $\le n$ such that 
\begin{equation*}
{\widetilde r}_{\overline P|P}(\pi,z)=\frac{r_1(\pi,z)}{r_2(\pi,z)}.
\end{equation*}

Furthermore, observe that by (\ref{5.2}) a complex
number $\eta$ is a zero of ${\widetilde r}_{\overline P|P}(\pi,z)$ if and
only if $-\overline\eta$ is a pole of ${\widetilde r}_{\overline P|P}(\pi,z)$.
Thus by Hadamard's factorization theorem there exists a polynomial $Q(z)$ of 
degree $\le n$ such that
\begin{equation}\label{5.3}
{\widetilde r}_{\overline P|P}(\pi,z)=e^{Q(z)}
\frac{
 \displaystyle{\prod_\eta \left[ (1-\frac{z}{\eta}) 
\exp \left( \sum^n_{k=1}
\frac{1}{k} \left(\frac{z}{\eta}\right)^k \right)\right]^{a(\eta)}}
}
{\displaystyle{\prod_\eta\left[ (1+\frac{z}{\ov\eta}) \exp \left(\sum^n_{k=1}\frac{1}{k}
\left(-\frac{z}{\ov\eta}\right)^k\right)\right]^{a(\eta)}}},
\end{equation}

where $\eta$ runs over all zeros of ${\widetilde r}_{\overline P|P}(\pi,z)$ and $a(\eta)$ denotes the order of the zero $\eta$.

Let $D(\pi)$ denote the set of all poles and zeros of ${\widetilde r}_{\overline P|P}(\pi,z)$. Given $\eta\in D(\pi)$, we denote by $m(\eta)$ the order of
$\eta$, i.e., $m(\eta)$ is the integer such that 
$(z-\eta)^{-m(\eta)}$ $\widetilde r_{\overline P|P}(\pi,z)$ 
is holomorphic in a neighborhood of $\eta$ and does not vanish at
$z=\eta$. For $\eta\in\C^\times$ we define the function $h_\eta(z)$ by 
\begin{equation*}
h_\eta(z)=-\frac{1}{\eta}\sum^{n-1}_{k=0}\left(\frac{z}{\eta}\right)^k,
\quad z\in\C.
\end{equation*}

Then it follows from Corollary \ref{c4.4} that the series 
\begin{equation}\label{5.4}
f(\pi,z)=\sum_{\eta\in D(\pi)\atop|\eta|\le 1}\frac{m(\eta)}{z-\eta}+
\sum_{\eta\in D(\pi)\atop|\eta|>1}
m(\eta)
\left\{ \frac{1}{z-\eta}-h_\eta(z)\right\}
\end{equation}

is absolutely convergent on compact subsets of $\C\setminus D(\pi)$ and the 
resulting function $f(\pi,z)$ is a meromorphic function on $\C$ whose
set of poles equals $D(\pi)$. Differentiating (\ref{5.3}), we get
\begin{equation*}
{\widetilde r}_{\overline P|P}(\pi,z)^{-1}\frac{d}{dz}{\widetilde r}_{\overline P|P}(\pi,z)=
f(\pi,z)-\sum_{ \eta\in D(\pi)\atop|\eta|\le1}
m(\eta)h_\eta(z)+Q^\prime(z).
\end{equation*}

Thus there is a polynomial $g(\pi,z)$ of degree $\le n-1$ such that
\begin{equation}\label{5.5}
{\widetilde r}_{\overline P|P}(\pi,z)^{-1}\frac{d}{dz}{\widetilde r}_
{\overline P|P}(\pi,z)=f(\pi,z)+g(\pi,z).
\end{equation}

We begin with the investigation of $g(\pi,z)$. 

\begin{prop}\label{p5.2}
Let $m=\dim G$. There exist $C,c>0$ such that
\begin{equation*}
|g(\pi,z)|\le C(1+|z|^2+\Lambda^2_\pi)^{18m}
\end{equation*}

for all $\pi\in\Pi_{\di}(M(\A),K_f)$ and all $z\in\C$ with $|\Re(z)|\le c$. 
\end{prop}

\begin{proof}
Let $c>0$ be the constant occurring in Proposition \ref{p5.1}. First assume
that $|\Re(z)|=c$.
By Proposition \ref{p5.1} it suffices to estimate $f(\pi,z)$. Referring again
to Proposition \ref{p5.1}, it follows that $D(\pi)$ is contained in the strip
$|\Re(z)|\le c$. Hence  we may assume that $c>0$ has been chosen so that
for all $\pi\in\Pi_{\di}(M(\A),K_f)$, the zeros and poles of 
${\widetilde r}_{\overline P|P}(\pi,z)$ are contained in the strip
$|\Re(z)|<c-\delta$, where $\delta>0$ is independent of $\pi$. Hence the poles
of $f(\pi,z)$ are contained in $|\Re(z)|<c-\delta$. Let $\eta\in\C$ be a
pole of $f(\pi,z)$. Then for $|\Re(z)|\ge c$ we get
\begin{equation}\label{5.6}
|z-\eta|\ge|\Re(z-\eta)|\ge c-|\Re(\eta)|\ge\delta.
\end{equation}

Furthermore, from the definition of $f(\pi,z)$ it follows that
\begin{equation*}
\begin{split}
|f(\pi,z)|   \le&\sum_{|\eta|\le2|z|}\frac{|m(\eta)|}{|z-\eta|}  +
\sum_{ 1\le|\eta|\le2|z|} |m(\eta)|  |h_\eta(z)|\\
&  + \sum_{|\eta|>2|z|} |m(\eta)|
\Bigg|\frac{1}{z-\eta}-h_\eta(z)\Bigg|. 
\end{split}
\end{equation*}

Using (\ref{5.6}) and Proposition \ref{p4.2}, we can estimate the first sum 
as follows
\begin{equation*}
\begin{split}
\sum_{|\eta|\le2|z|} \frac{|m(\eta)|}{|z-\eta|}\le \frac{1}{\delta}
\sum_{|\eta|\le2|z|} |m(\eta)|
&\le \frac{2}{\delta}n_P(\pi,2|z|)\\
&\le C(1+|z|^2+\Lambda^2_\pi)^{8m}.
\end{split}
\end{equation*}

Again by Proposition \ref{p4.2}, we obtain for the second sum
\begin{equation}\label{5.7}
\begin{split}
\sum_{1\le|\eta|\le2|z|}|m(\eta)| |h_\eta(z)|&  \le 
\sum_{1\le|\eta|\le2|z|}|m(\eta)| \sum^{n-1}_{k=0}
\frac{|z|^k}{|\eta|^{k+1}}\\
&\le C|z|^{n-1}n_P(\pi,2|z|)\\
&   \le C_1 |z|^{n-1}(1+|z|^2+\Lambda^2_\pi)^{8m}.
\end{split}
\end{equation}

Finally, by Corollary \ref{c4.4} we get
\begin{equation}\label{5.8}
\begin{split}
\sum_{|\eta|>2|z|} |m(\eta)| \Bigg| \frac{1}{z-\eta}-h_\eta(z)\Bigg|
&  \le
2|z|^n \sum_{|\eta|>2|z|} \frac{|m(\eta)|}{|\eta|^{n+1}} \\
  & \le C(1+\Lambda^2_\pi)^{8m}|z|^n \\
 & \le C(1+|z|^2+\Lambda^2_\pi)^{8m}.
\end{split}
\end{equation}

Putting our estimates together, it follows that there exists $C>0$ such 
that
\begin{equation*}
|f(\pi,z)|\le C(1+|z|^2+\Lambda^2_\pi)^{18m}
\end{equation*}

for $|\Re(z)|\ge c$. Hence by Proposition \ref{p5.1}, there exists $C>0$ such
that
\begin{equation}\label{5.9}
|g(\pi,z)| \le C(1+|z|^2+\Lambda^2_\pi)^{18m}
\end{equation}

for all $\pi\in\Pi_{\di}(M(\A),K_f)$ and all $z\in\C$ with $|\Re(z)|=c$.
Now the proposition follows from the Phragmen-Lindel\"of theorem.
\end{proof}

Note that Proposition \ref{p5.2} gives  an upper bound for $g(\pi,z)$ on the 
imaginary axis. 

We shall now investigate $f(\pi,z)$. From the
definition of $f(\pi,z)$ by (\ref{5.3}) it is clear that the growth of
$f(\pi,z)$ along the imaginary axis depends on the distance of the poles
and zeros of $\tilde r_{\ov P|P}(\pi,z)$ from the imaginary axis. Therefore,
without any further information about the distribution of the poles and zeros
we cannot expect to get any estimates for $f(\pi,i\lambda)$ as $|\lambda|\to
\infty$. However, what  we can hope for is to obtain  estimates for integrals
involving $f(\pi,i\lambda)$.

To this end, we decompose  $f(\pi,z)$ as follows
\begin{equation*}
\begin{split}
f(\pi,z)=&\sum_{\eta\in D(\pi)\atop
|\eta|\le 2|z|}m(\eta)\frac{1}{z-\eta}
-\sum_{\eta\in D(\pi)\atop
1\le |\eta|\le 2|z|}m(\eta)h_\eta(z)\\
&+\sum_{\eta\in D(\pi)\atop|\eta|> 2|z|}m(\eta)\left\{\frac{1}{z-\eta}
-h_\eta(z)\right\}.
\end{split}
\end{equation*}
As for the second and the third sum, we observe that the estimations 
(\ref{5.7}) and (\ref{5.8}) are uniform in $z\in\C$. It remains to consider the first sum which we
denote by $f_1(\pi,z)$. Let
$$D_\pm(\pi)=\{\eta\in D(\pi)\mid\pm m(\eta)>0\}.$$
Then the map $\eta\to -\ov\eta$ is a bijection of 
$D_+(\pi)$ onto $D_-(\pi)$, and therefore, $f_1(\pi,z)$ can be written as
\begin{equation*}
f_1(\pi,z)=\sum_{\eta\in D_+(\pi)\atop
|\eta|\le 2|z|}m(\eta)\left\{\frac{1}{z-\eta}-
\frac{1}{z+\ov\eta}\right\}.
\end{equation*}

In particular, for $\lambda\in\R\setminus\{0\}$ we get
\begin{equation*}
f_1(\pi,i\lambda)=-\sum_{\eta\in D_+(\pi)\atop
|\eta|\le 2|\lambda|}m(\eta) 
\frac{2\Re(\eta)}{\Re(\eta)^2+(\lambda-\Im(\eta))^2}.
\end{equation*}

Let $\zeta\in C^\infty(\R)$ be such that $0\le\zeta\le1$, $\zeta(u)=0$ for 
$|u|\ge 3$ and
$\zeta(u)=1$ for $|u|\le 2$. Then it follows that
\begin{equation*}
\begin{split}
|f_1(\pi,i\lambda)|\le&\sum_{\eta\in D_+(\pi)\atop
|\eta|\le 2|\lambda|}m(\eta) 
\frac{2|\Re(\eta)|}{\Re(\eta)^2+(\lambda-\Im(\eta))^2}\\
&\le\sum_{\eta\in D_+(\pi)}\zeta\left(\frac{|\eta|}{|\lambda|}\right)m(\eta)
\frac{2|\Re(\eta)|}{\Re(\eta)^2+(\lambda-\Im(\eta))^2}.
\end{split}
\end{equation*}
Thus we have proved that for $\lambda\in\R\setminus\{0\}$ the following 
inequality holds
\begin{equation}\label{5.10}
\begin{split}
|f(\pi,i\lambda)|\le &\sum_{\eta\in D_+(\pi)}\zeta\left(\frac{|\eta|}{|\lambda|}\right)m(\eta)
\frac{2|\Re(\eta)|}{\Re(\eta)^2+(\lambda-\Im(\eta))^2}\\
&+C(1+\lambda^2+\Lambda_\pi^2)^{8m}.
\end{split}
\end{equation}
Put
\begin{equation*}
F(\lambda):=\begin{cases}
 \displaystyle{\sum_{\eta\in D_+(\pi)}\zeta\left(\frac{|\eta|}{|\lambda|}\right)m(\eta)
\frac{2|\Re(\eta)|}{\Re(\eta)^2+(\lambda-\Im(\eta))^2}} ,&\lambda\not=0;\\
0,&\lambda=0.
\end{cases}
\end{equation*}

Note that $0\notin D(\pi)$. Therefore,
on any finite interval $[-a,a]$, $F(\lambda)$ is the sum of finitely many
smooth and nonnegative functions. Hence $F(\lambda)$ is a smooth and
nonnegative function. We shall now estimate the integral of $F(u)$ over
a finite interval. Using Proposition \ref{p4.2} and the properties of
$\zeta$, we obtain

\begin{equation}\label{5.11}
\begin{split}
\int_0^\lambda F(u)\,du\le \sum_{\eta\in D_+(\pi)\atop
|\eta|\le 3\lambda}&m(\eta)\int_0^\lambda\frac{2|\Re(\eta)|}{\Re(\eta)^2+(u-\Im(\eta))^2}\,du\\
&\le 2\pi n_P(\pi,3\lambda)\\
&\le C(1+\lambda^2+\Lambda_\pi^2)^{8m}.
\end{split}
\end{equation}
Let $N\ge 8m+2$ and $R>0$. Using integration by parts and (\ref{5.11}), we 
obtain
\begin{equation*}
\begin{split}
\bigg|\int_{-R}^RF(u)(1+&u^2)^{-N}\,du\bigg|\\
&=\bigg|\int_{-R}^R\left(\int_0^u F(t)\,dt\right)
\frac{d}{du}(1+u^2)^{-N}\,du\bigg|\\
&\hskip20pt+\bigg|\int_{-R}^RF(t)\,dt\bigg|(1+R^2)^{-N}\\
&\le C(1+\Lambda_\pi^2)^{8m}\int_\R(1+u^2)^{8m}\bigg|\frac{d}{du}(1+u^2)^{-N}
\bigg|\,du\\
&\le C_N(1+\Lambda_\pi^2)^{8m}.
\end{split}
\end{equation*}
Here we have used that by (\ref{5.11}) the boundary term is bounded by a
constant independent of $R$.
Since $F\ge0$, this inequality implies that $F(u)$ is integrable with
respect to the measure $(1+u^2)^{-N}du$. 
Putting our estimates together, we obtain the following theorem.

\begin{theo}\label{th5.3}
Let $M\in\cL(M_0)$ and assume that $\dim \af_M/\af_G=1$. Let $P\in\cP(M)$ and
let $m=\dim G(\R)$. 
For every $N\ge 8m+2$ there exists $C_N>0$ such that for all 
$\pi\in\Pi_{\di}(M(\A);K_f)$ the following inequality holds
\begin{equation*}
\int_\R\bigg|\tilde r_{\ov P|P}(\pi,iu)^{-1}\frac{d}{du}\tilde r_{\ov P|P}(\pi,iu)\bigg|(1+u^2)^{-N}\,du\le C_N(1+\Lambda_\pi^2)^{18m}.
\end{equation*}
\end{theo}

Now suppose that  $M\in\cL(M_0)$ is arbitrary.  Then we 
have to consider the multidimensional logarithmic derivatives of the 
normalizing factors defined by Arthur in \cite{A4}. For this purpose we will 
use the notion of a $(G,M)$ family introduced by Arthur in Section 6 of 
\cite{A5}. For the convenience of the reader we  recall the  definition
of a $(G,M)$ family and explain some of its properties.

For each $P\in\cP(M)$, let $c_P(\lambda)$ be a smooth function on $i\af_M^*$. 
Then the set 
$$\{c_P(\lambda)\mid P\in\cP(M)\}$$
is called a $(G,M)$ family if the following holds: Let $P,P^\prime\in\cP(M)$
be adjacent parabolic groups and suppose that $\lambda$ belongs to the 
hyperplane spanned by the common wall of the chambers of $P$ and $P^\prime$.
Then
$$c_P(\lambda)=c_{P^\prime}(\lambda).$$
Let
\begin{equation}\label{5.12}
\theta_P(\lambda)=\vol\left(\af_P^G/\Z(\Delta_P^\vee)\right)^{-1}
\prod_{\alpha\in\Delta_P}\lambda(\alpha^\vee),\quad\lambda\in i\af_P^*,
\end{equation}
where $\Z(\Delta_P^\vee)$ is the lattice in $\af_P^G$ generated by the 
co-roots
$$\{\alpha^\vee\mid \alpha\in\Delta_P\}.$$
Let $\{c_P(\lambda)\}$ be a $(G,M)$ family. Then by Lemma 6.2 of \cite{A5},
the function
\begin{equation}\label{5.13}
c_M(\lambda)=\sum_{P\in\cP(M)}c_P(\lambda)\theta_P(\lambda)^{-1},
\end{equation}
which is defined on the complement of a finite set of hyperplanes,
extends to a smooth function on $i\af_M^*$. The value of $c_M(\lambda)$ at
$\lambda=0$ is of particular importance in connection with the spectral side
of the trace formula. It can be computed as follows. Let $p=\dim(A_M/A_G)$.
Set $\lambda=t\Lambda$, $t\in\R$, $\Lambda\in\af_M^*$, and let $t$ tend to 0.
Then
\begin{equation}\label{5.14}
c_M(0)=\frac{1}{p!}\sum_{P\in\cP(M)}\left(\lim_{t\to0}\left(\frac{d}{dt}\right)^p c_P(t\Lambda)\right)\theta_P(\Lambda)^{-1}
\end{equation}
\cite[(6.5)]{A5}. This expression is of course  independent of $\Lambda$.

For any $(G,M)$ family $\{c_P(\lambda)\mid P\in\cP(M)\}$ and any $L\in\cL(M)$
there is associated a natural $(G,L)$ family which is defined as follows.
 Let $Q\in\cP(L)$ and suppose that $P\subset Q$. 
The  function
$$\lambda\in i\af_L^*\mapsto c_P(\lambda)$$
depends only on $Q$. We will  denote it by $c_Q(\lambda)$. Then
$$\{c_Q(\lambda)\mid Q\in\cP(L)\}$$
is a $(G,L)$ family. We write
$$c_L(\lambda)=\sum_{Q\in\cP(L)}c_Q(\lambda)\theta_Q(\lambda)^{-1}$$
for the corresponding function (\ref{5.13}). 

Let $Q\in\cP(L)$ be fixed. If $R\in\cP^L(M)$, then 
$Q(R)$ is the unique group in $\cP(M)$ such that $Q(R)\subset Q$ and
$Q(R)\cap L=R$. Let $c_R^Q$ be the function on $i\af_M^*$ which is defined by
$$c_R^Q(\lambda)=c_{Q(R)}(\lambda).$$
Then $\{c_R^Q(\lambda)\mid R\in\cP^L(M)\}$ is an $(L,M)$ family. Let 
$c_M^Q(\lambda)$ be the function (\ref{5.13}) associated to this $(L,M)$ 
family.

We consider now special $(G,M)$ families defined by the global 
normalizing factors. Fix $P\in\cP(M)$, $\pi\in\Pi_{\di}(M(\A))$ and 
$\lambda\in i\af_M^*$. Define
\begin{equation}\label{5.15}
\nu_Q(P,\pi,\lambda,\Lambda):=r_{Q|P}(\pi,\lambda)^{-1}
r_{Q|P}(\pi,\lambda+\Lambda),\quad Q\in\cP(M).
\end{equation}
This set of functions is a $(G,M)$ family \cite[p.1317]{A4}. 
It is of a special form. Given $\beta\in\Sigma^r_P\cap\Sigma^r_{\ov Q}$, set
$$r_\beta(\pi,z)=\tilde r_{\ov P_\beta|P_\beta}(\pi,z),\quad z\in\C.$$
Then by (\ref{2.22}) we have
$$\nu_Q(P,\pi,\lambda,\Lambda)=\prod_{\beta\in\Sigma^r_Q\cap\Sigma^r_{\ov P}}
 r_\beta(\pi,\lambda(\beta^\vee))^{-1}
 r_\beta(\pi,\lambda(\beta^\vee)+\Lambda(\beta^\vee)).$$
Suppose that $L\in\cL(M)$, $L_1\in\cL(L)$ and $S\in\cP(L_1)$. Let 
$$\{\nu_{Q_1}^S(P,\pi,\lambda,\Lambda)\mid Q_1\in\cP^{L_1}(L)\}$$ 
be the associated $(L_1,L)$ family and let $\nu_L^S(P,\pi,\lambda,\Lambda)$
be the function (\ref{5.13}) defined by this family. Set
$$\nu_L^S(P,\pi,\lambda):=\nu_L^S(P,\pi,\lambda,0).$$
If $\beta$ is any root in $\Sigma^r(G,A_M)$, let $\beta_L^\vee$ denote the 
projection of $\beta^\vee$ onto $\af_L$. If $F$ is a subset of 
$\Sigma^r(G,A_M)$, let $F_L^\vee$ be the disjoint union of all the vectors
$\beta_L^\vee$, $\beta\in F$. 
Then by Proposition 7.5 of \cite{A4} we have
\begin{equation}\label{5.16}
\begin{split}
\nu_L^S(P,\pi,\lambda)=&\sum_F\vol(\af_L^{L_1}/\Z(F_L^\vee))\\
&\cdot\left(\prod_{\beta\in F} r_\beta(\pi,\lambda(\beta^\vee))^{-1}
r_\beta^\prime(\pi,\lambda(\beta^\vee))\right),
\end{split}
\end{equation}
where $F$ runs over all subsets of $\Sigma^r(L_1,A_M)$ such that $F_L^\vee$ 
is a basis of $\af_L^{L_1}$.
Let $N\in\N$. Then by (\ref{5.16}) we get
\begin{equation*}
\begin{split}
\int_{i\af_L^*/\af_G^*}|&\nu_L^S(P,\pi,\lambda)|
(1+\parallel\lambda\parallel^2)^{-N}\,d\lambda\le
\sum_F\vol(\af_L^{L_1}/\Z(F_L^\vee))\\
&\cdot\int_{i\af_L^*/\af_G^*}\prod_{\beta\in F}\Big| r_\beta(\pi,\lambda(\beta^\vee))^{-1}
r_\beta^\prime(\pi,\lambda(\beta^\vee))\Big|(1+\parallel\lambda\parallel^2)^{-N}\,d\lambda.
\end{split}
\end{equation*}
Here $F$ runs over all subsets of $\Sigma^r(L_1,A_M)$ such that $F_L^\vee$ is 
a basis of $\af_L^{L_1}$. Fix such a subset $F$. Let 
$$\{\tilde\omega_\beta\mid \beta\in F\}$$
be the basis of $(\af_L^{L_1})^*$ which is dual to $F_L^{L_1}$. We
 can write $\lambda\in i\af_L^*/i\af_G^*$ as
$$\lambda=\sum_{\beta\in F}z_\beta\tilde\omega_\beta+\lambda_1,\quad z_\beta\in
i\R,\quad\lambda_1\in i\af_{L_1}^*/i\af_G^*.$$
Observe that $\lambda(\beta^\vee)=z_\beta$. Suppose that $N>2\dim(A_{L_1}/A_G)+2$. Then there exists $C_N>0$, independent of $\pi$,  such that 
\begin{equation*}
\begin{split}
\int_{i\af_L^*/\af_G^*}\prod_{\beta\in F}&\Big| r_\beta(\pi,\lambda(\beta^\vee))^{-1}
r_\beta^\prime(\pi,\lambda(\beta^\vee))\Big|(1+\parallel\lambda\parallel^2)^{-N}\,d\lambda\\
&\le C_N\prod_{\beta\in F}\int_{i\R}\Big| r_\beta(\pi,z_\beta)^{-1}
r_\beta^\prime(\pi,z_\beta)\Big|(1+|z_\beta|^2)^{-N/2}\,dz_\beta.
\end{split}
\end{equation*}
Combined with Theorem \ref{th5.3} we obtain 
\begin{theo}\label{th5.4} 
Let $M\in\cL(M_0)$, $L\in\cL(M)$, $L_1\in\cL(L)$
and $S\in\cP(L_1)$. Let $m=\dim G(\R)/K_\infty$. For every $N\ge 8m+2$ there
exists $C_N>0$ such that
\begin{equation*}
\int_{i\af^*_L/\af^*_G}|\nu^S_L(P,\pi,\lambda)|(1+\parallel\lambda\parallel^2
)^{-N}d\lambda\le C_N(1+\lambda^2_\pi+\lambda^2_\sigma)^{8m^2}
\end{equation*}

for all $\pi\in\Pi_{\di}(M(\A),K_f)$ and any minimal $K_\infty$-type 
$\sigma$ of $I^G_P(\pi_\infty)$.
\end{theo}

\section{Absolute convergence of the spectral side}
\setcounter{equation}{0}

In this section we prove Theorem \ref{0.1} and Theorem \ref{th0.2}. For this purpose we have to study
 the multidimensional 
logarithmic derivatives of the global intertwining operators that are the main
ingredients of the spectral side. First we explain the structure of the 
spectral side in more detail. 
Let $M\in\cL(M_0)$. Fix $P\in\cP(M)$  and 
$\lambda\in i\af_M^*$. For $Q\in\cP(M)$ and $\Lambda\in i\af_M^*$ define
$$\mM_Q(P,\lambda,\Lambda)=M_{Q|P}(\lambda)^{-1}M_{Q|P}(\lambda+\Lambda).$$
Then
\begin{equation}\label{6.1}
\{\mM_Q(P,\lambda,\Lambda)\mid \Lambda\in i\af_M^*,\;Q\in\cP(M)\}
\end{equation}
is a $(G,M)$ family with values in the space of operators on $\cA^2(P)$
\cite[p.1310]{A4}. 

Let $L\in\cL(M)$. Then as above, the $(G,M)$ family (\ref{6.1}) has an 
associated $(G,L)$ family
$$\{\mM_{Q_1}(P,\lambda,\Lambda)\mid \Lambda\in i\af_L^*,\;Q_1\in\cP(L)\}
$$
and
$$\mM_L(P,\lambda,\Lambda)=\sum_{Q_1\in\cP(L)}\mM_{Q_1}(P,\lambda,\Lambda)
\theta_{Q_1}(\Lambda)^{-1}$$
extends to a smooth function on $i\af_L^*$. Put
$$\mM_L(P,\lambda)=\mM_L(P,\lambda,0).$$
For $s\in W(\af_M^*)$ let
$$M(P,s)=M_{P|P}(s,0).$$
The spectral side is a sum of distributions
$$\sum_{\chi\in\mX}J_\chi$$
on $G(\A)^1$. By Theorem 8.2 of \cite{A4}, the distribution $J_\chi$ can be
described as follows. 
Let $\chi\in \mX$, $\pi\in\Pi(M(\A)^1)$ and $h\in C^\infty_c(G(\A)^1)$. 
Note that $\mM_L(P,\lambda)$ and $\rho_{\chi,\pi}(P,\lambda,h)$ both act in 
the Hilbert space $\H_P(\pi)_\chi$. Let $W^L(\af_M)_{\reg}$ be the set of
elements $s\in W(\af_M)$ such that $\{H\in\af_M\mid sH=H\}=\af_L$.
 Then $J_\chi(f)$ equals the sum over $M\in\cL(M_0)$, $L\in\cL(M)$,
$\pi\in\Pi(M(\A)^1)$ and $s\in W^L(\af_M)_{\reg}$ of the product of 
$$|W_0^M|\;|W_0|^{-1}|\det(s-1)_{\af_M^L}|^{-1}$$
with
$$
\int_{i\af_L^*/i\af_G^*}|\cP(M)|^{-1}\sum_{P\in\cP(M)}\tr\left(\mM_L(P,\lambda)M(P,s)\rho_{\chi,\pi}(P,\lambda,h)\right)d\lambda.$$
Our goal is to determine the conditions under which  the 
integral-series obtained by summing this expression  over $\chi\in\mX$, is
absolutely convergent. 
 Since $M(P,s)$ is unitary, 
we have to estimate the 
integral
\begin{equation}\label{6.2}
\int_{i\af_L^*/i\af_G^*}\parallel\mM_L(P,\lambda)\rho_{\chi,\pi}(P,\lambda,h)\parallel_1d\lambda,
\end{equation}
where $\parallel\cdot\parallel_1$ denotes the trace norm. 

We shall now assume that $h\in\Co^1(G(\A)^1)$. 
Let $N_{Q|P}(\pi,\lambda)$, $P,Q\in\cP(M)$, be the normalized intertwining 
operator which by (\ref{2.19}) is defined as
$$N_{Q|P}(\pi,\lambda):=r_{Q|P}(\pi,\lambda)^{-1}M_{Q|P}(\pi,\lambda),\quad
\lambda\in\af_{M,\C}^*,$$
Let $P\in\cP(M)$ and $\lambda\in i\af_M^*$ be fixed.
For $Q\in\cP(M)$ and  $\Lambda\in i\af_M^*$ define
\begin{equation}\label{6.3}
\mN_Q(P,\pi,\lambda,\Lambda)=N_{Q|P}(\pi,\lambda)^{-1}N_{Q|P}(\pi,\lambda+\Lambda),
\end{equation}
Then as functions of $\Lambda\in i\af_M^*$, 
$$\{\mN_Q(P,\pi,\lambda,\Lambda)\mid Q\in\cP(M)\}$$
is  a $(G,M)$ family. The verification is the same as in the case of
the unnormalized intertwining operator \cite[p.1310]{A4}. For $L\in\cL(M)$, 
let
$$\{\mN_{Q_1}(P,\pi,\lambda,\Lambda)\mid \Lambda\in i\af_L^*,\;Q_1\in\cP(L)\}
$$
be the associated $(G,L)$ family. 

Let $\mM_{Q_1}(P,\pi,\lambda,\Lambda)$ be the restriction of 
$\mM_{Q_1}(P,\lambda,\Lambda)$ to $\H_P(\pi)_\chi$. 
Then by (\ref{2.19}) and (\ref{5.15}) it 
follows that
\begin{equation}\label{6.4}
\mM_{Q_1}(P,\pi,\lambda,\Lambda)=\mN_{Q_1}(P,\pi,\lambda,\Lambda)
\nu_{Q_1}(P,\pi,\lambda,\Lambda)
\end{equation}
for all $\Lambda\in i\af_L^*$ and all $Q_1\in\cP(L)$. 

For $Q\supset P$ let
$\hat L_P^Q\subset \af_P^Q$ be the lattice generated by 
$\{\tilde\omega^\vee\mid \tilde\omega\in\hat\Delta_P^Q\}$. Define
$$\hat\theta_P^Q(\lambda)=\vol(\af_P^Q/\hat L_P^Q)^{-1}
\prod_{\tilde\omega\in\hat\Delta_P^Q}\lambda(\tilde\omega^\vee).$$
For $S\in\cF(L)$ put
\begin{equation}\label{6.5}
\begin{split}
\mN_S^\prime&(P,\pi,\lambda)\\
&=\lim_{\Lambda\to0}\sum_{\{R\mid R\supset S\}}(-1)^{\dim(A_S/A_R)}
\hat\theta_S^R(\Lambda)^{-1}\mN_R(P,\pi,\lambda,\Lambda)\theta_R(\Lambda)^{-1}.
\end{split}
\end{equation}
Let $\mM_L(P,\pi,\lambda)$ be the restriction of 
$\mM_L(P,\lambda)$ to $\H_P(\pi)_\chi$. 
Then by Lemma 6.3 of \cite{A5} we have
\begin{equation}\label{6.6}
\mM_L(P,\pi,\lambda)=\sum_{S\in\cF(L)}\mN_S^\prime(P,\pi,\lambda)
\nu_L^S(P,\pi,\lambda).
\end{equation}
Hence the integral (\ref{6.2}) can be estimated by
$$
\sum_{S\in\cF(L)}\int_{i\af_L^*/i\af_G^*}\parallel \mN_S^\prime(P,\pi,\lambda)
\rho_{\chi,\pi}(P,\lambda,h)\parallel_1|\nu_L^S(P,\pi,\lambda)|\;d\lambda.$$

We shall now study the integral in more detail. Let $\Omega$ and $\Omega_K$ be
 the Casimir operators of $G(\R)$ and $K_\infty$ respectively. Set
$$\Delta=\Id-\Omega+2\Omega_K.$$
Then $\Delta$ acts on $\cA^2_{\chi,\pi}(P)$ through each of the representations
$\rho_{\chi,\pi}(P,\lambda)$. Let $K_f$ be an open compact  subgroup of 
$G(\A_f)$ and let $\sigma\in\Pi(K_\infty)$. Then the operators
$$\rho_{\chi,\pi}(P,\lambda,\Delta),\quad \lambda\in i\af_P^*,$$
have $\cA^2_{\chi,\pi}(P)_{K_f}$ and $\cA^2_{\chi,\pi}(P)_{K_f,\sigma}$ as 
invariant 
subspaces. We shall denote the restriction of 
$\rho_{\chi,\pi}(P,\lambda,\Delta)$ to $\cA^2_{\chi,\pi}(P)_{K_f}$ and 
$\cA^2_{\chi,\pi}(P)_{K_f,\sigma}$, respectively,
by $\rho_{\chi,\pi}(P,\lambda,\Delta)_{K_f}$ and
$\rho_{\chi,\pi}(P,\lambda,\Delta)_{K_f,\sigma}$, respectively. Recall that 
by (\ref{2.3}),
$\rho_{\chi,\pi}(P,\lambda)$ is equivalent to $I^G_P(\pi_\lambda)\otimes\Id$.
Let $\lambda_\pi$ and $\lambda_\sigma$  denote the Casimir eigenvalues of 
$\pi_\infty$ and $\sigma$, respectively. Then it follows from 
Proposition 8.22 of \cite{Kn} that
\begin{equation}\label{6.7}
\rho_{\chi,\pi}(P,\lambda,\Delta)_{K_f,\sigma}=
\left(1+\parallel\lambda\parallel^2-\lambda_\pi+2\lambda_\sigma\right)\Id.
\end{equation}

To estimate the right hand side we use the following lemma.

\begin{lem}\label{l6.1}
For all $\pi\in\Pi_{\di}(M(\A)^1;K_f)$ and $\sigma\in\Pi(K_\infty)$, 
one has
\begin{equation}\label{6.8}
-\lambda_\pi+\lambda_\sigma\ge0, \quad\mbox{if}\quad
\dim\cA^2_{\pi}(P)_{K_f,\sigma}\not=\{0\}.
\end{equation}
\end{lem}

\begin{proof}
The lemma is a consequence of a more general result. Let 
$\pi_\infty\in \Pi(M(\R))$ and suppose that $\sigma\in\Pi(K_\infty)$ 
occurs in $I_{P(\R)}^{G(\R)}(\pi_\infty)|_{K_\infty}$. Let 
$$\pi_\infty|_{K_\infty\cap M(\R)}=\sum_{\omega\in\Pi(K_\infty\cap M(\R))}
n_\omega \omega.$$
Then 
$$[I_{P(\R)}^{G(\R)}(\pi_\infty)|_{K_\infty}:\sigma]=\sum_{\omega\in
\Pi(K_\infty\cap M(\R))}n_\omega[\sigma|_{K_\infty\cap M(\R)}:\omega]$$
\cite[p.208]{Kn}. Hence there exists $\omega\in\Pi(K_\infty\cap M(\R))$ such
that
$$[\sigma|_{K_\infty\cap M(\R)}:\omega]>0\quad\mbox{and}\quad
[\pi_\infty|_{K_\infty\cap M(\R)}:\omega]>0.$$
By \cite[(5.15)]{Mu2}, the first condition implies that the Casimir 
eigenvalues $\lambda_\omega$ of $\omega$ and $\lambda_{\sigma}$ of
$\sigma$ satisfy $\lambda_\omega\le\lambda_{\sigma_\infty}$. On the
other hand, since $\omega$ occurs in $\pi_\infty|_{K_\infty\cap M(\R)}$ it 
follows that $-\lambda_{\pi_\infty}+\lambda_\omega\ge 0$ 
\cite[Lemma 2.6]{DH}. This completes the proof.
\end{proof}

Using (\ref{6.7}) and (\ref{6.8}), it follows that 
\begin{equation}\label{6.9}
\begin{split}
\parallel\rho_{\chi,\pi}(P,\lambda,\Delta)_{K_f,\sigma}\parallel^2&\ge
(1+\parallel\lambda\parallel)^2+(-\lambda_\pi+2\lambda_\sigma)^2\\
&\ge\frac{1}{4}(1+\parallel\lambda\parallel^2+\lambda_\pi^2+\lambda_\sigma^2).
\end{split}
\end{equation}
Let $S\in\cF(L)$ be fixed.
Given an open compact subgroup $K_f$ of $G(\A_f)$ and 
$\sigma\in\Pi(K_\infty)$, let 
$\mN_S^\prime(P,\pi,\lambda)_{K_f,\sigma}$ 
denote the
restriction of $\mN_S^\prime(P,\pi,\lambda)$ to $\cA^2_\pi(P)_{K_f,\sigma}$. 

\begin{lem}\label{l6.2}
Let $K_f$ be an open compact subgroup of $G(\A_f)$ and let 
$h\in \Co^1(G(\A)^1)$ be bi-invariant under $K_f$. 
Suppose that there exist $N\in\N$ and $C>0$ such that 
\begin{equation}\label{6.10}
\parallel\mN_S^\prime(P,\pi,\lambda)_{K_f,\sigma}\parallel\le 
C(1+\parallel\lambda\parallel^2+\lambda_\pi^2+\lambda_\sigma^2)^N
\end{equation}
for all $\pi\in\Pi_{\di}(M(\A),K_f)$, $\sigma\in\Pi(K_\infty)$ and 
$\lambda\in i\af_L^*$. Then for every $k\in\N$ there exists $C_k>0$ such that
$$\int_{i\af_L^*/i\af_G^*}\parallel \mN_S^\prime(P,\pi,\lambda)
\rho_{\chi,\pi}(P,\lambda,h)\parallel_1|\nu_L^S(P,\pi,\lambda)|\;d\lambda
\le C_k(1+\Lambda_\pi)^{-k}$$
for all $\chi\in\mX$ and $\pi\in\Pi(M(\A)^1)$.
\end{lem}

\begin{proof}
Since $h$ is bi-invariant under $K_f$, $\rho_{\chi,\pi}(P,\lambda,h)$ maps 
the Hilbert space $\ov \cA^2_{\chi,\pi}(P)$ into the subspace 
$\ov \cA^2_{\chi,\pi}(P)_{K_f}$. Moreover 
$\ov \cA^2_{\chi,\pi}(P)_{K_f}$ 
is an invariant subspace for $\rho_{\chi,\pi}(P,\lambda,h)$. Hence 
$\rho_{\chi,\pi}(P,\lambda,h)=0$, unless $\pi$ belongs to 
$\Pi_{\di}(M(\A),K_f)$. So we may assume that $\pi$ belongs to
$\Pi_{\di}(M(\A),K_f)$. 
Then for each $k\in\N$ we get

\begin{equation}\label{6.11}
\begin{split}
\parallel \mN_S^\prime(P,&\pi,\lambda)
\rho_{\chi,\pi}(P,\lambda,h)\parallel_1\\
&=\parallel \mN_S^\prime(P,\pi,\lambda)_{K_f}
\rho_{\chi,\pi}(P,\lambda,h)_{K_f}\parallel_1\\
&\le\parallel \mN_S^\prime(P,\pi,\lambda)_{K_f}
\rho_{\chi,\pi}(P,\lambda,\Delta^{2k})_{K_f}^{-1}\parallel_1\\
&\hskip82pt\cdot\parallel\rho_{\chi,\pi}(P,\lambda,\Delta^{2k}h)\parallel.
\end{split}
\end{equation}
Furthermore, using (\ref{6.9}) and (\ref{6.10}) we get
\begin{equation}\label{6.12}
\begin{split}
\parallel \mN_S^\prime&(P,\pi,\lambda)_{K_f}
\rho_{\chi,\pi}(P,\lambda,\Delta^{2k})_{K_f}^{-1}\parallel_1\\
&\le\sum_{\sigma\in\Pi(K_\infty)}
\parallel\mN_S^\prime(P,\pi,\lambda)_{K_f,\sigma}\parallel\cdot
\parallel\rho_{\chi,\pi}(P,\lambda,\Delta^{2k})_{K_f,\sigma}^{-1}\parallel_1\\
&\le 2C\sum_{\sigma\in\Pi(K_\infty)}\frac{\dim\cA^2_{\chi,\pi}(P)_{K_f,\sigma}}
{(1+\parallel\lambda\parallel^2+\lambda_\pi^2+\lambda_\sigma^2)^{k-N}}.
\end{split}
\end{equation}
By Lemma 6.1 of \cite{Mu3} there exist $C_1>0$ and $N_1\in\N$ such that
$$\dim\cA^2_\chi(P)_{K_f,\sigma}\le 
C_1(1+\lambda_\chi^2+\lambda_\sigma^2)^{N_1}$$
for all $\chi\in\mX$ and $\sigma\in\Pi(K_\infty)$. Actually in \cite{Mu3} we 
considered the space $\cA^2(P,\chi,\sigma)$, where $\sigma$ is an irreducible
 representation of $K$. The two spaces are not equal, but they are 
closely related. Moreover 
$\lambda_\chi$ was denoted by $\mu_\chi$ in \cite{Mu3}.
If $\cA^2_{\chi,\pi}(P)\not=0$,
it follows from Langlands' construction of $\cA^2_\chi(P)$ in terms of 
iterated residues of cuspidal Eisenstein series that 
\begin{equation}\label{6.13}
|\lambda_\chi-\lambda_\pi|\le c
\end{equation}
with $c>0$ independent of $\chi$ and $\pi$ (see (4.21) of \cite{Mu3}). Hence
there exist $C_2>0$ and $N_1\in\N$ such that
\begin{equation}\label{6.14}
\dim\cA^2_{\chi,\pi}(P)_{K_f,\sigma}
\le C_2(1+\lambda_\pi^2+\lambda_\sigma^2)^{N_1}
\end{equation}
for all $\chi\in\mX$, $\pi\in\Pi_{\di}(M(\A),K_f)$ and 
$\sigma\in\Pi(K_\infty)$. 
Set 
$$N_2=\frac{1}{2}(N+N_1).$$
Now observe that there exists $n_0\in\N$ such that 
$$\sum_{\sigma_\infty\in\Pi(K_\infty)}(1+\lambda_{\sigma})^{-n}<\infty.$$
the right hand side is finite for $n\ge n_0$.
Let $\Lambda_\pi$ be the number defined by (\ref{4.5}). 
Then by (\ref{6.12}) and (\ref{6.14}) it follows that for every 
$k>2(n_0+N_2)$ there exists $C_k>0$ such that 
\begin{equation}\label{6.15}
\begin{split}
\parallel \mN_S^\prime(P,\pi,\lambda)_{K_f}&
\rho_{\chi,\pi}(P,\lambda,\Delta^{4k})_{K_f}^{-1}\parallel_1\\
&\le C_k(1+\parallel\lambda\parallel^2+\Lambda_\pi^2)^{N_2-k}\\
&\le C_k(1+\parallel\lambda\parallel^2)^{(N_2-k)/2}
(1+\Lambda_\pi^2)^{(N_2-k)/2}
\end{split}
\end{equation}
for all $\pi\in\Pi_{\di}(M(\A),K_f)$ and $\chi\in\mX$. Next observe that
for $\lambda\in i\af^*_L$ the operator $\rho_{\chi,\pi}(P,\lambda,g)$ is
unitary. Hence it follows that
\begin{equation}\label{6.16}
\parallel\rho_{\chi,\pi}(P,\lambda,\Delta^{4k}h)\parallel\le
\parallel \Delta^{4k}h\parallel_{L^1(G(\A)^1)}
\end{equation}
for all $\pi\in\Pi_{\di}(M(\A),K_f)$ and $\chi\in\mX$. 
Combing (\ref{6.11}), (\ref{6.15}) and (\ref{6.16}), it follows that for
every $n\in\N$ there exists $C_n>0$ such that
$$\parallel \mN_S^\prime(P,\pi,\lambda)\rho_{\chi,\pi}(P,\lambda,h)\parallel_1
\le C_n(1+\parallel\lambda\parallel^2)^{-n}(1+\Lambda_\pi^2)^{-n}$$
for all $\chi\in\mX$ and $\pi\in\Pi_{\di}(M(\A),K_f)$.
Combined with Theorem \ref{th5.4} the claimed estimation of the integral 
follows.
\end{proof}

\smallskip
\noindent
{\it Proof of Theorem \ref{th0.1}:}
Let $h\in\Co^1(G(\A)^1)$ be
 bi-invariant under $K_f$. As observed in the proof of Lemma \ref{l6.2}, 
it follows that
$\rho_{\chi,\pi}(P,\lambda,h)=0$, unless $\pi\in\Pi_{\di}(M(\A),K_f)$. Let
$L_{\di}^2(M(\Q)A_P(\R)^0\ba M(\A))$ be the largest closed subspace of the
Hilbert space 
$L^2(M(\Q)A_P(\R)^0\ba M(\A))$ which decomposes discretely under the regular
representation of $M(\A)$. Then 
$$L_{\di}^2(M(\Q)A_P(\R)^0\ba M(\A))=\bigoplus_{\pi\in\Pi(M(\A))}
m(\pi)\H_\pi,$$
and each multiplicity $m(\pi)$ is finite.
Thus, if the assumption  (\ref{6.10}) of Lemma \ref{l6.2} is satisfied,  it
 follows from Lemma \ref{l6.2} that for every 
$n\in\N$ there exists $C_n>0$ such that 
\begin{equation}\label{6.17}
\begin{split}
\sum_{\chi\in\mX}\sum_{\pi\in\Pi(M(\A)^1)}
\int_{i\af_L^*/i\af_G^*}\parallel \mN_S^\prime(&P,\pi,\lambda)
\rho_{\chi,\pi}(P,\lambda,h)\parallel_1|\nu_L^S(P,\pi,\lambda)|\;d\lambda\\
&\le C_n\sum_{\pi\in\Pi_{\di}(M(\A)^1,K_f)}m(\pi)(1+\Lambda_\pi)^{-n}.
\end{split}
\end{equation}
It remains to investigate the sum on the right hand side. Let 
$$K_{M,f}=K_f\cap M(\A_f).$$
 Then there exist arithmetic subgroups $\Gamma_{M,i}\subset
M(\R)$, $i=1,...,l$, such that 
$$M(\Q)\backslash M(\A)/K_{M,f}\cong\bigcup_{i=1}^l(\Gamma_{M,i}\backslash
M(\R))$$
(cf. Section 9 of \cite{Mu1}). Therefore we get
\begin{equation}\label{6.18}
\begin{split}
L^2(A_M(\R)^0M(\Q)&\backslash M(\A))^{K_{M,f}}\\
&\cong\bigoplus_{i=1}^l L^2\bigl(A_M(\R)^0\Gamma_{M,i}\backslash M(\R)\bigr)
\end{split}
\end{equation}
as $M(\R)$-modules. For each $i$, $i=1,...,l$, let 
$L^2_{\di}(A_P(\R)^0\Gamma_{M,i}\backslash M(\R))$ be the discrete subspace of
the regular representation of $M(\R)$ in 
$L^2(A_P(\R)^0\Gamma_{M,i}\backslash M(\R))$. Then it follows from (\ref{6.18})
that
\begin{equation}\label{6.19}
\begin{split}
L^2_{\di}(A_P(\R)^0M(\Q)&\backslash M(\A))^{K_{M,f}}\\
&\cong\bigoplus_{i=1}^l L^2_{\di}(A_P(\R)^0\Gamma_{M,i}\backslash M(\R))
\end{split}
\end{equation}
as $M(\R)$ modules. For $i$, $1\le i\le l$, and $\pi_\infty\in\Pi(M(\R))$ 
denote by $m_{\Gamma_{M,i}}(\pi_\infty)$ the multiplicity of $\pi_\infty$ 
in the regular representation of $M(\R)$ in $L^2_{\di}(A_P(\R)^0\Gamma_{M,i}
\backslash M(\R))$. Then by (\ref{6.19}) we get
\begin{equation}\label{6.20}
\begin{split}
&\sum_{\pi\in\Pi_{\di}(M(\A)^1,K_f)}m(\pi)(1+\Lambda_\pi)^{-n}\\
&\hskip38pt\le\sum_{i=1}^l\sum_{\pi_\infty\in\Pi(M(\R))}m_{\Gamma_{M,i}}(\pi_\infty)
(1+\Lambda_{\pi_\infty})^{-n}
\end{split}
\end{equation}
Let $\sigma\in\Pi(K_\infty)$ be a minimal $K_\infty$-type occurring in
$I^G_P(\pi_\infty)$ with Casimir eigenvalue $\lambda_\sigma$.
 Let $K_{M,\infty}=
M(\R)\cap K_\infty$. By (5.15) of \cite{Mu2} we have that
$\lambda_\sigma\ge\lambda_\tau$ for any irreducible constituent 
$\tau\in\Pi(K_{M,\infty})$ of $\sigma_\infty|K_{M,\infty}$. Thus the right
hand side of (\ref{6.20}) is bounded by
$$\sum_{i=1}^l\sum_{\tau\in\Pi(K_{M,\infty})}\sum_{\pi_\infty\in\Pi(M(\R))}
m_{\Gamma_{M,i}}(\pi_\infty)
\frac{\dim(\H(\pi_\infty)\otimes V_\tau)^{K_{M,\infty}}}
{(1+\lambda_{\pi_\infty}^2+\lambda_\tau^2)^{n/2}}.$$
By Corollary 0.3 of \cite{Mu2} this sum is finite for $n$ sufficiently large.
Thus we proved
\begin{prop}\label{p6.3}
Let $K_f$ be an open compact subgroup of $G(\A_f)$ and let 
$h\in \Co^1(G(\A)^1)$ be bi-invariant under $K_f$. 
Suppose that there exist $N\in\N$ and $C>0$ such that 
\begin{equation}\label{6.21}
\parallel\mN_S^\prime(P,\pi,\lambda)_{K_f,\sigma}\parallel\le 
C(1+\parallel\lambda\parallel^2+\lambda_\pi^2+\lambda_\sigma^2)^N
\end{equation}
for all $\pi\in\Pi_{\di}(M(\A),K_f)$, $\sigma\in\Pi(K_\infty)$ and 
$\lambda\in i\af_L^*$. Then 
\begin{equation}\label{6.22}
\sum_{\chi\in\mX}\sum_{\pi\in\Pi(M(\A)^1)}\int_{i\af_L^*/i\af_G^*}\parallel\mM_L(P,\lambda)\rho_{\chi,\pi}(P,\lambda,h)\parallel_1d\lambda<\infty
\end{equation}
\end{prop}

Let $h\in \Co^1(G(\A)^1)$. Then there exists an open compact subgroup
$K_f$ of $G(\A_f)$ such that $h$ is bi-invariant under $K_f$. Using the 
 observations made at the beginning of this section, it follows that
 (\ref{6.22}) implies that the spectral side of the trace
formula is absolutely convergent.

We shall now continue by investigating  condition (\ref{6.21}) in detail. To
 calculate 
$\mN_S^\prime(P,\pi,\lambda)$, let
$\Lambda\in i\af_M^*$. By \cite[p.37]{A5} $\mN_S^\prime(P,\pi,\lambda)$ equals
\begin{equation*}
\begin{split}
\frac{1}{q!}\sum_{\{R\mid R\supset S\}}
(-1)^q\hat\theta_S^R(\Lambda)^{-1}\left(\lim_{t\to0}
\left(\frac{d}{dt}\right)^q\mN_R(P,\pi,\lambda,t\Lambda)\right)\theta_R(\Lambda)^{-1},
\end{split}
\end{equation*}
where $q=\dim(A_S/A_R)$. Since $N_{Q|P}(\pi,\lambda)$ is unitary for 
$\lambda\in i\af_M^*$, it follows from  (\ref{6.3}) that we have to estimate
the norm of
\begin{equation}\label{6.23}
\lim_{t\to0}\left(
\frac{d}{dt}\right)^qN_{Q|P}(\pi,\lambda+t\Lambda)_\sigma,
\quad\lambda\in i\af_M^*.
\end{equation}
To this end, we may use  (\ref{2.16}) and (\ref{2.17}) to
replace $N_{Q|P}(\pi,\lambda)$ by 
$$R_{Q|P}(\pi,\lambda)=\otimes_vR_{Q|P}(\pi_v,\lambda).$$
 Next note that any compact open subgroup 
$K_f=\prod_{v<\infty}K_v$
 of $G(\A_f)$ is such that $K_v$ is a hyperspecial compact subgroup for 
almost all $v$. 
Hence, by (\ref{2.11}) there exists a finite set of places $S_0$, including
 the 
Archimedean one, such that    we have
$$R_{Q|P}(\pi_v,\lambda)_{K_v}=\Id,\quad v\notin S_0,\,
\pi\in\Pi_{\di}(M(\A),K_f).$$
Let $D_\Lambda$ denote the directional derivative on $i\af_M^*$ in the 
direction of $\Lambda$. 
Then it follows that there exists $C>0$ such that
\begin{equation}\label{6.24}
\begin{split}
\parallel\mN_S^\prime(P,\pi,\lambda)_{K_f,\sigma}\parallel\le&
C\Biggl(\sum_{v\in S_0\setminus\{\infty\}}\sum_{k=1}^q\parallel D_\Lambda^k 
R_{Q|P}(\pi_v,\lambda)_{K_v}\parallel\\
&\quad\quad+\sum_{k=1}^q\parallel D_\Lambda^k 
R_{Q|P}(\pi_\infty,\lambda)_{\sigma}\parallel
\Biggr)
\end{split}
\end{equation}
for all $\lambda\in i\af_M^*$, $\sigma\in\Pi(K_\infty)$ and $\pi\in\Pi(M(\A))$.
Together with Proposition \ref{p6.3} this implies Theorem \ref{0.1}.

\hfill$\Box$

\smallskip
\noindent
{\it Proof of Theorem \ref{th0.2}:}
The proof of Theorem \ref{th0.2} is similar to the proof of Theorem \ref{0.1}.
We only have to modify some of the arguments.
Given an open compact subgroup $K_f$ of $G(\A_f)$ and $\sigma\in\Pi(K_\infty)$, let $\Pi_{K_f,\sigma}$ denote the orthogonal projection
of the Hilbert space $\ov\cA^2_{\chi,\pi}(P)$ onto the finite dimensional
subspace $\cA^2_{\chi,\pi}(P)_{K_f,\sigma}$. Let $h\in\Co^1(G(\A)^1)$ be $K$-finite.
Then there exists an open compact subgroup $K_f$ of $G(\A_f)$ such that $h$ is
left and right invariant under $K_f$. Furthermore, there exist 
$\sigma_1,...,\sigma_m\in\Pi(K_\infty)$ such that
\begin{equation}\label{6.25}
\rho_{\chi,\pi}(P,\lambda,h)=\sum_{i,j=1}^m\Pi_{K_f,\sigma_i}\circ
\rho_{\chi,\pi}(P,\lambda,h)\circ\Pi_{K_f,\sigma_j}
\end{equation}
for all $\pi\in\Pi(M(\A)^1)$ and $\chi\in\mX$. 
Let $k\in\N$. Then by (\ref{6.25}) we
get
\begin{equation}\label{6.26}
\begin{split}
\parallel \mN_S^\prime&(P,\pi,\lambda)\rho_{\chi,\pi}(P,\lambda,h)\parallel_1\\
&\le\sum_{i=1}^m\parallel \mN_S^\prime(P,\pi,\lambda)_{K_f,\sigma_i}
\parallel\cdot
\parallel\rho_{\chi,\pi}(P,\lambda,\Delta^{2k})^{-1}_{K_f,\sigma_i}
\parallel_1\\
&\hskip5truecm\cdot\parallel\rho_{\chi,\pi}(P,\lambda,\Delta^{2k}h)\parallel.
\end{split}
\end{equation}
Here we assume, of course, that $\cA^2_{\chi,\pi}(P)_{K_f,\sigma_i}\not=0$, 
$i=1,...,m$. Then it follows from (6.9)  that
\begin{equation}\label{6.27}
\parallel \rho_{\chi,\pi}(P,\lambda,\Delta^{2k})^{-1}_{K_f,\sigma_i} 
\parallel_1
\le C\frac{\dim\cA^2_{\chi,\pi}(P)_{K_f,\sigma_i}}
{(1+\parallel\lambda\parallel^2+\lambda_\pi^2)^k}
\end{equation}

for $i=1,...,m$. Given $\sigma\in\Pi(K_\infty)$, let 
\begin{equation*}
\begin{split}
\Pi_{\di}&(M(\A)^1)_{K_f,\sigma}\\
&=\bigl\{\pi\in\Pi_{\di}(M(\A)^1;K_f)\;|\;
[I_P^G(\pi_\infty)|_{K_\infty}:\sigma]>0\bigr\}.
\end{split}
\end{equation*}

Then we proceed as above to show that for every $n\in\N$ 
there exists $C_n>0$ such that
\begin{equation}\label{6.28}
\begin{split}
\sum_{\chi\in\mX}\sum_{\pi\in\Pi(M(\A)^1)}&
\int_{i\af_L^*/i\af_G^*}\parallel \mN_S^\prime(P,\pi,\lambda)
\rho_{\chi,\pi}(P,\lambda,h)\parallel_1|\nu_L^S(P,\pi,\lambda)|\;d\lambda\\
&\le C_n \sum_{i=1}^m\sum_{\pi\in\Pi_{\di}(M(\A)^1)_{K_f,\sigma_i}}m(\pi)
(1+\lambda_\pi^2)^{-n}.
\end{split}
\end{equation}
To estimate the right hand side, we fix $\sigma\in \Pi(K_\infty)$. Then
as in (\ref{6.20}) we get
\begin{equation*}
\begin{split}
&\sum_{\pi\in\Pi_{\di}(M(\A)^1)_{K_f,\sigma}}\frac{m(\pi)}{(1+\lambda_\pi^2)^{n}}\\
&\hskip8pt\le\sum_{i=1}^l\sum_{\pi_\infty\in\Pi(M(\R))}
m_{\Gamma_{M,i}}(\pi_\infty)\frac{\dim(\H(\pi_\infty)\otimes V_\sigma)^{K_{M,\infty}}}{(1+\lambda_{\pi_\infty}^2)^{n}}.
\end{split}
\end{equation*}
It follows from Theorem 0.1 of \cite{Mu1} that for sufficiently large $n$,
this series is convergent. This completes the proof of Theorem \ref{th0.2}.

\hfill$\Box$

\smallskip
We observe that for tempered representations, the existence of estimates
like (\ref{0.2}), (\ref{0.3}) and (\ref{0.4})
follows from results of Arthur  \cite[p.51]{A5} and 
\cite[Lemma 2.1]{A8}. Let $\Pi_{\tm}(M(\A)^1)$ be the subspace of all $\pi$ in
$\Pi(M(\A)^1)$ such that the local constituents  $\pi_v$ of $\pi$ are
tempered for all $v$. Then  we  obtain 
\begin{prop}\label{p6.4}
For every $M\in{\L}(M_0)$, $L\in{\L}(M)$ and $P\in{\cP}(M)$ we have
\begin{equation*}
\sum_{\chi\in\mX}\sum_{\pi\in\Pi_{\tm}(M(\A)^1)}\int_{i\af_L^*/i\af_G^*} 
\parallel
\mM_L(P,\pi,\lambda)\rho_{\chi,\pi}(P,\lambda,h)\parallel_1d\lambda<\infty.
\end{equation*}
\end{prop}

\section{The example of $\GL_n$}
\setcounter{equation}{0}

In this section we shall briefly discuss the case where $G=\GL_n$. Let $P_0$
be the subgroup of upper triangular matrices of $G$. This is the minimal 
standard parabolic subgroup of $G$. Its Levi subgroup $M_0$ is the group of 
diagonal matrices. Let $P$ be a parabolic subgroup of $G$ defined over $\Q$, 
and let $M$ be the unique Levi component of $P$ which contains $M_0$. Then 
$$M\cong\GL_{n_1}\times\cdots\times\GL_{n_r}.$$
We shall identify $\af_M$ with $\R^r$. Let $e_1,...,e_r$ denote the standard 
basis of $(\R^r)^*$. Then the roots $\Sigma_P$ are given by
$$\Sigma_P=\{e_i-e_j\mid 1\le i<j\le r\}.$$
Let $v$ be a place of $\Q$. Fix a nontrivial continuous character $\psi_v$ of
the additive group $\Q_v^+$ of $\Q_v$ and equip $\Q_v$ with the Haar measure
which is selfdual with respect to $\psi_v$. Given irreducible unitary
representations $\pi_{1,v}$ and $\pi_{2,v}$ of $\GL_{n_1}(\Q_v)$ and 
$\GL_{n_2}(\Q_v)$,
 respectively, let $L(s,\pi_{1,v}\times\pi_{2,v})$ and 
$\epsilon(s,\pi_{1,v}\times\pi_{2,v},\psi_v)$ denote the Rankin-Selberg 
$L$-factor and the $\epsilon$-factor 
defined by Jacquet, Piateski-Shapiro, and Shalika \cite{JPS}, \cite{JS1}.

Let $P_1,P_2\in\cP(M)$. Then there exist permutations 
$\sigma_1,\sigma_2\in S_r$ such that the set of roots of $(P_i,A_i)$ is
given by
$$\Sigma_{P_k}=\{e_i-e_j\mid \sigma_k(i)<\sigma_k(j)\}.$$
Put
$$I(\sigma_1,\sigma_2)=\{(i,j)\mid 1\le i,j\le r,\; \sigma_1(i)<\sigma_1(j),
\sigma_2(i)>\sigma_2(j)\}.$$
Then
$$\Sigma_{P_1}\cap\Sigma_{\ov P_2}=\{e_i-e_j\mid (i,j)\in I(\sigma_1,\sigma_2)
\}.$$
Let $\pi_v=\pi_{1,v}\otimes\cdots\otimes\pi_{r,v}$, where $\pi_{i,v}\in\Pi(\GL_{n_i}(\Q_v))$,
 $i=1,...,r$. Given $\bs=(s_1,...,s_r)\in\C^r$, set
\begin{equation*}
\begin{split}
r_{P_2|P_1}&(\pi_v,\bs):=\\
&\prod_{(i,j)\in I(\sigma_1,\sigma_2)}
\frac{L(s_i-s_j,\pi_{i,v}\times\widetilde\pi_{j,v})}{L(1+s_i-s_j,\pi_{i,v}\times\widetilde\pi_{j,v})\epsilon(s_i-s_j,\pi_{i,v}\times\widetilde\pi_{j,v},\psi_v)}.
\end{split}
\end{equation*}
As explained in \cite[p.87]{AC}, the meromorphic functions 
$r_{P_2|P_1}(\pi,\bs)$ satisfy all the properties of Theorem 2.1 of \cite{A7}
and  they are the natural choice of normalizing factors in the case of $\GL_n$.
 We note that they do not coincide with the normalizing factors used 
in the previous sections. They differ, however, only by a factor which can be
expressed in terms of the $\epsilon$-factors. 

Now let $\pi_1$ and $\pi_2$ be automorphic representations of $\GL_{n_1}(\A)$
and $\GL_{n_2}(\A)$, respectively. Then the global Rankin-Selberg $L$-function
$L(s,\pi_1\times\pi_2)$ is defined by
$$L(s,\pi_1\times\pi_2)=\prod_vL(s,\pi_{1,v}\times\pi_{2,v}),$$
where the product is over all places $v$ of $\Q$ and $\pi_i=
\otimes_i\pi_{i,v}$. The product converges absolutely in a half-plane 
$\Re(s)\gg0$. If $\pi_1$ and $\pi_2$ belong to the discrete spectrum of 
$\GL_{n_1}(\A)$ and $\GL_{n_2}(\A)$, respectively, then $L(s,\pi_1\times\pi_2)$
admits a meromorphic extension to the whole complex plane. 

To define the 
global $\epsilon$-factor $\epsilon(s,\pi_1\times\pi_2)$ one has to pick a
non-trivial continuous character $\psi:\A^+\to\C^\times$ of the additive
 group $\A^+$ of $\A$. Then $\psi=\otimes_v\psi_v$ and $\epsilon(s,\pi_{1,v}\times\pi_{2,v},\psi_v)=1$ for almost all places $v$. Hence the product
$$\epsilon(s,\pi_1\times\pi_2,\psi)=\prod_v\epsilon(s,\pi_{1,v}\times\pi_{2,v},\psi_v)$$
exists for all $s\in\C$ and defines an entire function. The global 
$\epsilon$-factor is  independent of $\psi$ and therefore, will be denoted by
$\epsilon(s,\pi_1\times\pi_2)$. 

Let $\pi\in\Pi_{\di}(M(\A))$. Then $\pi=\pi_1\otimes\cdots\otimes\pi_r$ with
$\pi_i\in\Pi_{\di}(\GL_{n_i}(\A))$ and for $\bs\in\C^r$, the global
 normalizing factor is defined by
$$r_{P_2|P_1}(\pi,\bs):=\prod_{(i,j)\in I(\sigma_1,\sigma_2)}
\frac{L(s_i-s_j,\pi_{i}\times\widetilde\pi_{j})}{L(1+s_i-s_j,\pi_{i}\times\widetilde\pi_{j})\epsilon(s_i-s_j,\pi_{i}\times\widetilde\pi_{j})}.$$
Theorem \ref{th5.3} is closely related to the estimation of the winding
numbers
$$\int_1^\lambda\frac{L'(1+it,\pi_1\times\widetilde\pi_2)}{L(1+it,\pi_1\times\widetilde\pi_2)}\;dt$$
with upper bounds depending on the Casimir eigenvalues of $\pi_{1,\infty}$ and
$\pi_{2,\infty}$
in the same way as in Theorem \ref{th5.3}. In the present case, such
estimates can be obtained using standard methods of analytic number theory.
In fact, the bounds can be improved considerably.

Next we discuss the conditions (\ref{0.2}) and (\ref{0.4}). As mentioned in
the introduction, for $\GL_n$ it is possible to prove that (\ref{0.2}) and
(\ref{0.4}) hold. We shall briefly indicate the main steps of the proof.
Let $\rho$ be a representation
of $\GL_m(\Q_v)$ and $s\in\C$. Then we denote by $\rho[s]$ the representation
 of
$\GL_m(\Q_v)$ defined by
$$\rho[s](g)=|\det g|^s\rho(g),\quad g\in\GL_m(\Q_v).$$
Let $\pi$ be a cuspidal automorphic representation of $\GL_m(\A)$. Then it is
 known that each local component $\pi_v$ of $\pi$ is generic \cite{Sk} and
therefore, by \cite{JS2}  it follows that $\pi_v$ is
 equivalent to a full induced representation, i.e.,
\begin{equation}\label{7.1}
\pi_v\cong I^{\GL_m}_P(\tau_1[t_1]\otimes\cdots\otimes\tau_r[t_r]),
\end{equation}
where $P$ is a standard parabolic subgroup of $\GL_m$ with Levi component
$\GL_{m_1}\times\cdots\times\GL_{m_r}$, $\tau_i$ is a tempered
 representation of
$\GL_{m_i}(\Q_v)$ and the $t_i$'s are real numbers satisfying
$$t_1> t_2>\cdots >t_r,\quad |t_i|<1/2,\; i=1,...,r.$$
For $\pi_v$ unramified, Luo, Rudnick and Sarnak \cite{LRS} proved that the 
parameters $t_i$ satisfy the following nontrivial bound:
\begin{equation}\label{7.2}
\max_i|t_i|<\frac{1}{2}-\frac{1}{m^2+1}.
\end{equation}
Using the same method, one can show that (\ref{7.2})  holds at all places. 
Now let $P$ be a standard parabolic subgroup of $\GL_n$ with Levi component 
$\GL_{n_1}\times\cdots\times\GL_{n_r}$ and let $\pi_v$ be the local
 $v$-component
of a cuspidal automorphic representation of $M(\A)$. Then
 $\pi_v=\otimes_i\pi_{i,v}$
and each $\pi_{i,v}$ is a full induced representation of the form (\ref{7.1})
 with parameters $t_{ij}$ satisfying (\ref{7.2}). Using induction in
 stages, it follows that for each $i$  there exist a parabolic subgroup 
$R_i$ of $\GL_{n_i}(\Q_v)$ of type $(n_{i1},...,n_{il_i})$, a discrete series
 representation $\delta_{i,v}$ of $M_{R_i}(\Q_v)$  and  $\bt_i
=(t_{i1},...,t_{il_i})\in\R^{l_i}$ 
satisfying 
\begin{equation}\label{7.3}
\max_i|t_{ij}|<\frac{1}{2}-\frac{1}{n^2_i+1},
\end{equation}

such that $\pi_{i,v}\cong I^M_{R_i}(\delta_{i,v},\bt_i)$. Put 
$l=l_1+\cdots +l_r$, 
$$\delta_v=\otimes_i\delta_{i,v},\quad \bt=(t_{11},...,t_{1l_1},...,t_{r1},
...,t_{rl_r}).$$

Generalizing property (R.2) of 
\cite[p.172]{A8}, we get
\begin{equation}\label{7.4}
R_{Q|P}(\pi_v,\bs)=R_{Q(R)|P(R)}(\delta_v,\bs+\bt),\quad \bs\in\C^l,
\end{equation}
where $\bs$ is identified with an element in $\C^l$ with respect to the
 embedding
which corresponds to the canonical embedding $\af_M^*\subset\af_{P(R)}^*.$ 
This leads to an immediate reduction of the problem. We can assume that $\pi_v$
is square integrable. However, now we have to estimate the norm  of the
derivatives of $R_{Q|P}(\pi_v,\bs)_{K_v}$  
(resp. $R_{Q|P}(\pi_v,\bs)_{\sigma_v}$)  in the domain
$$\{\bs\in\C^r\mid |\Re(s_i)|<1/2-1/(n^2+1),\;i=1,...,r\}.$$
The important point is that for $\pi_v$ square integrable, 
$R_{Q|P}(\pi_v,\bs)$ is holomorphic in the domain
$$\{\bs\in\C^r\mid \Re(s_i-s_j)>-1,\; 1\le i<j\le r\}$$ 
\cite{MW}.
Using the product formula for normalized intertwining operators, the above
problem
can be further reduced to the case where $P$ is maximal parabolic and
$Q=\ov P$. Then $M=\GL_{n_1}\times\GL_{n_2},$
$\pi_v=\pi_{1,v}\otimes\pi_{2,v}$, and we may regard the intertwining 
operator as a function $R_{\ov P|P}(\pi_v,s)$ of one complex variable.
Now we distinguish two cases.

{\bf 1.} $v<\infty$. 

Let $K_v\subset \GL_n(\Q_v)$ be an open compact subgroup. We may assume that
$K_v$ is a congruence subgroup. Then we have to estimate the norm of 
derivatives
of $R_{\ov P|P}(\pi_v,s)_{K_v}$ in the strip $|\Re(s)|<1-2/(n^2+1)$. Let 
$$K_{M,v}=K_v\cap M(\Q_v).$$
Then $K_{M,v}$ is an open compact subgroup of $M(\Q_v)$. Let ${\bf 1}$ denote
the trivial representation of $K_{M,v}$ and let $\Pi_2(M(\Q_v);K_{M,v})$ be 
the set of all $\pi_v\in\Pi_2(M(\Q_v))$ such that 
$[\pi_v|_{K_{M,v}}:{\bf 1}]>0$.
By Theorem 10 of \cite{HC2}, $\Pi_2(M(\Q_v);K_{M,v})$ is a compact subset
of $\Pi_2(M(\Q_v))$. Furthermore, $\af_M^*\cong\R^2$ acts on 
$\Pi_2(M(\Q_v))$ by
$$\pi_{1,v}\otimes\pi_{2,v}\mapsto\pi_{1,v}[iu_1]\otimes\pi_{2,v}[iu_2],\quad
(u_1,u_2)\in\R^2.$$
The stabilizer  of a given representation $\pi_v$ is a lattice $L\subset\R^2$
so that the orbit $\ho_{\pi_v}$ of $\pi_v$ is a compact torus $\R^2/L$. Thus
there exist $\delta_1,...,\delta_l\in\Pi_2(M(\Q_v);K_{M,v})$ such that
$$\Pi_2(M(\Q_v);K_{M,v})=\ho_{\delta_1}\sqcup\cdots\sqcup\ho_{\delta_l}.$$
Since 
$$R_{\ov P|P}(\pi_{1,v}[iu_1]\otimes\pi_{2,v}[iu_2],s)=R_{\ov P|P}(\pi_{1,v}\otimes\pi_{2,v},s+i(u_1+u_2)),$$
it suffices to consider a fixed discrete series representation $\pi_v$. Now
recall that $R_{\ov P|P}(\pi_v,s)$ is holomorphic in the strip $|\Re(s)|<1$.
Furthermore by Theorem 2.1 of \cite{A7}, $R_{\ov P|P}(\pi_v,s)_{K_v}$ is a
finite rank matrix whose entries are rational functions of $p^{-s}_v$. Hence
for every $u\in\R$, $R_{\ov P|P}(\pi_v,u+iw)_{K_v}$ is a periodic function
of $w\in\R$. From these
observations it follows immediately that for every $k\in\N_0$ there exists
$C>0$ such that
\begin{equation}\label{7.5}
\parallel D^k_s R_{\ov P|P}(\pi_v,s)_{K_v}\parallel\le C
\end{equation}
for all $s\in\C$ in the strip $|\Re(s)|\le 1-2/(n^2+1).$

{\bf 2.} $v=\infty$.

Let $\sigma_v\in\Pi(\rO(n))$. Then we have to estimate the norm of derivatives of
$R_{\ov P|P}(\pi_v,s)_{\sigma_v}$ in the strip $|\Re(s)|<1-2/(n^2+1)$. First
note that
$$M(\R)\cong(\R^*)^2\times(\SL_{n_1}(\R)\times\SL_{n_2}(\R)).$$
Furthermore the set of discrete series representations of $\SL_{n_i}(\R)$
containing a fixed ${\SO}(n_i)$-type is finite \cite[p.398]{Wa2}. Hence 
in the same way as above, it follows that we can fix the discrete series
representation $\pi_v$. Again $R_{\ov P|P}(\pi_v,s)$ is holomorphic in the
strip $|\Re(s)|<1$ and by Theorem 2.1 of \cite{A7}, 
$R_{\ov P|P}(\pi_v,s)_{\sigma_v}$ is a rational function of $s\in\C$. This
implies that for every $k\in\N_0$ there exist $C>0$ and $N\in\N$ such that
\begin{equation}\label{7.6}
\parallel D^k_s R_{\ov P|P}(\pi_v,s)_{\sigma_v}\parallel\le C(1+|s|)^N
\end{equation}
for all $s\in\C$ with $|\Re(s)|<1-2/(n^2+1)$.

Combining (\ref{7.5}) and (\ref{7.6}) with the various 
steps of the reduction it follows that (\ref{0.2}) and (\ref{0.4})
hold for all local components $\pi_v$ of cuspidal automorphic representations.

It remains to deal with local components of automorphic forms in the
 residual spectrum. For this purpose we use the description of the residual
spectrum given by M{\oe}glin and Waldspurger \cite{MW}. First we recall the
notion of a Speh representation \cite[I.5]{MW}. Let $k|m$, $d=m/k$ and $R$
a standard parabolic subgroup of $\GL_m$ of type $(d,...,d)$. Let $\delta$ be
a discrete series representation of $\GL_d(\Q_v)$. Then the induced 
representation
$$I^{\GL_m}_R(\delta[(k-1)/2]\otimes\delta[(k-3)/2\otimes\cdots\otimes\delta
[(1-k)/2])$$
has a unique irreducible quotient which we denote by $J(\delta,k)$. It follows
from Theorem D of \cite{Ta} and \cite{Vo} that for every 
$\pi_v\in\Pi(\GL_m(\Q_v))$ there exist a standard parabolic subgroup $P$ of
type $(m_1,...,m_r)$, $k_i|m_i$, discrete series representations $\delta_i$
of $\GL_{d_i}(\Q_v)$, $d_i=m_i/k_i$, and real numbers $t_1,...,t_r$ 
satisfying $|t_i|<1/2$ such that
\begin{equation*}
\pi_v\cong I^{\GL_m}_P(J(\delta_1,k_1)[t_1]\otimes\cdots\otimes J(\delta_r,k_r)[t_r]).
\end{equation*}

Now suppose that $\pi_v$ is a local component of an automorphic 
representation $\pi$ in the residual spectrum of $\GL_m(\A)$. By \cite{MW}
there exist a standard parabolic subgroup $Q$ of $\GL_m$ of type $(d,...,d)$
and a cuspidal automorphic representation $\mu$ of $\GL_d(\A)$ such that
$\pi_v$ is the unique irreducible quotient of the induced representation
$$I^{\GL_m}_Q(\mu_v[(k-1)/2]\otimes\mu_v[(k-3)/2]\otimes\cdots\otimes
\mu_v[(1-k)/2]),$$
where $\mu_v$ is the $v$-component of $\mu$. 
As explained above, 
$\mu_v$ is equivalent to an induced representation of the form (\ref{7.1})
with parameters $t_i$ satisfying (\ref{7.2}). Using induction in stages,
it follows that
$$\mu_v\cong I^{\GL_d}_R(\delta_1[t_1]\otimes\cdots\otimes\delta_r[t_r]),$$
where $R$ is a standard parabolic subgroup of $\GL_d$ of type $(d_1,...,d_r)$,
$\delta_i$ is a discrete series representations of $\GL_{d_i}(\Q_v)$, 
$i=1,...,r$, and the parameters $t_i$ satisfy 
$t_1\ge t_2\ge\cdots\ge t_r$ and (\ref{7.2}). Then
it follows from Proposition I.9 and Lemma I.8 of \cite{MW} that there is a 
standard parabolic subgroup $P$ of $\GL_m$ of type $(kd_1,...,kd_r)$ such that
\begin{equation}\label{7.7}
\pi_v\cong I^{\GL_m}_P(J(\delta_1,k)[t_1]\otimes\cdots\otimes J(\delta_r,k)[t_r])
\end{equation}
and 
\begin{equation}\label{7.8}
\max_i|t_i|<\frac{1}{2}-\frac{1}{m^2+1}.
\end{equation}
This is the extension of the results of \cite{LRS} to local components of
automorphic representations in the discrete spectrum.

Now we can  proceed in  the same way as in the 
cuspidal case. The only  difference is that we have to deal  with
the slightly more general Speh representations in place of the discrete series
representations. In this way one can establish 
(\ref{0.2}) and (\ref{0.4}). This implies that for $\GL_n$ the spectral side
of the Arthur trace formula is absolutely convergent. Details will appear in
\cite{MS}.

\end{document}